\documentclass[12pt]{amsart}

\usepackage{amscd,latexsym,amsthm,amsfonts,amssymb,amsmath,amsxtra}
\usepackage[mathscr]{eucal}
\renewcommand\theequation{\thesection.\arabic{equation}}

\newcommand{\BA}{{\mathbb {A}}}

\newcommand{\BC}{{\mathbb {C}}}
\newcommand{\BD}{{\mathbb {D}}}

\newcommand{\BR}{{\mathbb {R}}}

\newcommand{\BZ}{{\mathbb {Z}}}

\newcommand{\CA}{{\mathcal {A}}}

\newcommand{\CE}{{\mathcal {E}}}
\newcommand{\CF}{{\mathcal {F}}}

\newcommand{\CI}{{\mathcal {I}}}

\newcommand{\CK}{{\mathcal {K}}}

\newcommand{\CN}{{\mathcal {N}}}
\newcommand{\CO}{{\mathcal {O}}}

\newcommand{\CS}{{\mathcal {S}}}
\newcommand{\CT}{{\mathcal {T}}}

\newcommand{\RR}{{\mathrm {R}}}

\newcommand{\RU}{{\mathrm {U}}}

\newcommand{\Asai}{{\mathrm{Asai}}}

\newcommand{\cusp}{{\mathrm{cusp}}}

\newcommand{\disc}{{\mathrm{disc}}}

\newcommand{\el}{{\mathrm{ell}}}

\newcommand{\fj}{{\mathcal{FJ}}}

\newcommand{\Gal}{{\mathrm{Gal}}}
\newcommand{\GL}{{\mathrm{GL}}}

\newcommand{\GSp}{{\mathrm{GSp}}}

\newcommand{\Mat}{{\mathrm{Mat}}}

\renewcommand{\Re}{{\mathrm{Re}}}

\newcommand{\simp}{{\mathrm{sim}}}
\newcommand{\SL}{{\mathrm{SL}}}

\newcommand{\SO}{{\mathrm{SO}}}

\newcommand{\Sp}{{\mathrm{Sp}}}
\newcommand{\Stab}{{\mathrm{Stab}}}
\newcommand{\st}{{\mathrm{st}}}

\newcommand{\tr}{{\mathrm{tr}}}

\newcommand{\ovl}{\overline}
\newcommand{\udl}{\underline}
\newcommand{\wt}{\widetilde}

\def\bks{{\backslash}}

\def\diag{{\rm diag}}

\def\lam{{\lambda}}

\def\sym{{\rm sym}}
\def\sig{{\sigma}}

\def\tilpi{{\widetilde{\pi}}}

\def\tilsp{{\widetilde{\Sp}}}

\newtheorem{thm}{Theorem}[section]

\newtheorem {conj}[thm]{Conjecture}

\newtheorem {ques/conj}[thm]{Question/Conjecture}

\newtheorem{rmk}[thm]{Remark}

\makeatletter

\newcommand{\Rmnum}[1]{\expandafter\@slowromancap\romannumeral #1@}
\makeatother

\begin{document}
\renewcommand{\theequation}{\arabic{equation}}
\numberwithin{equation}{section}

\title[Automorphic Integral Transforms for Classical Groups I]{Automorphic Integral Transforms for Classical Groups I: Endoscopy Correspondences}

\author{Dihua Jiang}
\address{School of Mathematics\\
University of Minnesota\\
Minneapolis, MN 55455, USA}
\email{dhjiang@math.umn.edu}


\date{December, 2012}


\thanks{}


\begin{abstract}
A general framework of constructions of endoscopy correspondences via automorphic integral transforms for classical groups is formulated in terms of
the Arthur classification of the discrete spectrum of square-integrable automorphic forms. This suggests a principle, which is called
the $(\tau,b)$-theory of automorphic forms of classical groups, to reorganize and extend the series of work of Piatetski-Shapiro, Rallis,
Kudla and others on standard $L$-functions of classical groups and theta correspondence.
\end{abstract}

\maketitle


\section{Introduction}

Automorphic forms are transcendental functions with abundant symmetries and fundamental objects to arithmetic and geometry.
In the theory of automorphic forms, it is an important and difficult problem to construct explicitly automorphic forms with specified properties,
in particular, to construct cuspidal automorphic forms of various types.

Ilya Piatetski-Shapiro is a pioneer to use representation theory to study automorphic forms. The main idea is to construct
certain models of representation-theoretic nature to obtain cuspidal automorphic forms, that is, to explicitly realize cuspidal automorphic forms
in the space of functions over certain geometric spaces.
The celebrated theorem of Gelfand and Piatetski-Shapiro shows that all cuspidal automorphic forms are rapidly
decreasing on a fundamental domain $\BD$ when approaching the cusps, and hence can be realized discretely in the space $L^2(\BD)$
of all square integrable functions on $\BD$. Therefore, it is fundamental to understand the space $L^2(\BD)$ or more precisely
the discrete spectrum of $L^2(\BD)$.

In the modern theory of automorphic forms, this problem is formulated as follows. Let $G$ be a reductive algebraic
group defined over a number field (or a global field) $F$ and consider the quotient space
$$
X_G=Z_G(\BA)G(F)\bks G(\BA)
$$
where $\BA$ is the ring of adeles of $F$ and $Z_G$ is the center of $G$. As a consequence of the reduction theory in
the arithmetic theory of algebraic groups, this space $X_G$ has finite volume with respect to the canonical
quotient measure. Geometrically $X_G$ is a modern replacement of the classical locally symmetric spaces which are
the land for classical automorphic forms. It is fundamental to understand the space $L^2(X_G)$ of all
square integrable functions. It is clear that under the right translation by $G(\BA)$, the space $L^2(X_G)$ is
unitary representation of $G(\BA)$ with respect to the $L^2$-inner product. The theorem of Gelfand and Piatetski-Shapiro
says that all irreducible cuspidal automorphic representations of $G(\BA)$ (with trivial central character) can
be realized as irreducible submodules of $L^2(X_G)$. Combining with the Langlands theory of Eisenstein series, one
can state the following fundamental result in the theory of automorphic forms.

\begin{thm}[Gelfand, Piatetski-Shapiro; Langlands]\label{gpsl}
The discrete spectrum of the space of square integrable automorphic functions has the following decomposition
$$
L^2_\disc(X_G)
=
\oplus_{\pi\in \widehat{G(\BA)}}m_\disc(\pi)V_\pi
$$
with finite multiplicity $m_\disc(\pi)$ for all irreducible unitary representations $(\pi,V_\pi)$ of $G(\BA)$,
where $\widehat{G(\BA)}$ is the unitary dual of $G(\BA)$.
\end{thm}

One of the remaining fundamental problems concerning $L^2_\disc(X_G)$ is to understand the multiplicity $m_\disc(\pi)$
in general.

When $G=\GL_n$, the celebrated theorems of Shalika (\cite{Sl74})
and of Piatetski-Shapiro (\cite{PS79}) for cuspidal spectrum; and of Moeglin and Waldspurger (\cite{MW89})
for residual spectrum assert that for any irreducible unitary representation $\pi$ of $\GL_n(\BA)$
has multiplicity $m_\disc(\pi)$ at most one. However, when $G\neq\GL_n$, the situation could be very different.

For instance, when $G=\SL_n$ with $n\geq 3$, D. Blasius finds a family of cuspidal automorphic representations
with higher multiplicity, i.e. $m_\disc(\pi)>1$ (\cite{Bl94}). See also the work of E. Lapid (\cite{Lp99}). However,
when $n=2$ the multiplicity is one, which is proved by D. Ramakrishnan (\cite{Rm00}) based on an earlier work of Labesse and Langlands (\cite{LL79}).
More interestingly, when $G=G_2$, the $F$-split exceptional group of type $G_2$,
Gan, Gurevich and the author show that there exist a family of cuspidal automorphic representations of $G_2(\BA)$,
whose cuspidal multiplicity can be as high as possible (\cite{GGJ02}), and see \cite{Gn05} for more complete result.
From the argument of \cite{GGJ02}, such a result can also be expected for other exceptional groups. On the other hand, it does not happen for $F$-quasisplit
classical groups following Arthur's theorem in his forthcoming book (\cite{Ar12}), but when classical groups are not $F$-quasisplit,
the multiplicity could still be higher according to the work of J.-S. Li (\cite{Li95} and \cite{Li97}).
Now, following the Arthur multiplicity formula (\cite{Ar12} and \cite{Mk12}), the multiplicity in the discrete
spectrum of $F$-quasisplit classical groups is at most one in general, except for certain cases of even special orthogonal
groups where the discrete multiplicity could be two.

When $G=\GL_n$, the discrete spectrum was completely determined by Moeglin and Waldspurger in 1989 (\cite{MW89}),
which confirms a conjecture of Jacquet in 1983. For $G$ to be a $F$-quasisplit classical groups,
the stable trace formula of Arthur is able to give a complete description of
the discrete spectrum, up to global Arthur packets, in terms of a certain family of automorphic representations of $\GL$, based on the
fundamental work of B.-C. Ngo (\cite{N10}) and of G. Laumon, P.-H. Chaudouard, J.-L. Waldspurger,  and others
(see the references within \cite{Ar12} and \cite{Mk12}).
If $G$ is not $F$-quasisplit (i.e. general inner forms), much more work needs to be done (\cite[Chapter 9]{Ar12}).
It is far from reach if $G$ is of exceptional type.

One of the great outcomes from the trace formula approach to understand the structure of the discrete spectrum of classical groups is the
discovery of the endoscopy structure of automorphic forms by R. Langlands (\cite{L79} and \cite{L83}), the twisted version of which is referred to
the book of R. Kottwitz and D. Shelstad (\cite{KtS99}). The trace formula method confirms the existence of the endoscopy
structure of automorphic forms on classical groups, which leads a uniform proof of the existence of the weak functorial transfer from
$F$-quasisplit classical groups to the general linear groups for automorphic representations occurring in the discrete spectrum of classical
groups, among other fundamentally important results (\cite{Ar12}, and also \cite{Mk12}), based on the fundamental work of
Laumon, Ngo, Waldspurger, Chaudouard, and others.

Endoscopy transfers are the Langlands functorial transfers, which indicates that automorphic representations $G(\BA)$ may come from
its endoscopy groups. This is a special case of the general Langlands functoriality conjecture (\cite{L70}). In addition to the trace formula approach,
there have been other methods to prove the existence or to construct explicitly various cases of the Langlands functorial transfers.

One of the major methods is to establish the Langlands functorial transfers from classical groups to general linear groups via the Converse Theorem
of Cogdell and Piatetski-Shapiro (\cite{CPS94}) combined with the analytic properties of tensor product $L$-functions of
classical groups and general linear groups. This method has been successful for generic cuspidal automorphic representations of
$F$-quasisplit classical groups (\cite{CKPSS04}, \cite{CPSS11}, \cite{KK05}, and also see \cite{AS06} for spinor group case)
and has been extended to prove the existence of the $n$-symmetric power liftings of $\GL(2)$
for $n=2,3,4$ (\cite{GlJ78}, \cite{KmS02} and \cite{Km03}). As outlined by D. Soudry in \cite{Sd06},
this method can be essentially used to establish the
Langlands functorial transfers (the strong version, i.e. compatible with the local Langlands transfer at all local places)
from classical groups to general linear groups for all cuspidal automorphic representations, as long as
the analytic properties of tensor product $L$-functions of classical groups and general linear groups can be established in this
generality. Some recent progress in this aspect is referred to the work of Pitale, Saha, and Schmidt (\cite{PSS11}) and of
Lei Zhang and the author (\cite{JZ12} and also \cite{Sn12} and \cite{Sn13}).

Another method to establish certain cases of the existence of the Langlands functoriality is to use the relative trace formula approach
which was introduced by H. Jacquet. See \cite{Jc97} for discussions in this aspect and \cite{Jc04} and \cite{MR10} for
more recent progress.

For the Langlands functorial transfers beyond the endoscopy theory, it is referred to the work of Langlands
(\cite{L04}) and the more recent work of Frenkel, Langlands and Ngo (\cite{FLN10} and \cite{FN11}). An alternative
approach is taken by L. Lafforgue (\cite{Lf10} and \cite{Lf12}, and also \cite{Jc12}).

The remaining methods are more constructive in nature.

First, the {\bf theta correspondence} for reductive dual pairs in the sense of Howe is a method
to construct automorphic forms using the theta functions built from the Weil representation on the Schr\"odinger model.
The method of using classical theta functions to construct automorphic forms goes back at least to
the work of C. Siegel (\cite{Sg66}) and
the work of A. Weil (\cite{Wl64} and \cite{Wl65}). The fundamental work of Howe (\cite{Hw79})
started the representation-theoretic approach in
this method guided by the idea from Invariant Theory.

This idea (and the method) has been extended a similar constructions for exceptional groups and reductive dual pairs in the 1990's, and is usually
called the {\bf exceptional theta correspondence}. The main work includes that of Kazhdan and Savin (\cite{KzS90}), that of Ginzburg,
Rallis, and Soudry (\cite{GRS97}), and many others.

More recently, a striking method to construct more general automorphic forms was carried by Ginzburg, Rallis, and Soudry in 1999
(called the {\bf automorphic descent method}), which extends the classical examples
of Eichler and Zagier (\cite{EZ85}) and of I. Piatetski-Shapiro on $\Sp(4)$ (\cite{PS83})
to general $F$-quasisplit classical groups and establishes the relation
between their constructions and the Langlands functoriality between classical groups and general linear groups.
This construction is discussed in general in their
recent book (\cite{GRS11}) and will be reviewed with some details in Section 6, including some extension of the methods through the
work of Ginzburg, Soudry, Baiying Liu, Lei Zhang and the author (\cite{GJS12}, \cite{GJS13}, \cite{JL13}, \cite{JZ13} and \cite{Liu13}).

The objective here is to reformulate and extend the previous known cases of constructions and the related theory in a much more general
framework in terms of the Arthur
classification of the discrete spectrum of classical groups, and suggest a new theory, called $(\tau,b)$-theory, of automorphic forms
for further exploration and investigation.

One of the key ingredients of the $(\tau,b)$-theory is to formulate
explicit constructions of the elliptic endoscopy correspondences (including transfers and descents),
via {\bf automorphic integral transforms} with the kernel functions. Such general constructions extends almost all previous known constructions
in the theory of automorphic forms. Some particular cases were formulated by Ginzburg for $F$-split classical groups
in \cite{G08} and \cite{G12} for constructions of special CAP automorphic representations in the sense of Piatetski-Shapiro,
without referring to the Arthur classification theory (\cite{Ar12}). However, it is much better understood if they are formulated in terms of
Arthur's classification of the discrete spectrum of classical groups (\cite{Ar12} and \cite{Mk12}).
This is the main point of the current paper, which extends
the case of certain orthogonal groups considered in the Oberwolfach report of the author (\cite{Jn11}).

In this introduction, $G$ is assumed to be an $F$-quasisplit classical group, which is not a unitary group. The details about unitary groups
will be covered in the rest of the paper.
Let $\psi\in\wt{\Psi}_2(G)$ be a global Arthur parameter for $G$, following the notation in \cite{Ar12}.
Assume that $\psi$ is not simple. This means that $\psi$ can
be expressed as a formal sum
\begin{equation}\label{apsi}
\psi=\psi_1\boxplus\psi_2,
\end{equation}
and there are two classical groups $G_0=G(\psi_1)$ and $H=G(\psi_2)$ determined by $\psi_1$ and $\psi_2$, respectively, such that
$\psi_1\in\wt{\Psi}_2(G_0)$ and $\psi_2\in\wt{\Psi}_2(H)$ and $G_0\times H$ is an elliptic endoscopy group of $G$. The key point here is to use
$\psi_1$ (and the structures of $H$ and $G$) to construct a family of {\bf automorphic kernel functions}
$\CK_{\psi_1;H,G}(h,g)$, which are automorphic functions on $H(\BA)\times G(\BA)$.

\begin{conj}[Main Conjecture on Endoscopy Correspondences]\label{mcec}
Assume that $\psi=\psi_1\boxplus\psi_2$ is a global Arthur parameter for $G$ such that $G_0=G(\psi_1)$ and $H=G(\psi_2)$ are twisted endoscopy
groups associated to $\psi_1$ and $\psi_2$, respectively, so that $G_0\times H$ is an elliptic endoscopy group of $G$ associated to $\psi$.
Let $\sigma$ and $\pi$ be irreducible automorphic representations of $H(\BA)$ and $G(\BA)$ occurring in the discrete spectrum respectively.
Then there exists a family of automorphic kernel functions $\CK_{\psi_1;H,G}(h,g)$ built from $\psi_1$, and the structures of $H$ and $G$,
such that for $\varphi_\sigma\in\sigma$ and $\varphi_\pi\in\pi$, if the following integral
\begin{equation}\label{mci}
\int_{H(F)\bks H(\BA)}\int_{G(F)\bks G(\BA)}\CK_{\psi_1;H,G}(h,g)\varphi_\sigma(h)\overline{\varphi_\pi}(g)dhdg
\end{equation}
is convergent and nonzero, then $\sigma\in\wt\Pi_{\psi_2}(\epsilon_{\psi_2})$ if and only if
$\pi\in\wt\Pi_{\psi}(\epsilon_\psi)$ where $\wt\Pi_{\psi_2}(\epsilon_{\psi_2})$ and $\wt\Pi_{\psi}(\epsilon_\psi)$ are
the corresponding automorphic $L^2$-packets, the definition of which is given in Section 2.1.
\end{conj}

It is clear that the integral in \eqref{mci} converges absolutely if $\sigma$ and $\pi$ are cuspidal; otherwise, a suitable
regularization (by using Arthur truncation, for instance) is needed. This integral also produces two-way constructions or
integral transforms as follows.

Let $\sigma$ be an irreducible cuspidal automorphic representation of $H(\BA)$. Define
\begin{equation}\label{ethg}
\CT_{\psi_1,H}(\varphi_\sigma)(g)
:=
\int_{H(F)\bks H(\BA)}\CK_{\psi_1;H,G}(h,g)\varphi_\sigma(h)dh.
\end{equation}
Then this integral is absolutely convergent and defines an automorphic function on $G(\BA)$. Denote by
$\CT_{\psi_1,H}(\sigma)$ the space generated by all automorphic functions $\CT_{\psi_1,H}(\varphi_\sigma)(g)$ and call
$\CT_{\psi_1,H}(\sigma)$ the {\bf endoscopy transfer} from $H$ to $G$ via $\psi_1$ or the
{\bf endoscopy transfer} from $G_0\times H$ to $G$. When $\sigma$ is not cuspidal, the integral in \eqref{ethg} may not
be convergent and certain regularization is needed.

Conversely, let $\pi$ be an irreducible cuspidal automorphic
representation of $G(\BA)$. Define
\begin{equation}\label{edgh}
{\mathcal D}_{\psi_1,G}(\varphi_\pi)(h)
:=
\int_{G(F)\bks G(\BA)}\CK_{\psi_1;H,G}(h,g)\overline{\varphi_\pi}(g)dg.
\end{equation}
Then this integral is absolutely convergent and defines an automorphic function on $H(\BA)$. Denote by
${\mathcal D}_{\psi_1,G}(\pi)$ the space generated by all automorphic functions ${\mathcal D}_{\psi_1,G}(\varphi_\pi)(h)$
and call ${\mathcal D}_{\psi_1,G}(\pi)$ the {\bf endoscopy descent} from $G$ to $H$ via $\psi_1$.
When $\pi$ is not cuspidal, the integral in \eqref{edgh} may not
be convergent and certain regularization is needed.

The goal is to investigate the automorphic representations $\CT_{\psi_1,H}(\sigma)$ of $G(\BA)$ and
${\mathcal D}_{\psi_1,G}(\pi)$ of $H(\BA)$, respectively. From Conjecture \ref{mcec}, one may expect to
produce automorphic representations in the global Arthur packet $\wt\Pi_{\psi_1\boxplus\psi_2}$ of $G(\BA)$
from the endoscopy transfer $\CT_{\psi_1,H}(\sigma)$ if $\sigma$ runs in the automorphic $L^2$-packet $\wt\Pi_{\psi_2}(\epsilon_{\psi_2})$;
and to produce automorphic representations in the global Arthur packet $\wt\Pi_{\psi_2}$ of $H(\BA)$ from
the endoscopy descent ${\mathcal D}_{\psi_1,G}(\pi)$ if $\pi$ runs in the automorphic $L^2$-packet
$\wt\Pi_{\psi_1\boxplus\psi_2}(\epsilon_{\psi_1\boxplus\psi_2})$.

According to the structure of global Arthur parameters, any $\psi\in\wt{\Psi}_2(G)$ can be expressed as a formal sum
\begin{equation}\label{ap1}
\psi=(\tau_1,b_1)\boxplus(\tau_2,b_2)\boxplus\cdots\boxplus(\tau_r,b_r)
\end{equation}
of simple global Arthur parameters $(\tau_1,b_1), (\tau_2,b_2),\cdots$, and $(\tau_r,b_r)$, where $\tau_i$ are irreducible unitary self-dual
cuspidal automorphic representations of $\GL_{a_i}(\BA)$ and $b_i\geq 1$ are integers representing the $b_i$-dimensional irreducible
representation $\nu_i$ of $\SL_2(\BC)$ for all $i$. In \cite{Ar12}, $(\tau_i,b_i)$ are replaced by $\tau_i\boxtimes\nu_i$, respectively.
It is natural to start the investigation with $\psi_1=(\tau,b)$ in Conjecture \ref{mcec}, i.e.
$$
\psi=(\tau,b)\boxplus\psi_2.
$$

The $(\tau,b)$-theory of automorphic forms is to understand Conjecture \ref{mcec} with $\psi_1=(\tau,b)$ and related problems, based on
or guided by the Arthur classification of the discrete spectrum for classical groups. The core problem is to understand the spectral
structures of $\CT_{\psi_1,H}(\sigma)$ and ${\mathcal D}_{\psi_1,G}(\pi)$; and to characterize the non-vanishing of
$\CT_{\psi_1,H}(\sigma)$ and ${\mathcal D}_{\psi_1,G}(\pi)$ in terms of certain basic invariants of automorphic forms.
It is not clear how to use the basic invariants from the trace formula, i.e., the stable orbital integrals,  to the study of
Conjecture \ref{mcec}, because technically, it is not clear how to relate the integrals in Conjecture \ref{mcec} to the
corresponding stable orbital integrals, although this is the most natural way to think of in the current setting.

When $\tau$ is a quadratic character $\chi$ of $\GL_1(\BA)$, the $(\chi,b)$-theory is essentially about the theta
correspondence for reductive dual pairs. The basic invariants used by Rallis in his program (\cite{Rl87} and \cite{Rl91}) are the poles of
the twisted standard $L$-functions through the doubling integrals of Piatetski-Shapiro and Rallis (\cite{GPSR87}). For the general
$(\tau,b)$-theory, the first basic invariant considered is the structure of Fourier coefficients of automorphic forms and
automorphic representations. As the theory develops, other invariants like poles or special values of certain family of
automorphic $L$-functions, and various periods of automorphic forms will be employed. In this generality, the
Fourier coefficients of automorphic forms are defined in terms of nilpotent orbits in the corresponding Lie algebra and
are natural invariants of automorphic representations close to the type of integrals in \eqref{mci}.
Section \ref{fcno} is to discuss the basic theory in this aspect.

The Arthur classification of the discrete spectrum for $F$-quasisplit classical groups is reviewed briefly in Section \ref{ap}, where
the basic endoscopy structures of each global Arthur packet is discussed. This leads to the explicit constructions of
the {\bf automorphic kernel functions}. More precisely, the constructions are formulated for the case when
$\psi_1$ is a simple global Arthur parameter, say, of form $(\tau,b)$, where $\tau$ is an irreducible self-dual unitary
cuspidal automorphic representation of $\GL_a(\BA)$ and $b>1$ is an integer. Because of the parity structure of simple Arthur
parameters, this leads to two different types of constructions: one uses the generalized Bessel-Fourier coefficients in Section \ref{ec(1)}
and the other uses the generalized Fourier-Jacobi coefficients in Section \ref{ec(2)}, based on the explicit classification of
endoscopy structure for all $F$-quasisplit classical groups in Section 3.
In Section 6, the discussions are given to address how the theory of automorphic descent method can be extended
from the viewpoint of
the $(\tau,b)$-theory to obtain refined structure about the global Arthur packets of simple type. This includes the work of
Ginzburg, Soudry and the author (\cite{GJS12}), the thesis of Baiying Liu (\cite{Liu13})
and a more recent work of Liu and the author (\cite{JL13}). The particular case of $(\tau,b)$-theory
with $\tau$ being a quadratic character $\chi$ of $\GL_1$ and its relation to the theory of theta correspondence will be discussed in Section 7.
Section 8 will address the general formulation of the $(\tau,b)$-theory and basic problems. The final section is to address other related
constructions and other related topics in the theory of automorphic forms.

The author is strongly influenced by the idea of the automorphic descent method of Ginzburg, Rallis and Soudry on the topics
discussed in the paper, and benefits greatly from his long term collaboration with Ginzburg and Soudry on various research projects
related to topics discussed here. The main idea of the author contributing to the theory is to formulate the general framework of
the constructions of endoscopy correspondences via automorphic integral transforms in terms of the Arthur classification of
the discrete spectrum of $F$-quasisplit classical groups. This idea came up to the author after reading the Clay lecture notes of
Arthur (\cite{Ar05}), who sent it to me before it got published. In particular, this idea of the author got developed during the
two-month visit in Paris in 2006, which was arranged by B.-C. Ngo for a one-month period in University of Paris-Sud and by M. Harris
for a one-month period in University of Paris 7. Since then, the author has been invited to talk about the constructions and the related theory
at many places, including Paris, Jerusalem, Toronto, New York, Vienna, Oberwolfach, Beijing, Philadelphia, Chicago, Boston, New Haven, and so on.
The author would like to take this opportunity to express his appreciation to all the invitations and supports he received over the years.
In particular, the author would like to thank J. Arthur, M. Harris, C. Moeglin, B.-C. Ngo, and D. Vogan for helpful conversations and suggestions,
to thank D. Ginzburg and D. Soudry for productive collaboration on relevant research projects over the years, and to thank Baiying Liu
and Lei Zhang for joining this project recently, which brings progress in various aspects of the project.
The guidance of I. Piatetski-Shapiro and S. Rallis which brought the author to this subject can be clearly noticed throughout the whole paper and
will continue to influence him in his future work. The work of the author has been supported in part by the NSF Grant
DMS--0400414, DMS--0653742, and DMS--1001672, over the years. Finally, the author would like to thank R. Howe for his leadership in
the organization of the conference: Automorphic Forms and Related Geometry: assessing the legacy of I.I. Piatetski-Shapiro, at Yale
University, April, 2012, which makes it possible to write up this paper for the conference proceedings.

\section{Arthur Parametrization of Discrete Spectrum}\label{ap}

Let $F$ be a number field and $E$ be a quadratic extension of $F$, whose Galois group is denoted by $\Gamma_{E/F}=\{1,c\}$.
Take $G_n$ to be either the $F$-quasisplit unitary groups
$\RU_{2n}=\RU_{E/F}(2n)$ or $\RU_{2n+1}=\RU_{E/F}(2n+1)$ of hermitian type, the $F$-split special orthogonal group $\SO_{2n+1}$, the symplectic group
$\Sp_{2n}$, or the $F$-quasisplit even special orthogonal group $\SO_{2n}$.
Let $F'$ be a number field, which is $F$ if $G_n$ is not a unitary group, and
is $E$ if $G_n$ is a unitary group. Denote by $\RR_{F'/F}(n):=\RR_{F'/F}(\GL_{n})$ the Weil restriction of the $\GL_{n}$ from $F'$ to $F$.

Following Chapter 1 of Arthur's book (\cite{Ar12}) and also the work of C.-P. Mok (\cite{Mk12}),
take the closed subgroup $G(\BA)^1$ of $G(\BA)$ given by
\begin{equation}
G(\BA)^1:=\{x\in G(\BA)\mid |\chi(x)|_\BA=1, \forall \chi\in X(G)_F\}
\end{equation}
where $X(G)_F$ is the group of all $F$-rational characters of $G$. Since $G(F)$ is a discrete subgroup of $G(\BA)^1$ and
the quotient $G(F)\bks G(\BA)^1$ has finite volume with respect to the canonical Haar measure, consider the
the space of all square integrable functions on $G(F)\bks G(\BA)^1$, which is denoted by
$L^2(G(F)\bks G(\BA)^1)$. It has the following embedded, right $G(\BA)^1$-invariant Hilbert subspaces
\begin{equation}
L^2_{\cusp}(G(F)\bks G(\BA)^1)\subset L^2_{\disc}(G(F)\bks G(\BA)^1)\subset L^2(G(F)\bks G(\BA)^1)
\end{equation}
where $L^2_{\disc}(G(F)\bks G(\BA)^1)$ is called the discrete spectrum of $G(\BA)^1$, which decomposes into
a direct sum of irreducible representations of $G(\BA)^1$; and $L^2_{\cusp}(G(F)\bks G(\BA)^1)$ is called the
cuspidal spectrum of $G(\BA)^1$, which consists of all cuspidal automorphic functions on $G(\BA)^1$ and is
embedded into $L^2_{\disc}(G(F)\bks G(\BA)^1)$.

As in Section 1.3 of \cite{Ar12}, denote by $\CA_\cusp(G)$, $\CA_2(G)$, and $\CA(G)$ respectively,
the set of irreducible unitary representations of $G(\BA)$ whose restriction to $G(\BA)^1$ are irreducible constituents of
$L^2_{\cusp}(G(F)\bks G(\BA)^1)$, $L^2_{\disc}(G(F)\bks G(\BA)^1)$, and $L^2(G(F)\bks G(\BA)^1)$, respectively.
It is clear that
\begin{equation}
 \CA_\cusp(G)\subset \CA_2(G)\subset \CA(G).
\end{equation}
In particular, when $G=\RR_{F'/F}(N)=\RR_{F'/F}(\GL_N)$, the Weil restriction, the notation is simplified as follows: $\CA_\cusp(N):=\CA_\cusp(G)$, $\CA_2(N):=\CA_2(G)$, and
$\CA(N):=\CA(G)$.

\subsection{Discrete spectrum of classical groups}
Recall more notation from Sections 1.3 and 1.4 of \cite{Ar12} and Section 2.3 of \cite{Mk12}
in order to state Arthur's theorem in more explicit way.
$\CA_\cusp^+(G)$ and $\CA_2^+(G)$ are defined in the same way as $\CA_\cusp(G)$ and $\CA_2(G)$, respectively, except
that the irreducible representations $\pi$ of $G(\BA)$ may not be unitary.
The set $\CA^+(N)$ can be easily defined in terms of the set $\CA(N)$. The elements in $\CA(N)$ can be explicitly given as follows.

Let $N=\sum_{i=1}^rN_i$ be a any partition of $N$ and $P=M_PN_P$ be the corresponding standard (upper-triangular)
parabolic subgroup of $\GL_N$. For each $i$, take $\pi_i\in\CA_2(N_i)$ and form an induced representation
\begin{equation}\label{gind}
\pi:=\CI_P(\pi_1\otimes\pi_2\otimes\cdots\otimes\pi_r),
\end{equation}
which is irreducible and unitary. The automorphic realization of $\pi$ is given by the meromorphic continuation of the Eisenstein series
associated to the induced representation.
The set $\CA(N)$ consists of all elements of type \eqref{gind}. Then the set
$\CA^+(N)$ consists of all elements of type as in \eqref{gind}, but with $\pi_i$ belonging to $\CA_2^+(N_i)$ for
$i=1,2,\cdots,r$. Hence one has
\begin{equation}\label{an}
\CA_\cusp(N)\subset\CA_2(N)\subset\CA(N)\subset\CA^+(N).
\end{equation}
The corresponding sets of global Arthur parameters are given by
\begin{equation}\label{gapn}
\Psi_\cusp(N)\subset\Psi_2(N)=\Psi_\simp(N)\subset\Psi_\el(N)\subset\Psi(N)\subset\Psi^+(N).
\end{equation}
More precisely, the set $\Psi_2(N)=\Psi_\simp(N)$ consists of all pairs $(\tau,b)$ with $N=ab$ ($a,b\geq 1$) and
$\tau\in\CA_\cusp(a)$. The subset $\Psi_\cusp(N)$ consists of all pairs $(\tau,1)$, i.e. with $b=1$ and $N=a$ and
$\tau\in\CA_\cusp(N)$. The set $\Psi(N)$ consists of all possible formal (isobaric) sum
\begin{equation}\label{psin}
\psi:=\psi_1\boxplus\psi_2\boxplus\cdots\boxplus\psi_r,
\end{equation}
where $\psi_i$ belongs to $\Psi_\simp(N_i)$ with any partition $N=\sum_{i=1}^rN_i$ and hence $\psi_i=(\tau_i,b_i)$ with
$N_i=a_ib_i$ and $\tau_i\in\CA_\cusp(a_i)$. The elements of $\Psi^+(N)$ are of the form \eqref{psin} with
$\tau_i\in\CA_\cusp^+(a_i)$. The set $\Psi_\el(N)$ consists of all elements of the form \eqref{psin} with the property that
$\psi_i$ and $\psi_j$ are not equal if $i\neq j$. The equality of $\psi_i$ and $\psi_j$ means that $b_i=b_j$ and $\tau_i\cong\tau_j$.

For $\psi\in\Psi(N)$ as given in \eqref{psin}, define the dual of $\psi$ by
\begin{equation}\label{dgapn}
\psi^*:=\psi_1^*\boxplus\psi_2^*\boxplus\cdots\boxplus\psi_r^*,
\end{equation}
where $\psi_i^*:=(\tau_i^*,b_i)$ with $\tau_i^*$ being the dual $(\tau_i^c)^\vee$ of the Galois twist $\tau_i^c$ of $\tau_i$,
for $i=1,2,\cdots,r$. Denote
by $\wt{\Psi}(N)$ the subset of $\Psi(N)$ consisting of all conjugate self-dual parameters, i.e. $\psi^*=\psi$. Hence one has
\begin{eqnarray}
\wt{\Psi}_\cusp(N)&:=&\wt{\Psi}(N)\cap\Psi_\cusp(N),\\
\wt{\Psi}_\simp(N)&:=&\wt{\Psi}(N)\cap\Psi_\simp(N),\\
\wt{\Psi}_\el(N)&:=&\wt{\Psi}(N)\cap\Psi_\el(N).
\end{eqnarray}
By Theorem 1.4.1 of \cite{Ar12} and Theorem 2.4.2 of \cite{Mk12},
for each conjugate self-dual, simple, generic global Arthur parameter $\psi$ in $\wt{\Psi}_\cusp(N)$,
there is a unique (up to isomorphism class of endoscopy data) endoscopy group $G_\psi$ in $\wt{\CE}_\el(N)$, the
set of twisted elliptic endoscopy groups of $\RR_{F'/F}(N)$, with
the weak Langlands functorial transfer from $G_\psi$ to $\RR_{F'/F}(N)$.
Note that $G_\psi$ is a simple algebraic group defined over $F$ which is either an $F$-quasisplit unitary group, a symplectic group
or an $F$-quasisplit orthogonal group. The set of all these simple classical groups is denoted by $\wt{\CE}_\simp(N)$.
The theme of \cite{Ar12} and \cite{Mk12} is to classify the discrete spectrum for these classical groups $G$ in $\wt{\CE}_\simp(N)$.

Note that when $F'=E$ and $G_\psi=\RU_{E/F}(N)$, the endoscopy data contains the $L$-embedding $\xi_{\chi_\kappa}$ from the
$L$-group $^LG_\psi$ to the $L$-group of $\RR_{E/F}(N)$, which depends on the sign $\kappa=\pm$. The two embeddings are not equivalent.
This will be more specifically discussed in Section where the endoscopy groups of $\RU_{E/F}(N)$ are discussed.
The general theory of twisted endoscopy is referred to \cite{KtS99}.

For each $G\in\wt{\CE}_\simp(N)$, denote by $\wt{\Psi}(G)$ the set of global conjugate self-dual Arthur parameters $\psi$ in
$\wt{\Psi}(N)$, which factor through the Langlands dual group $^LG$ in the sense of Arthur (Page 31 in \cite{Ar12} and also \cite{Mk12}).
When $\psi$ belongs to $\wt{\Psi}_\el(N)$, the above property defines $\wt{\Psi}_2(G)$.
Similarly, for each $G\in\wt{\CE}_\el(N)$ (the twisted elliptic endoscopy groups of $\RR_{F'/F}(N)$), one defines the set
of global Arthur parameters $\wt{\Psi}_2(G)$. As remarked in Pages 31--32 of \cite{Ar12} (and also in \cite{Mk12})
the elements of $\wt{\Psi}(G)$
may be more precisely written as a pair $(\psi,\wt{\psi}_G)$ and the projection
$$
(\psi,\wt{\psi}_G)\mapsto \psi
$$
is not injective in general. However, as explained in Page 33 of \cite{Ar12} (and also in \cite{Mk12}), for $\psi\in\wt{\Psi}_\el(N)$, there exists
only one $G\in\wt{\CE}_\el(N)$ and one $\psi_G\in\wt{\Psi}_2(G)$ projects to $\psi$, and hence the set $\wt{\Psi}_2(G)$ can be
viewed as a subset of $\wt{\Psi}_\el(N)$. Moreover, one has the following disjoint union:
$$
\wt{\Psi}_\el(N)=\cup_{G\in\wt{\CE}_\el(N)}\wt{\Psi}_2(G).
$$
Note that the set $\wt{\Psi}_2(G)$ is identified as a subset of $\wt{\Psi}_\el(N)$ via the $L$-embedding $\xi_{\chi_\kappa}$ in the
endoscopy data $(G,\xi_{\chi_\kappa})$.

Theorem 1.5.2 of \cite{Ar12} and Theorem 2.5.2 of \cite{Mk12} give
the following decomposition of the discrete spectrum of $G\in\wt{\CE}_\simp(N)$, which is a refinement of Theorem \ref{gpsl}:
\begin{equation}\label{arthur}
L^2_\disc(G(F)\bks G(\BA))
\cong
\oplus_{\psi\in\wt{\Psi}_2(G)}\oplus_{\pi\in\wt{\Pi}_\psi(\epsilon_\psi)}m_\psi\pi,
\end{equation}
where $\epsilon_\psi$ is a linear character defined as in Theorem 1.5.2 of \cite{Ar12} in terms of $\psi$ canonically, and
$\wt{\Pi}_\psi(\epsilon_\psi)$ is the subset of the global Arthur packet $\wt{\Pi}_\psi$ consisting of
all elements $\pi=\otimes_v\pi_v$ whose characters are equal to $\epsilon_\psi$, and finally, $m_\psi$ is the multiplicity
which is either $1$ or $2$, depending only on the global Arthur parameter $\psi$. See Theorem 1.5.2 of \cite{Ar12}
and Theorem 2.5.2 of \cite{Mk12} for details.

Note that in the case of unitary groups, the set $\wt{\Psi}_2(G)$ depends on
a given embedding $\xi_{\chi_\kappa}$ of $L$-groups, as described in \cite[Section 2.1]{Mk12}. Hence one should use the notation
$\wt{\Psi}_2(G,\xi_{\chi_\kappa})$ for $\wt{\Psi}_2(G)$. If $\xi_{\chi_\kappa}$ is fixed in a discussion, it may be dropped from the
notation if there is no confusion.

It is clear that the subset $\wt{\Pi}_\psi(\epsilon_\psi)$ of the global Arthur packet $\wt{\Pi}_\psi$ consists of all members in
$\wt{\Pi}_\psi$ which occur in the discrete spectrum of $G$. It is discussed in \cite{Ar11} and \cite[Section 4.3]{Ar12} that
$\wt{\Pi}_\psi(\epsilon_\psi)$ is in fact the subset of \cite[Section 4.3]{Ar12} consisting of all members which may occur in the whole
automorphic $L^2$-spectrum of $G$. As remarked in \cite[Section 4.3]{Ar12}, it is a long term study problem to
understand other possible automorphic representations in the given global Arthur packets \cite[Section 4.3]{Ar12}.
In this paper, for simplicity, the set $\wt{\Pi}_\psi(\epsilon_\psi)$ is called the {\bf automorphic $L^2$-packet} attached to the global Arthur parameter $\psi$.

\subsection{Arthur-Langlands transfers}
For each $G\in\wt{\CE}_\simp(N)$, by the definition of twisted elliptic (simple type) endoscopy data, there is $L$-homomorphism
$$
\xi_{\chi_\kappa}=\xi_\psi:\ {^LG}\rightarrow{^L\RR_{F'/F}(N)},
$$
where $\RR_{F'/F}(N)=\RR_{F'/F}(\GL_N)$.
The corresponding Langlands functorial transfer from $G$ to $\RR_{F'/F}(N)$ is given as follows.

Take $\psi\in\wt{\Psi}_2(G)$. The theory of stable trace formula (\cite{Ar12}) produces the automorphic $L^2$-packet
$\wt{\Pi}_\psi(\epsilon_\psi)$, which consists of members in the global Arthur packet $\wt{\Pi}_\psi$ that belong to $\CA_2(G)$.
On the other hand, $\psi$ can be expressed as in \eqref{psin}
$$
\psi:=\psi_1\boxplus\psi_2\boxplus\cdots\boxplus\psi_r,
$$
where $\psi_i$ belongs to $\wt{\Psi}_\simp(N_i)$ with a partition $N=\sum_{i=1}^rN_i$ and hence $\psi_i=(\tau_i,b_i)$ with
$N_i=a_ib_i$ and $\tau_i\in\CA_\cusp(a_i)$ being conjugate self-dual. For each $\psi_i\in\wt{\Psi}_\simp(N_i)$, by the theorem of Moeglin and
Waldspurger (\cite{MW89} and Theorem 1.3.3 of \cite{Ar12}), there is a unique (up to isomorphism)
irreducible residual representation
$\pi_i$ of $\RR_{F'/F}(N_i)(\BA)$ which is generated by the Speh residues of the Eisenstein series with cuspidal support
$$
(\RR_{F'/F}(a_i)^{\times b_i},\tau_i^{\otimes b_i}),
$$
and is sometimes denoted by $\Delta(\tau_i,b_i)$.
Then one defines the Arthur representation of $\RR_{F'/F}(N)(\BA)$,
\begin{equation}\label{pipsi}
\pi_\psi:=\CI_P(\pi_1\otimes\pi_2\otimes\cdots\otimes\pi_r),
\end{equation}
which is induced from the parabolic subgroup of $\RR_{F'/F}(N)$ attached to the
partition $N=\sum_{i=1}^rN_i$.
It is an irreducible unitary automorphic representation of $\RR_{F'/F}(N)(\BA)$, whose automorphic realization is given by the Langlands
theory on meromorphic continuation of residual Eisenstein series on $\RR_{F'/F}(N)(\BA)$.

The Arthur-Langlands transfer is a Langlands functorial transfer corresponding to $\xi_\psi$, which takes the set $\wt{\Pi}_\psi(\epsilon_\psi)$
to the automorphic representation $\pi_\psi$, in the sense that for each $\pi\in\wt{\Pi}_\psi(\epsilon_\psi)$, the Satake parameter at each
unramified local place $v$, $c(\pi_v)$, matches the corresponding Satake parameter $c(\pi_{\psi,v})$, i.e.
$$
\xi_{\psi_v}(c(\pi_v))=c(\pi_{\psi,v}).
$$
This gives the functorial interpretation of the Arthur decomposition \eqref{arthur} for the discrete spectrum $L^2_\disc(G(F)\bks G(\BA))$
for each $G\in\wt{\CE}_\simp(N)$, which can be presented by the following diagram:
$$
\begin{matrix}
\\
              &  & &\wt{\Psi}_2(G)& & &\\
               & &&&&&\\
             &   &&\psi&&&\\
             &   &&&&&\\
             &   &\swarrow& &\searrow&&\\
             &   &&&&&\\
\CA_2(G) & \wt{\Pi}_{\psi}(\epsilon_\psi)&       & {\longrightarrow}        &           & \pi_\psi& \CA(N)\\
\\
\end{matrix}
$$

\begin{rmk}
It is an important problem to study when the Arthur-Langlands transfer from any member $\pi$ in the
automorphic $L^2$-packet $\wt{\Pi}_{\psi}(\epsilon_\psi)$
to the Arthur representation $\pi_\psi$ is compatible with the local Langlands functorial transfers at all local places. In particular, when the global Arthur
parameter $\psi\in \wt{\Psi}_2(G)$ is not generic, it is in general to expect that the image under the Arthur-Langlands
transfer of all members in $\wt{\Pi}_{\psi}(\epsilon_\psi)$ contains automorphic representations other than $\pi_\psi$. It is an even more
subtle problem to figure out what automorphic representations should be expected in the transfer image of the automorphic $L^2$-packet
$\wt{\Pi}_{\psi}(\epsilon_\psi)$ for a given non-generic global Arthur parameter $\psi\in \wt{\Psi}_2(G)$.

On the other hand, it is also important to figure out what kind automorphic representations are expected to be in an automorphic $L^2$-packet
$\wt{\Pi}_{\psi}(\epsilon_\psi)$ for a given global Arthur parameter $\psi\in \wt{\Psi}_2(G)$.
Moeglin provides conditions and conjectures in \cite{Mg11} on whether an automorphic $L^2$-packet
$\wt{\Pi}_{\psi}(\epsilon_\psi)$ contains nonzero residual representations, while the thesis
of Paniagua-Taboada in \cite{PT11}
provides cases that an automorphic $L^2$-packet $\wt{\Pi}_{\psi}(\epsilon_\psi)$ may contains no cuspidal members for $\SO_{2n}$.
The work of Liu and Zhang joint with the author in \cite{JLZ12} constructs residual representations in certain automorphic
$L^2$-packets $\wt{\Pi}_{\psi}(\epsilon_\psi)$ for $F$-quasisplit classical groups.
\end{rmk}

\section{Endoscopy Structures for Classical Groups}

For each $(G,\xi_{\chi_\kappa})\in\wt{\CE}_\simp(N)$ and for each $\psi\in\wt{\Psi}_2(G)$, the decomposition as in \eqref{psin}
$$
\psi:=\psi_1\boxplus\psi_2\boxplus\cdots\boxplus\psi_r
$$
also implies an endoscopy structure for $\psi$ and hence for the automorphic representations in $\wt{\Pi}_\psi(\epsilon_\psi)$.

Recall that $\psi_i$ belongs to $\wt{\Psi}_\simp(N_i)$ with the partition $N=\sum_{i=1}^rN_i$ and that $\psi_i=(\tau_i,b_i)$ with the property that
$N_i=a_ib_i$ and $\tau_i\in\CA_\cusp(a_i)$ is {\it conjugate self-dual}, i.e. $\tau_i^*=(\tau_i^c)^\vee=\tau_i$.
For each $\psi_i\in\wt{\Psi}_\simp(N_i)$, there is a unique twisted elliptic (simple type)
endoscopy group $G_i$ (as part of the endoscopy datum) in $\wt{\CE}_\simp(N_i)$. Hence for $\psi$, one may consider
$$
G_1\times G_2\times\cdots\times G_r.
$$
One may call this product of groups an (generalized) elliptic endoscopy group of $G$ determined by $\psi$.

When $F'=E$ and $G=\RU(N)=\RU_{E/F}(N)$, one should add the $L$-embedding
$$
\zeta_{\udl{\chi}_{\udl{\kappa}}}\ :\ {^L(\RU(N_1)\times\cdots\times\RU(N_r))}\rightarrow{^L\RU(N)}
$$
as given in \cite[(2.1.13)]{Mk12}, and the sign $\udl{\kappa}=(\kappa_1,\cdots,\kappa_r)$ is defined as in
\cite[(2.1.12)]{Mk12} with $\kappa_i=(-1)^{N-N_i}$ for $i=1,2,\cdots,r$. The $L$-embeddings from the two $L$-groups to
the $L$-group $^L\RR_{E/F}(N)$ are given as in \cite[(2.1.14)]{Mk12}.

The main constructions of endoscopy correspondences considered in this paper are to take the following decomposition of $\psi$:
\begin{equation}\label{psi2}
\psi=(\tau,b)\boxplus\psi_{N-ab}
\end{equation}
with $\tau\in\CA_\cusp(a)$ being conjugate self-dual, and $\psi_{N-ab}\in\wt{\Psi}_\el(N-ab)$. Note that $(\tau,b)\in\wt{\Psi}_\simp(ab)$.
Take $G_0\in\wt{\CE}_\simp(ab)$ and $H\in\wt{\CE}_\simp(N-ab)$ and consider the endoscopy structure
\begin{equation}\label{end2}
G_0\times H \rightarrow G,
\end{equation}
with $(\tau,b)\in\wt{\Psi}_2(G_0)$ and $\psi_{N-ab}\in\wt{\Psi}_2(H)$.
The explicit constructions for the endoscopy transfers associated to \eqref{psi2} will be discussed below case by case and
are expected to confirm Conjecture \ref{mcec} for those cases.
For $\tau\in\CA_\cusp(a)$, one may use the notation that $a=a_\tau$ to indicate the relation between
$\RR_{F'/F}(a)$ and $\tau$. The explicit constructions for the endoscopy correspondences are distinguished
according to the types of the functorial transfer of $\tau$ and the parity of $a_\tau$, which will be discussed with details below.

When $F'=F$, a self-dual $\tau\in\CA_\cusp(a)$ is called {\it of symplectic type} if the (partial) exterior square $L$-function
$L^S(s,\tau,\wedge^2)$ has a (simple) pole at $s=1$; otherwise, $\tau$ is called {\it of orthogonal type}. In the latter case,
the (partial) symmetric square $L$-function $L^S(s,\tau,\sym^2)$ has a (simple) pole at $s=1$.

When $F'=E$, the quadratic field extension of $F$, a conjugate self-dual $\tau\in\CA_\cusp(a)$ yields a simple generic parameter
$(\tau,1)$ in $\wt\Psi_2(a)$ which produces a sign $\kappa_a$ for the endoscopy data $(\RU_{E/F}(a),\xi_{\chi_{\kappa_a}})$. By Theorem
2.5.4 of \cite{Mk12}, the (partial) $L$-function
$$
L(s,(\tau,1),\Asai^{\eta_{(\tau,1)}})
$$
has a (simple) pole at $s=1$ with the sign $\eta_{(\tau,1)}=\kappa_a\cdot(-1)^{a-1}$
(see also \cite[Theorem 8.1]{GGP12} and \cite[Lemma 2.2.1]{Mk12}).
Then the irreducible cuspidal automorphic representation $\tau$ or equivalently the simple generic parameter $(\tau,1)$
is called {\it conjugate orthogonal} if $\eta_{(\tau,1)}=1$ and {\it conjugate symplectic} if $\eta_{(\tau,1)}=-1$, following the terminology of
\cite[Section 3]{GGP12} and \cite[Section 2]{Mk12}.
Here
$L^S(s,(\tau,1),\Asai^+)$ is the (partial) Asai $L$-function of $\tau$ and $L^S(s,(\tau,1),\Asai^-)$ is the (partial) Asai $L$-function of
$\tau\otimes\omega_{E/F}$, where $\omega_{E/F}$ is the quadratic character associated to $E/F$ by the global classfield theory.

The sign of $(\tau,b)\in\wt\Psi_2(ab)$ can be calculated following \cite[Section 2.4]{Mk12}. Fix the sign $\kappa_a$ as before for the
endoscopy data $(\RU_{E/F}(a),\xi_{\chi_{\kappa_a}})$, the sign of $(\tau,1)$ is
$$
\eta_{(\tau,1)}=\eta_\tau=\kappa_a(-1)^{a-1}.
$$
Hence the sign of $(\tau,b)$ is given by
$$
\eta_{(\tau,b)}=\kappa_a(-1)^{a-1+b-1}=\kappa_a(-1)^{a+b}=\eta_\tau(-1)^{b-1}.
$$
As in \cite[(2.4.9)]{Mk12}, define
$$
\kappa_{ab}:=\kappa_a(-1)^{ab-a-b+1}.
$$
Then $\kappa_{ab}(-1)^{ab-1}=\eta_\tau(-1)^{b-1}=\eta_{(\tau,b)}$ and hence $\kappa_{ab}=\eta_\tau(-1)^{(a-1)b}$, which
gives the endoscopy data $(\RU_{E/F}(ab),\xi_{{\kappa_{ab}}})$.

\subsection{Endoscopy structure for $\SO_{2n+1}$}
Assume that $G\in\wt{\CE}_\simp(N)$ is a $F$-rational form of $\SO_{2n+1}$. Then $G$ must be the $F$-split $\SO_{2n+1}$, $N=2n$ and $G^\vee=\Sp_{2n}(\BC)$.

In this case, the Langlands dual group $^LG$ is a direct product of
the complex dual group $G^\vee$ and the Galois group $\Gamma_F=\Gal(\overline{F}/F)$.
For each $\psi\in\wt{\Psi}_2(\SO_{2n+1})$, which is of form
\begin{equation}\label{psioso}
\psi=(\tau,b)\boxplus\psi_{2n-ab},
\end{equation}
with $\tau\in\CA_\cusp(a)$ being self-dual and $b>0$, such that $(\tau,b)$ is of symplectic type, and with
$\psi_{2n-ab}$ being a $(2n-ab)$-dimensional elliptic global Arthur parameter of symplectic type, there are following
two cases, which have to be distinguished when constructing the endoscopy correspondences by integral transforms.

{\bf Case (a=2e):}\ When $a=a_\tau=2e$, if $\tau$ is of symplectic type, then $b=2l+1$; and if $\tau$ is of orthogonal type,
then $b=2l$. The endoscopy group of $\SO_{2n+1}$ associated to the global elliptic Arthur parameter
$\psi$ as in \eqref{psioso} is
\begin{equation}\label{b1}
\SO_{ab+1}(\BA)\times\SO_{2n-ab+1}(\BA)\rightarrow \SO_{2n+1}(\BA).
\end{equation}
Following the formulation of Wen-Wei Li (\cite{Liw11}), one can have the following variants:
$$
\SO_{ab+1}(\BA)\times\tilsp_{2n-ab}(\BA)\rightarrow\tilsp_{2n}(\BA),
$$
and
$$
\tilsp_{ab}(\BA)\times\tilsp_{an-ab}(\BA)\rightarrow\tilsp_{2n}(\BA).
$$

{\bf Case (a=2e+1):}\
When $a=a_\tau=2e+1$, then $\tau$ is of orthogonal type and $b=2l$. The endoscopy group of $\SO_{2n+1}$ associated to
the global elliptic Arthur parameter $\psi$ as in \eqref{psioso} is also
\begin{equation}\label{b21}
\SO_{ab+1}(\BA)\times\SO_{2n-ab+1}(\BA)\rightarrow \SO_{2n+1}(\BA).
\end{equation}
In fact, the explicit constructions of the integral transforms for this case use the following
endoscopy group of $\SO_{2n+1}$ associated to the global elliptic Arthur parameter $\psi$:
\begin{equation}\label{b22}
\tilsp_{ab}(\BA)\times\tilsp_{2n-ab}(\BA)\rightarrow \SO_{2n+1}(\BA).
\end{equation}
Note that this endoscopy structure is not given from \cite{Ar12}, but is given in \cite{Liw11}. By the duality in the
sense of Howe, one may consider in this case the endoscopy group of $\tilsp_{2n}(\BA)$:
\begin{equation}\label{b23}
\SO_{ab+1}(\BA)\times\SO_{2n-ab+1}(\BA)\rightarrow \tilsp_{2n}(\BA).
\end{equation}
Hence the constructions discussed for those cases would be natural extensions of the classical theta correspondences for reductive
dual pairs. More details will be given in Section 7.

\subsection{Endoscopy structure for $\Sp_{2n}$}
Assume that $G\in\wt{\CE}_\simp(N)$ is equal to $\Sp_{2n}$. Then $N=2n+1$ and $G^\vee=\SO_{2n+1}(\BC)$. In this case,
the Langlands dual group $^LG$ is a direct product of the complex dual group $G^\vee$ and the Galois group
$\Gamma_F=\Gal(\overline{F}/F)$. For each $\psi\in\wt{\Psi}_2(\Sp_{2n})$, which is of form
\begin{equation}\label{psisp}
\psi=(\tau,b)\boxplus\psi_{2n+1-ab},
\end{equation}
with $\tau\in\CA_\cusp(a)$ being self-dual and $b>0$, such that $(\tau,b)$ is of orthogonal type and with
$\psi_{2n+1-ab}$ is a $(2n-ab+1)$-dimensional elliptic global Arthur parameter of orthogonal type, there are following
two cases, which have to be distinguished when constructing the endoscopy correspondences by integral transforms.

{\bf Case (a=2e):}\ Assume that $(\tau,b)$ is of orthogonal type with even dimension. In this case, one must have that $a=a_\tau=2e$.
If $\tau$ is of symplectic type, then $b=2l$ must be even, and if $\tau$ is of orthogonal type, then $b=2l+1$.
The endoscopy group of $\Sp_{2n}$ associated to the global elliptic Arthur parameter $\psi$ as in \eqref{psisp} is
\begin{equation}\label{c1}
\SO_{ab}(\BA)\times\Sp_{2n-ab}(\BA)\rightarrow \Sp_{2n}(\BA).
\end{equation}
Note that the orthogonal group $\SO_{ab}$ may not be $F$-split.

{\bf Case (a=2e+1):}\ Assume that $(\tau,b)$ is of orthogonal type with odd dimension. In this case, one must have that $a=a_\tau=2e+1$,
$\tau$ is of orthogonal type and $b=2l+1$. The endoscopy group of $\Sp_{2n}$ associated to the global elliptic Arthur parameter
$\psi$ as in \eqref{psisp} is
\begin{equation}\label{c2}
\Sp_{ab-1}(\BA)\times\SO_{2n-ab+1}(\BA)\rightarrow \Sp_{2n}(\BA).
\end{equation}
Again the orthogonal group $\SO_{2n-ab+1}$ may not be $F$-split.

\subsection{Endoscopy structure for $\SO_{2n}$}
Assume that $G\in\wt{\CE}_\simp(N)$ is an $F$-quasisplit even special orthogonal group $\SO_{2n}(q_V)$, which
is denoted by $\SO_{2n}(\eta)$ in \cite{Ar12}. Then $N=2n$ and $G^\vee=\SO_{2n}(\BC)$.

Following \cite[Section 1.2]{Ar11}, for $G=\SO_{2n}(\eta)$, the set of elliptic endoscopy groups of $G$ is parameterized by
pairs of even integers $(N_1',N_2')$ with $0\leq N_1'\leq N_2'$ and $N=N_1'+N_2'$ and by pairs $(\eta_1',\eta_2')$ of quadratic
characters of the Galois group $\Gamma_F$ with $\eta_V=\eta_1'\eta_2'$. The corresponding elliptic endoscopy groups are $F$-quasisplit groups
\begin{equation}\label{soesee}
G'=\SO_{N_1'}(\eta_1')\times\SO_{N_2'}(\eta_2'),
\end{equation}
with dual groups
$$
G'^\vee=\SO_{N_1'}(\BC)\times\SO_{N_2'}(\BC)\subset\SO_N(\BC)=G^\vee.
$$

The group $G=\SO_{N}(\eta)$ has also twisted elliptic endoscopy groups, which are parameterized by pairs of odd integers $(\wt{N}_1',\wt{N}_2')$
with $0\leq \wt{N}_1'\leq \wt{N}_2'$ and $N=\wt{N}_1'+\wt{N}_2'$, and by pairs of quadratic characters $(\wt{\eta}_1',\wt{\eta}_2')$ of
the Galois group $\Gamma_F$ with $\eta=\wt{\eta}_1'\cdot\wt{\eta}_2'$. The corresponding twisted elliptic endoscopy groups are
\begin{equation}\label{soetee}
\wt{G}'=\Sp_{\wt{N}_1'-1}\times\Sp_{\wt{N}_2'-1}
\end{equation}
with dual groups
$$
\wt{G}'^\vee=\SO_{\wt{N}_1'}(\BC)\times\SO_{\wt{N}_2'}(\BC)\subset\SO_N(\BC)=G^\vee.
$$

For the purpose of this paper, the above general structures will be specialized as follows.
For each $\psi\in\wt{\Psi}_2(\SO_{2n})$, which is of form
\begin{equation}\label{psieso}
\psi=(\tau,b)\boxplus\psi_{2n-ab},
\end{equation}
with $\tau\in\CA_\cusp(a)$ being self-dual and $b>0$ such that $(\tau,b)$ is of orthogonal type, and with
$\psi_{2n-ab}$ is a $(2n-ab)$-dimensional elliptic global Arthur parameter of orthogonal type, there are following
two cases, which have to be distinguished when constructing the endoscopy correspondences by integral transforms.

{\bf Case (a=2e):}\ Assume that $(\tau,b)$ is of orthogonal type with even dimension. In this case, one must have that $a=a_\tau=2e$.
If $\tau$ is of symplectic type, then $b=2l$ must be even; and if $\tau$ is of orthogonal type, then $b=2l+1$.
The standard elliptic endoscopy group of $\SO_{2n}$ associated to the global elliptic Arthur parameter $\psi$ as in \eqref{psieso} is
\begin{equation}\label{d1}
\SO_{ab}(\BA)\times\SO_{2n-ab}(\BA)\rightarrow \SO_{2n}(\BA).
\end{equation}
Note that the orthogonal groups here may not be $F$-split.

{\bf Case (a=2e+1):}\ Assume that $(\tau,b)$ is of orthogonal type with odd dimension. In this case, one must have that $a=a_\tau=2e+1$,
$\tau$ is of orthogonal type and $b=2l+1$. The twisted elliptic endoscopy group of $\SO_{2n}$ associated to
the global elliptic Arthur parameter $\psi$ as in \eqref{psieso} is
$$
\Sp_{ab-1}(\BA)\times\Sp_{2n-ab-1}(\BA)\rightarrow \SO_{2n}(\BA).
$$
Note that the orthogonal groups here may not be $F$-split.

\subsection{Endoscopy structure for $\RU_{E/F}(N)$}
The (standard) elliptic endoscopy groups of the $F$-quasisplit unitary groups $\RU_{E/F}(N)=\RU(N)$
are determined by J. Rogawski in \cite[Section 4.6]{Rg90}.
The discussion here follows from \cite[Section 2]{Mk12} for consistence of notation.

Assume that $\RU(N)$ be the $F$-quasisplit unitary group defined by an non-degenerate hermitian form $q_V$ of $N$ variables associated to the
quadratic field extension $E/F$. As a (twisted) endoscopy data of $\RR_{E/F}(N)$, one needs to keep the sign $\kappa$ in
$(\RU(N),\xi_{\chi_\kappa})$.

The set of (standard) elliptic endoscopy groups (or data) of $\RU(N)$ are given as in \cite[Section 2.4]{Mk12} and also in \cite{JZ13} by
\begin{equation}\label{usee}
(\RU(N_1)\times\RU(N_2),\zeta_{\udl\chi_{\udl\kappa}})=(\RU(N_1)\times\RU(N_2),{{\udl\kappa}})
\end{equation}
with $N_1,N_2\geq 0$ and $N=N_1+N_2$; and $\udl\kappa=((-1)^{N-N_1},(-1)^{N-N_2})$.
Write $\psi\in\wt{\Psi}_2(\RU(N))$ in form
\begin{equation}\label{psiu}
\psi=(\tau,b)\boxplus\psi^{N-ab},
\end{equation}
with $\tau\in\CA_\cusp(a)$ being conjugate self-dual and $b>0$; and $(\tau,b)\in\wt{\Psi}_2(\RU(ab))$, and
$\psi^{N-ab}\in\wt{\Psi}_2(\RU(N-ab))$. In order for $(\tau,b)$ to be in $\wt{\Psi}_2(\RU(ab))$, the sign $\eta_{(\tau,b)}$ of
the global simple Arthur parameter $(\tau,b)$ has the property that
$$
\eta_{(\tau,b)}=\kappa_{ab}(-1)^{ab-1}=\eta_\tau(-1)^{b-1}=\kappa(-1)^{N-1},
$$
which implies that $\eta_\tau=\kappa(-1)^{N-b}$ and leads to the following cases:
\begin{itemize}
\item[(1)] {\bf Case (N=2n):}\ If $\tau$ or $(\tau,1)$ is of conjugate orthogonal type, then $\eta_\tau=\kappa_a(-1)^{a-1}=1$ and hence
one must have that $(-1)^b=\kappa$; and if $\tau$ or $(\tau,1)$ is of conjugate symplectic type, then $\eta_\tau=-1$ and hence
one must have that $(-1)^b=-\kappa$.
\item[(2)] {\bf Case (N=2n+1):}\ If $\tau$ or $(\tau,1)$ is of conjugate orthogonal type, then $\eta_\tau=1$ and hence
one must have that $(-1)^b=-\kappa$; and if $\tau$ or $(\tau,1)$ is of conjugate symplectic type, then $\eta_\tau=-1$ and hence
one must have that $(-1)^b=\kappa$.
\end{itemize}

The standard elliptic endoscopy group of $\RU_N$ associated to the global elliptic Arthur parameter $\psi$ as in \eqref{psiu} is
\begin{equation}\label{ed+}
(\RU(ab)\times\RU(N-ab),((-1)^{ab},(-1)^{ab}))
\end{equation}
if $N=2n$, and
\begin{equation}\label{ed-}
(\RU(ab)\times\RU(N-ab),((-1)^{ab+1},(-1)^{ab}))
\end{equation}
if $N=2n+1$.

\section{Fourier Coefficients and Nilpotent Orbits}\label{fcno}

In order to formulate the automorphic integral transforms which produce various types of endoscopy correspondences, it is the key step to define
the right kernel functions for a given setting of the endoscopy correspondences. In principle, the kernel functions are defined
by means of discrete series automorphic representations of simple type in the sense of Arthur and by means of the structures of Fourier coefficients
of automorphic forms defined in terms of nilpotent orbits. The goal of this section is to state and discuss a general conjecture on
the structure of Fourier coefficients of discrete series automorphic representations and their relations to Arthur classification of
the discrete spectrum, which can be viewed as a natural extension of the global Shahidi conjecture on the genericity of tempered $L$-packets.

\subsection{Nilpotent orbits}

Let $G\in\wt{\CE}_\simp(N)$ be an $F$-quasisplit classical group and $\frak{g}$ be the Lie algebra of the algebraic group $G$.
Let $\CN(\frak{g})$ be the set of all nilpotent elements of $\frak{g}$ which is an algebraic conic variety defined over $F$ and is called
the nilcone of $\frak{g}$. Under the adjoint action of $G$, which is defined over $F$, the nilcone $\CN(\frak{g})$ is an algebraic
$G$-variety over $F$. Over algebraic closure $\overline{F}$,  the $\overline{F}$-points $\CN(\frak{g})(\overline{F})$ decomposes into
finitely many adjoint $G(\overline{F})$-orbits $\CO$. It is a well-known theorem (Chapter 5 of \cite{CM93}) that the set of those finitely many nilpotent orbits are in ono-to-one correspondence with the set of partitions of $N_G$ with certain parity constraints (called $G$-partitions of
$N_G$), where $N_G=N+1$ if $G$ is of type $\SO_{2n+1}$,
$N_G=N-1$ if $G$ is of type $\Sp_{2n}$, and $N_G=N$ if $G$ is of type $\SO_{2n}$ or a unitary group $\RU_N=\RU(N)$.

\subsection{Fourier coefficients of automorphic forms}

Over $F$, the $F$-rational points $\CN(\frak{g})(F)$ decomposes into $F$-stable adjoint $G(F)$-orbits $\CO^\st$, which are parameterized by
the corresponding $G$-partitions of $N_G$. For $X\in\CN(\frak{g})(F)$, the theory of the $\frak{sl}_2$-triple (over $F$) asserts
that there is a one-parameter semi-simple subalgebra $\hslash(t)$ in $\frak{g}$ with $t\in F$ and a nilpotent element $Y\in\CN(\frak{g})(F)$
such that the set $\{X,\hslash(1),Y\}$ generates a $F$-Lie subalgebra of $\frak{g}$, which is isomorphic to $\frak{sl}_2(F)$.
This means that
$$
[X,Y]=\hslash(1),\ [\hslash(1),X]=2X,\ [\hslash(1),Y]=-2Y.
$$
Under the adjoint action of $\hslash(1)$, the Lie algebra $\frak{g}(F)$ decomposes into a direct sum of adjoint $\hslash(1)$-eigenspaces
\begin{equation}\label{sl2}
\frak{g}=\frak{g}_{-m}\oplus\cdots\oplus\frak{g}_{-2}\oplus\frak{g}_{-1}\oplus\frak{g}_0\oplus\frak{g}_1\oplus\frak{g}_2\oplus\cdots\frak{g}_m
\end{equation}
for some positive integral $m$,
where $\frak{g}_r:=\{A\in\frak{g}\mid [\hslash(1),A]=rA\}$. Then $X\in\frak{g}_2$ and $Y\in\frak{g}_{-2}$.
Let $L_0$ be the algebraic $F$-subgroup of $G$ such that the Lie algebra of $L_0$ is
$\frak{g}_0$. Then the adjoint $G(F)$-orbits in the stable orbit $\CO^\st_X$ of $X$ are parameterized by the adjoint $L_0(F)$-orbits of maximal
dimension in either $\frak{g}_2(F)$ or $\frak{g}_{-2}(F)$. For $F$-quasisplit classical groups, the classification of
the adjoint $G(F)$-orbits in the stable orbit $\CO^\st_X$ of any $X\in\CN(\frak{g})(F)$ is explicitly given by J.-L. Waldspurger in
\cite[Chapitre I]{W01}.

Let $V_X$ be the unipotent subgroup of $G$ whose Lie algebra is $\oplus_{i\geq 2}\frak{g}_i$. There is a canonical way to define a
nontrivial automorphic character $\psi_X$ on $V_X(\BA)$, which is trivial on $V_X(F)$, from a nontrivial automorphic character $\psi_F$ of
$F\bks\BA$. More precisely, it is given by: for any $v\in V_X(\BA)$,
$$
\psi_X(v):=\psi_F(\tr(Y\cdot\ln(v))).
$$
Let $\pi$ belong to $\CA(G)$. The {\bf Fourier coefficient} of $\varphi_\pi$ in the space of $\pi$ is defined by
$$
\CF^{\psi_X}(\varphi_\pi)(g):=\int_{V_X(F)\bks V_X(\BA)}\varphi_\pi(vg)\psi_X(v)^{-1}dv.
$$
Since $\varphi_\pi$ is automorphic, the nonvanishing of $\CF^{\psi_X}(\varphi_\pi)$ depends only on the $G(F)$-adjoint orbit $\CO_X$ of
$X$.

Define $\frak{n}(\varphi_\pi):=\{X\in\CN(\frak{g})\mid \CF^{\psi_X}(\varphi_\pi)\neq 0\}$, which is stable under the $G(F)$-adjoint action.
$\frak{n}^m(\varphi_\pi)$ is the union of maximal $G(F)$-adjoint orbits in $\frak{n}(\varphi_\pi)$, according to
the the partial ordering of the corresponding stable nilpotent orbits, i.e. that of corresponding partitions.
Define simply that $\frak{n}(\pi)$ to be the union of
$\frak{n}^m(\varphi_\pi)$ for all $\varphi_\pi$ running in the space of $\pi$ and define
$\frak{n}^m(\pi)$ to be the union of maximal $G(F)$-adjoint orbits in $\frak{n}(\pi)$, in the same way as $\frak{n}^m(\varphi)$.

\subsection{Example of $\GL_N$} When $G=\GL_N$, every $G(F)$-orbit is stable and the set of all nilpotent orbits of $\GL_N(F)$
in $\CN(\frak{g})(F)$ are parameterized by partitions of $N$. The well known theorem proved by J. Shalika and I. Piatetski-Shapiro,
independently, says that every $\pi$ in $\CA_{\cusp}(N)$ is generic, that is, there is
an automorphic form $\varphi_\pi$ in the space of $\pi$ having a nonzero Whittaker-Fourier coefficient. Following the definition of
Fourier coefficients above, the Whittaker-Fourier coefficients are attached to the regular nilpotent orbit and hence have the
partition $[N]$ of the integer $N$. Therefore, one has that $\frak{n}^m(\pi)=\{[N]\}$ for any $\pi\in\CA_{\cusp}(N)$.

What happens if $\pi$ belongs to $\CA_2(N)$ in general?
As discussed in Section 2.2, if $\pi$ belongs to $\CA_2(N)$, by a theorem of Moeglin and Waldspurger, $\pi$ is of form $\Delta(\tau,b)$, with
$\tau\in\CA_{\cusp}(a)$ and $N=ab$. The main theorem in \cite{JL12} shows that $\frak{n}^m(\Delta(\tau,b))=\{[a^b]\}$.
See also \cite{G06} for a sketch of different proof.

In general, for $\pi\in\CA(N)$, the result is described in the following conjecture.

\begin{conj}\label{fcgl}
Let $\pi\in\CA(N)$ be of form
$$
\pi_\psi:=\CI_P(\pi_1\otimes\pi_2\otimes\cdots\otimes\pi_r),
$$
with $N=\sum_{i=1}^rN_i$ and with $N_i=a_ib_i$, $\pi_i=\Delta(\tau_i,b_i)\in\CA_2(N_i)$, and $\tau_i\in\CA_{\cusp}(a_i)$ for
$i=1,2,\cdots,r$. Then
$$
\frak{n}^m(\pi)=[a_1^{b_1}]+[a_2^{b_2}]+\cdots+[a_r^{b_r}].
$$
\end{conj}

Note that this conjecture can be considered as a special and more precise version of a more general conjecture discussed by
Ginzburg in \cite{G06} and could be verified by extending the argument of \cite{JL12}.

\subsection{Fourier coefficients and Arthur parameters}
In this section, after an explanation of Conjecture \ref{fcgl} in terms of Arthur parametrization,
a general conjecture will be stated on relation between Fourier coefficients and Arthur parameters for classical groups.

For $\pi\in\CA(N)$ of form,
$
\pi_\psi:=\CI_P(\pi_1\otimes\pi_2\otimes\cdots\otimes\pi_r),
$
with $N=\sum_{i=1}^rN_i$ and with $N_i=a_ib_i$, $\pi_i=\Delta(\tau_i,b_i)\in\CA_2(N_i)$, and $\tau_i\in\CA_{\cusp}(a_i)$ for
$i=1,2,\cdots,r$, as in Conjecture \ref{fcgl}, the global Arthur parameter $\psi$ attached to $\pi$ belongs to $\Psi(N)$ and has the form
$$
\psi=(\tau_1,b_1)\boxplus(\tau_2,b_2)\boxplus\cdots\boxplus(\tau_r,b_r)
$$
with $(\tau_i,b_i)\in\Psi_2(a_ib_i)$ for $i=1,2,\cdots,r$. It is not hard to figure out that as a representation in the
$a_ib_i$-dimensional complex vector space, the irreducible representation $(\tau_i,b_i)$ should correspond to the partition
$[b_i^{a_i}]$ of the integer $a_ib_i$ for $i=1,2,\cdots,r$, and hence the partition of $N$ corresponds to the direct sum of
$(\tau_1,b_1)$, $(\tau_2,b_2), \cdots, (\tau_r,b_r)$ should be $[b_1^{a_1}b_2^{a_2}\cdots b_r^{a_r}]$.

Note that the partition
$
[a_1^{b_1}]+[a_2^{b_2}]+\cdots+[a_r^{b_r}]
$
in Conjecture \ref{fcgl} is the image of the partition
$[b_1^{a_1}b_2^{a_2}\cdots b_r^{a_r}]$ under the Barbasch-Vogan duality mapping $\eta_{\frak{g}^\vee,\frak{g}}$ from the dual
Lie algebra $\frak{g}^\vee$ of $\frak{g}$ to $\frak{g}$ with $\frak{g}=\frak{gl}_N$ (\cite{BV85} and \cite{Ac03}).
This suggests a conjecture for $F$-quasisplit classical groups, which is formulated in terms of global Arthur packets.

\begin{conj}[Upper Bound of Maximal Fourier Coefficients]\label{cubmfc}
For a $F$-quasisplit classical group $G$,  let $\wt{\Pi}_{\psi}(\epsilon_\psi)$ be the automorphic $L^2$-packet attached to a
global Arthur parameter $\psi\in\wt{\Psi}_2(G)$.
Assume that $\underline{p}(\psi)$ is the partition attached to $(\psi,G^\vee(\BC))$. Then for any $\pi\in\wt{\Pi}_{\psi}(\epsilon_\psi)$,
$\frak{n}^m(\pi)$ is bounded above by the partition $\eta_{{\frak{g}^\vee,\frak{g}}}(\underline{p}(\psi))$.
\end{conj}

It remains very interesting to figure out exactly what $\frak{n}^m(\pi)$ are for $\pi\in\wt{\Pi}_{\psi}(\epsilon_\psi)$.

When the global Arthur parameter $\psi\in\wt{\Psi}_2(G)$ is generic, the partition $\underline{p}(\psi)$
attached to $(\psi,G^\vee(\BC))$ must be the trivial partition, i.e. corresponding to the zero nilpotent orbit in $\frak{g}^\vee(\BC)$.
Then the Barbasch-Vogan duality $\eta_{{\frak{g}^\vee,\frak{g}}}(\underline{p}(\psi))$ must be the partition corresponding to the
regular nilpotent orbit in $\frak{g}(F)$ (\cite{BV85} and \cite{Ac03}). Following the existence of the Arthur-Langlands transfer
(\cite{Ar12} and \cite{Mk12}) and the automorphic descent construction of Ginzburg-Rallis-Soudry (\cite{GRS11}), there exists an
irreducible generic cuspidal automorphic representation $\pi$ in the automorphic $L^2$-packet $\wt{\Pi}_{\psi}(\epsilon_\psi)$, such
that $\frak{n}^m(\pi)=\eta_{{\frak{g}^\vee,\frak{g}}}(\underline{p}(\psi))$. By the generalized Ramanujan conjecture for the generic
Arthur parameter $\psi$, this is a global version of the Shahidi conjecture, which claims that any global tempered $L$-packet has
a generic member.

On the other hand, consider a global Arthur parameter of simple type,
$
\psi=(\tau,b)\in\wt{\Psi}_2(G),
$
with $\tau\in\CA_\cusp(a)$ being conjugate self-dual.

If take $G=\SO_{2n+1}$, then $2n=ab$. Hence $\tau$ is of orthogonal type if and only if $b=2l$; and
$\tau$ is of symplectic type if and only if $b=2l+1$. The partition attached to $\psi$ is
$\udl{p}(\psi)=[b^a]$, and the Barbasch-Vogan dual is
$$
\eta_{{\frak{g}^\vee,\frak{g}}}(\underline{p}(\psi))=
\begin{cases}
[(a+1)a^{b-2}(a-1)1]&\text{if}\ b=2l;\\
[(a+1)a^{b-1}]&\text{if}\ b=2l+1.
\end{cases}
$$

If take $G=\SO_{2n}$, then $2n=ab$. Hence $\tau$ is of orthogonal type if and only if $b=2l+1$; and
$\tau$ is of symplectic type if and only if $b=2l$. The partition attached to $\psi$ is
$\udl{p}(\psi)=[b^a]$, and the Barbasch-Vogan dual is
$$
\eta_{{\frak{g}^\vee,\frak{g}}}(\underline{p}(\psi))=
\begin{cases}
[a^b]&\text{if}\ b=2l;\\
[a^{b-1}(a-1)1]&\text{if}\ b=2l+1.
\end{cases}
$$

If take $G=\Sp_{2n}$, then $2n+1=ab$. Hence $\tau$ must be of orthogonal type and $b=2l+1$ and $a=2e+1$. The partition attached to $\psi$ is
$\udl{p}(\psi)=[b^a]$, and the Barbasch-Vogan dual is
$$
\eta_{{\frak{g}^\vee,\frak{g}}}(\underline{p}(\psi))=[a^{b-1}(a-1)].
$$

Following \cite{G03}, \cite{G06}, \cite{G08}, \cite{GRS11}, \cite{JL12} and \cite{JLZ12}, Liu and the author
are able to show in \cite{JL13} that for any simple global Arthur parameter $\psi$ for all $F$-split classical groups,
Conjecture \ref{cubmfc} holds and
there is a member $\pi\in\wt{\Pi}_{\psi}(\epsilon_\psi)$ such that
$\frak{n}^m(\pi)=\eta_{{\frak{g}^\vee,\frak{g}}}(\underline{p}(\psi))$.

These two extreme cases may suggest the technical complication on how to figure out what $\frak{n}^m(\pi)$ are for
$\pi\in\wt{\Pi}_{\psi}(\epsilon_\psi)$ or even the upper bound partition and the lower bound partition for the
automorphic $L^2$-packet $\wt{\Pi}_{\psi}(\epsilon_\psi)$. Some testing cases are considered in \cite{JL13}, which
support that there should always be a member $\pi\in\wt{\Pi}_{\psi}(\epsilon_\psi)$ such that
$\frak{n}^m(\pi)=\eta_{{\frak{g}^\vee,\frak{g}}}(\underline{p}(\psi))$ in general.

\begin{conj}[Maximal Fourier Coefficients]\label{cmfc}
For a $F$-quasisplit classical group $G$,  let $\wt{\Pi}_{\psi}(\epsilon_\psi)$ be the automorphic $L^2$-packet attached to a
global Arthur parameter $\psi\in\wt{\Psi}_2(G)$.
Assume that $\underline{p}(\psi)$ is the partition attached to $(\psi,G^\vee(\BC))$. Then there exists at least one member $\pi\in\wt{\Pi}_{\psi}(\epsilon_\psi)$, such that
$\frak{n}^m(\pi)$ is equal to the partition $\eta_{{\frak{g}^\vee,\frak{g}}}(\underline{p}(\psi))$, i.e. the upper bound of the
maximal Fourier coefficient in Conjecture \ref{cubmfc} can be achieved by members in $\wt{\Pi}_{\psi}(\epsilon_\psi)$
\end{conj}

For a non-generic global Arthur parameter $\psi\in\wt{\Psi}_2(G)$, it is an interesting problem to determine  automorphic
representations occurring as members in the automorphic $L^2$-packet $\wt{\Pi}_{\psi}(\epsilon_\psi)$ by means of the basic
invariants attached to automorphic forms, and this problem is closely related to the topic on CAP automorphic representations,
which was initiated by Piatetski-Shapiro through his work on Saito-Kurokawa lifting (\cite{PS83}).

Along this line of ideas, it is natural to take the structure of Fourier coefficients of irreducible cuspidal automorphic
representations as basic invariants and study the classification problem. The approach should go back to the root of
the theory of automorphic forms. However, the representation-theoretic approach towards this problem goes back to the
notion of the rank of automorphic forms or automorphic representations, introduced by R. Howe (\cite{Hw81} and \cite{Hw82}). The work of J.-S. Li
(\cite{Li95} and \cite{Li97}) give a complete characterization of singular automorphic representations of classical groups in terms of
theta correspondences. From the point of view of the Arthur classification of discrete spectrum, Moeglin published a series of papers
on classifying automorphic representations which are expected to have quadratic unipotent Arthur parameters (\cite{Mg91}
and \cite{Mg94}), which have deep impact to the understanding of automorphic forms related to the work of Kudla and Rallis
on their regularized Siegel-Weil formula (\cite{KR94}) and the further arithmetic applications of Kudla (\cite{Kd97}).
A family of interesting examples of cubic unipotent cuspidal automorphic representations of the exceptional group of type $G_2$ were
constructed through the exceptional theta correspondences by Gan, Gurevich and the author in \cite{GGJ02}, which
produces examples of cuspidal automorphic representations with multiplicity as high as one wishes.

As one will see in late sections, the $(\tau,b)$-theory of automorphic forms and the constructions of endoscopy correspondences
take the structures of Fourier coefficients as a basic invariants of automorphic forms. With the Arthur classification theory in hand,
it is natural to try an understanding of the structure of Fourier coefficients of cuspidal automorphic representations in terms of
the global Arthur parameters. The readers will see this point in the discussion in the rest of this paper and also in
\cite{GRS03}, \cite{G03}, \cite{G06}, \cite{G08}, \cite{Jn11}, \cite{JL12}, and \cite{G12}.  As a side remark, the
structure of Fourier coefficients of certain residual representations made important progress towards the global Gan-Gross-Prasad
conjecture (\cite{GGP12}) through the a series of work of Ginzburg, Rallis and the author (\cite{GJR04}, \cite{GJR05}, and \cite{GJR09}).
Some connections of the Fourier coefficients to arithmetic and even mathematical physics can be found in
\cite{JR97}, \cite{GGS02}, \cite{MS12}, \cite{GMRV10} and \cite{GMV11}.

Finally, it is important to mention that the local theory goes back to the work of Howe on the notion of the rank
for representations (\cite{Hw82}), the work of Bernstein and Zelevinsky on the notion of the derivatives for admissible
representations of $p$-adic $\GL(n)$ (\cite{BZ77}) and the work of A. Aizenbud, D. Gourevitch, and S. Sahi (\cite{AGS11}
for archimedean cases, the work of Moeglin and Waldspurger
on the degenerate Whittaker models of $p$-adic groups (\cite{MW87}) and the work of Gourevitch and Sahi in \cite{GS12} for archimedean cases, and
certain refined structures of the Fourier coefficients was obtained by Moeglin in \cite{Mg96} for $p$-adic tempered representations,
by B. Harris in \cite{H12} for real tempered representations, and by Baiying Liu and the author in \cite{JL13} for cuspidal
automorphic representations of symplectic groups.

\section{Constructions of the Automorphic Kernel Functions}

The {\bf automorphic kernel functions} are constructed in terms of automorphic forms in automorphic $L^2$-packets
$\wt\Pi_{\psi_0}(\epsilon_{\psi_0})$,
with Fourier coefficients associated to relatively small partitions, following the idea of the classical theta correspondence and
the exceptional theta correspondence, which use the minimality of the kernel functions. Section 5.1 recalls the main results of \cite{JLZ12},
which provide at least one concretely constructed member in these automorphic $L^2$-packets
$\wt\Pi_{\psi_0}(\epsilon_{\psi_0})$. Section 5.2 describes explicitly a family of Fourier coefficients which forms a key step of the
constructions. Sections 5.3 and 5.4 give explicit data of the constructions and conjectures on the types of endoscopy correspondences
one expects from those constructions. The last section gives the construction data for unitary group case, the details of which
will be given in a forthcoming work of Lei Zhang and the author (\cite{JZ13}). Note that the constructions discussed here essentially
assume that $d=a-1$ in the construction data. The general case (without this assumption) will be discussed in a forthcoming work of the author
(\cite{Jn13}). Along the way, all known cases to the conjectures will also be briefly discussed.

\subsection{Certain families of residual representations}

Recall from \cite{JLZ12} certain families of the residual representations of $F$-quasisplit classical groups, which may be
considered as candidates for the constructions of {\bf automorphic kernel functions}.

Following the notation in \cite{MW95}, a standard Borel subgroup $P_0=M_0N_0$ of $G_n$ is fixed and realized in the upper-triangular
matrices. Let $T_0$ be the maximal split torus of the center of $M_0$, which defines the root system $R(T_0,G_n)$
with the positive roots $R^+(T_0,G_n)$ and the set $\Delta_0$ of simple roots corresponding to $P_0$.
Let $P=MN$ be a standard parabolic subgroup of $G_n$ (containing $P_0$) and
$T_{M}$ be the maximal split torus in the center of $M$. The set of restricted roots is denoted by $R(T_M,G_n)$. Define
$R^+(T_M,G_n)$ and $\Delta_M$, accordingly. Furthermore, define $X_{M}=X^{G_n}_{M}$ to be the group of all continuous homomorphisms from
$M(\BA)$ into $\BC^\times$, which are trivial on $M(\BA)^1$. Then following Page 6 of \cite{MW95} for the
explicit realization of $X_{M}$, define the real part of $X_{M}$, which is denoted by $\Re X_{M}$.

In particular, for a partition $n=r+n_0$, take the standard maximal parabolic subgroup $P_{r}=M_{r}N_{r}$ of $G_{n}$,
whose Levi subgroup $M_{r}$ is isomorphic to $\RR_{F'/F}(r)\times G_{n_0}$.
For any $g\in \RR_{F'/F}(r)$, define $\hat{g}=w_{r}g^{t}w_{r}$ or $w_{r}c(g)^{t}w_{r}$ in the unitary groups cases,
where $w_{r}$ be the anti-diagonal symmetric matrix defined inductively by
$\begin{pmatrix}0&1\\ w_{r-1}&0\end{pmatrix}$ and $c\in\Gamma_{E/F}=\{1,c\}$.
Then any $g\in M_{r}$ can be expressed as $\diag\{t,h,\hat{t}^{-1}\}$ with $t\in \RR_{F'/F}(r)$ and $h\in G_{n_0}$.
Since $P_{r}$ is maximal, the space of characters $X^{G_{n}}_{M_{r}}$ is one-dimensional. Using the normalization in \cite{Sh10}, it is
identified with $\BC$ by $s\mapsto\lam_s$.

Let $\sigma$ be an irreducible generic cuspidal automorphic representation of $G_{n_0}(\BA)$. Write $r=ab$.
Let $\phi$ be an automorphic form in the space $A(N_{ab}(\BA)M_{ab}(F)\bks G_{n}(\BA))_{\Delta(\tau,b)\otimes\sigma}$. Following \cite{L76} and \cite{MW95},
an Eisenstein series is defined by
$$
E^n_{ab}(\phi_{\Delta(\tau,b)\otimes\sigma},s)=E(\phi_{\Delta\otimes\sigma},s)=\sum_{\gamma\in P_{ab}(F)\bks G_{n}(F)}\lam_s\phi(\gamma g).
$$
It converges absolutely for the real part of $s$ large and has meromorphic continuation to the whole complex plane $\BC$.
When $n_0=0$, $\sigma$ is assumed to disappear.

Assume that $n_0>0$ here. The case when $n_0=0$ is similar and is referred to \cite[Theorem 5.2]{JLZ12}.
To determine the location of possible poles (at $\Re(s)\geq 0$)
of this family of residual Eisenstein series, or
more precisely the normalized Eisenstein series, and basic properties of the corresponding residual representations,
one takes the expected normalizing factor $\beta_{b,\tau,\sigma}(s)$ to be of the Langlands-Shahidi type, which is
given by a product of relevant automorphic $L$-functions:
\begin{equation}\label{nf}
\beta_{b,\tau,\sigma}(s)
:=L(s+\frac{b+1}{2},\tau\times\sigma)
\prod_{i=1}^{\lceil \frac{b}{2} \rceil} L(e_{b,i}(s)+1,\tau,\rho) \prod_{i=1}^{\lfloor \frac{b}{2} \rfloor} L(e_{b,i}(s),\tau,\rho^{-}),
\end{equation}
where $e_{b,i}(s):=2s+b+1-2i$, and $\rho$ and $\rho^-$ are defined as follows:
\begin{equation}\label{rho}
\rho:=\begin{cases}
\Asai^+ & \text{ if }  G_{n}=\RU_{2n}\\
\Asai^- & \text{ if }  G_{n}=\RU_{2n+1}\\
\sym^{2} & \text{ if } G_{n}=\SO_{2n+1}\\
\wedge^{2} & \text{ if } G_{n}=\Sp_{2n} \text{ or } \SO_{2n},
\end{cases}
\end{equation}
and
\begin{equation}\label{rho-}
\rho^{-}:=\begin{cases}
\Asai^- & \text{ if } G_{n}=\RU_{2n}\\
\Asai^+ & \text{ if } G_n=\RU_{2n+1}\\
\wedge^{2} & \text{ if } G_{n}=\SO_{2n+1}\\
\sym^{2} & \text{ if } G_{n}=\Sp_{2n} \text{ or } \SO_{2n}.
\end{cases}
\end{equation}
For unitary groups, $\Asai^+$ denotes the Asai representation of the $L$-group of $R_{E/F}(a)$ and $\Asai^-$ denotes
the Asai representation of the $L$-group of $R_{E/F}(a)$ twisted by $\omega_{E/F}$,
the character associated to the quadratic extension $E/F$ via the class field theory. For symplectic or orthogonal groups,
$\sym^2$ and $\wedge^2$ denote the symmetric and
exterior second powers of the standard representation of $\GL_a(\BC)$, respectively. In addition, one has the
following identities (Remark (3), Page 21, \cite{GRS11}):
 $$
 L(s,\tau\times\tau^c)=L(s,\tau,\rho)L(s,\tau,\rho^{-}),
 $$
where $\tau^* = \tau$, if $F'=F$; and $\tau^* = \tau^c$, if $F'=E$.
The involution $c$ is the nontrivial element in the Galois group $\Gamma_{E/F}$.

The function $\beta_{b,\tau,\sigma}(s)$ is used to normalize the Eisenstein series by
\begin{equation}\label{eq:iNE}
E^{n,*}_{ab}(\phi_{\Delta(\tau,b)\otimes \sig},s):=\beta_{b,\tau,\sigma}(s)E^n_{ab}(\phi_{\Delta(\tau,b)\otimes \sig},s).
\end{equation}

In order to determine the location of the poles of $E^{*}(\phi_{\Delta(\tau,b)\otimes \sig},s)$,
one needs to consider the following four cases:
\begin{enumerate}
\item\label{icase1} $L(s,\tau,\rho)$ has a pole at $s=1$, and $L(\frac{1}{2},\tau\times\sig)\neq 0$;
\item\label{icase2} $L(s,\tau,\rho)$ has a pole at $s=1$, and $L(\frac{1}{2},\tau\times\sig)=0$;
\item\label{icase3} $L(s,\tau,\rho^{-})$ has a pole at $s=1$, and $L(s,\tau\times\sig)$ has a pole at $s=1$;
\item\label{icase4} $L(s,\tau,\rho^{-})$ has a pole at $s=1$, and $L(s,\tau\times\sig)$ is holomorphic at $s=1$.
\end{enumerate}

The sets of possible poles according to the four cases are:
$$
X^+_{b,\tau,\sigma}:=\begin{cases}
\{\hat{0},\dots,\frac{b-2}{2},\frac{b}{2}\}, & \text{if Case (1)};\\
\{\hat{0},\dots,\frac{b-4}{2},\frac{b-2}{2}\}, & \text{ if Case (2)};\\
\{\hat{0},\dots,\frac{b-1}{2},\frac{b+1}{2}\}, & \text{ if Case (3)}; \\
\{\hat{0},\dots,\frac{b-3}{2},\frac{b-1}{2}\}, & \text{ if Case (4)}.
\end{cases}
$$
Note that $0$ is omitted in the set $X^+_{b,\tau,\sigma}$. The main results of \cite{JLZ12} can be summarized as follows.

\begin{thm}\label{jlzm}
Assume that $n_0>0$.
Let $\sig$ be an irreducible generic cuspidal automorphic representation of $G_{n_0}(\BA)$, and
$\tau$ be an irreducible unitary self-dual cuspidal automorphic representation of $\RR_{F'/F}(a)(\BA)$.
Then the following hold.
\begin{enumerate}
\item The normalized Eisenstein series $E^{n,*}_{ab}(\phi_{\Delta(\tau,b)\otimes \sig},s)$ is holomorphic for $\Re(s)\geq 0$ except
at $s=s_0\in X^+_{b,\tau,\sigma}$ where it has possibly at most simple poles.
\item For $s_0\in(0,\frac{b+1}{2}]$, assume that
the normalized Eisenstein series $E^{n,*}_{ab}(\phi_{\Delta(\tau,b)\otimes \sig},s)$
has a simple pole at $s=s_{0}$. Then the residue of $E^{n,*}_{ab}(\phi_{\Delta(\tau,b)\otimes \sig},s)$
at $s_{0}$ is square-integrable except $s_{0}=\frac{b-1}{2}$ in {\rm Case (3)}.
\item The global Arthur parameter for each of those residual representations is given in \cite[Section 6.2]{JLZ12}.
\end{enumerate}
\end{thm}

Note that the $n_0=0$ case is given in \cite[Theorem 5.2]{JLZ12}.
In \cite{JLZ12}, the theorem was proved for the cases when $G_n$ is not unitary. However, it is expected that the same argument
works for unitary groups as long as some technical results from \cite{Ar12} are also valid for unitary groups (\cite{Mk12}). .

The residual representations at the right end points of $X^+_{b,\tau,\sigma}$ or of $X^+_{b,\tau}$ if $n_0=0$
are the candidates currently used
to produce the {\bf automorphic kernel functions} in the explicit constructions of endoscopy correspondences for classical groups.
It is the point of \cite{JLZ12} that the other residual representations may also be able to produce {\bf automorphic kernel
functions} for other kinds of Langlands functorial transfers. One example of this kind is referred to \cite{GJS11}.

\subsection{Certain families of Fourier coefficients}
Fourier coefficients of automorphic forms on classical groups considered here are attached to a particular family of partitions
of type $[d^c1^*]$. The pair of $d$ and $c$ satisfies a parity condition according to the type of the classical group $G_m$.
More precisely one has the following cases.
\begin{enumerate}
\item When $d=2k-1$, $c=2f$ if $G_m$ is a symplectic group; and $c=2f$ or $c=2f+1$ if $G_m$ is an orthogonal group or unitary group.
\item When $d=2k$, $c=2f$ if $G_m$ is an orthogonal group; and $c=2f$ or $c=2f+1$ if $G_m$ is a symplectic group or unitary group.
\end{enumerate}
As in \cite[Chapitre I]{W01} the set of the representatives of the adjoint $G_n(F)$-orbits in the given stable orbit $\CO^\st_{[d^c1^*]}$
associated to the partition $[d^c1^*]$can be explicitly given.

{\bf Symplectic Group Case.}\
Assume that $G_m$ is the symplectic group $\Sp_{2m}$. If $d=2k-1$, then the integer $c=2f$ must be even,
By \cite[Chapitre I]{W01}, there is only one
$F$-rational $G_m(F)$-adjoint orbit in the given stable orbit $\CO^\st_{[d^c1^{2m-cd}]}$
associated to the partition $[d^c1^{2m-cd}]$, which corresponds to the pair $([d^c1^{2m-cd}], (\cdot,\cdot))$. The notation
$\udl{q}:=(\cdot,\cdot)$ means that there is no quadratic forms involved in the parametrization of
$F$-rational $G_m(F)$-adjoint orbits in the given stable orbit $\CO^\st_{[d^c1^{2m-cd}]}$. To indicate the dependence on
$\udl{q}$, $V_{X_{\udl{q}}}$ is used for $V_X$, and $\psi_{X_{\udl{q}}}$ for $\psi_X$, respectively.
The stabilizer in $L_0$ of the character $\psi_{X_{\udl{q}}}$ is isomorphic to
$$
\Sp_{c}\times\Sp_{2m-cd}.
$$
Hence, for any automorphic form $\varphi$ on $G_m(\BA)$, the $\psi_{X_{\udl{q}}}$-Fourier coefficient
$\CF^{\psi_{X_{\udl{q}}}}(\varphi)$ of $\varphi$, as defined in Section 4.2,
is automorphic on
$$
\Sp_{c}(\BA)\times\Sp_{2m-cd}(\BA).
$$

\begin{rmk} For a general automorphic form $\varphi$ on $G_m(\BA)$, one does not expect that the $\psi_{X_{\udl{q}}}$-Fourier coefficients
$\CF^{\psi_{X_{\udl{q}}}}(\varphi)(h,g)$ produces a reasonable (functorial) relation between automorphic forms on $\Sp_c(\BA)$ and automorphic forms
on $\Sp_{2m-cd}(\BA)$, by simply assuming that the following integral
$$
\int_{\Sp_c(F)\bks\Sp_c(\BA)}\int_{\Sp_{2m-cd}(F)\bks\Sp_{2m-cd}(\BA)}
\CF^{\psi_{X_{\udl{q}}}}(\varphi)(h,g)\varphi_\sigma(h)\ovl{\varphi_\pi(g)}dhdg
$$
converges and is non-zero for some choice of data, where $\sigma\in\CA_\cusp(\Sp_{c})$ and $\pi\in\CA_\cusp(\Sp_{2m-cd})$.

However, if one follows the idea of the classical theta correspondence using the minimality of the classical theta functions
or the Weil representations, the $\psi_{X_{\udl{q}}}$-Fourier coefficients $\CF^{\psi_{X_{\udl{q}}}}(\varphi)(h,g)$ do produce, via the
integral above,
reasonable functorial relations when one puts enough restrictions on the automorphic forms $\varphi$. This will be discussed in more detail
in Sections 5.3 and 5.4.
\end{rmk}

If $d=2k$, then $c=2f$ or $2f+1$ is any positive integer. Following \cite[Chapitre I]{W01}, the set of the $F$-rational
$G_m(F)$-adjoint orbits in the stable orbit $\CO^\st_{[d^{c}1^{(2m-cd)}]}$ are parameterized by the pairs
$([d^{c}1^{(2m-cd)}],(q_{d},\cdot))$, where $q_{d}$ runs through the set of all $F$-equivalence classes of the $c$-dimensional
non-degenerate quadratic forms over $F$. In this case, put $\udl{q}:=(q_d,\cdot)$.
The stabilizer in $L_0$ of the character $\psi_{X_{\udl{q}}}$ is isomorphic to
$$
\SO_c(q_{d})\times\Sp_{2m-cd}.
$$
Note that in this case, the space $\frak{g}_1$ in (4.1) is nonzero, which leads the consideration of the {\bf Fourier-Jacobi} coefficients in
the construction of {\bf automorphic kernel functions}, instead of using the $\psi_{X_{\udl{q}}}$-Fourier coefficients. This will
be discussed in detail below.

{\bf Even Special Orthogonal Group Case.}\
Assume that $G_m$ is an $F$-quasisplit even special orthogonal group $\SO(q_V)$
defined by the $2m$-dimensional non-degenerate quadratic form $q_V$.

If $d=2k$ is an even positive integer, then
$c=2f$ must be an even positive integer. Following \cite[Chapitre I]{W01}, there is only one $F$-rational $G_m(F)$-adjoint
orbit in the given stable orbit $\CO^\st_{[d^c1^{(2m-cd)}]}$
associated to the partition $[d^c1^{(2m-cd)}]$, corresponding to the pair $([d^c1^{(2m-cd)}],(\cdot,q_1))$, where $q_1$ is
the $(2m-cd)$-dimensional non-degenerate quadratic form over $F$ sharing the $F$-anisotropic kernel with $q_V$. Put
$\udl{q}:=(\cdot,q_1)$. Then the stabilizer in $L_0$ of the character $\psi_{X_{\udl{q}}}$ is isomorphic to
$$
\Sp_{c}\times\SO_{2m-cd}(q_1).
$$
Note that in this case, the space $\frak{g}_1$ in (4.1) is nonzero, which leads the consideration of the {\bf Fourier-Jacobi} coefficients in
the construction of {\bf automorphic kernel functions}, instead of using the $\psi_{X_{\udl{q}}}$-Fourier coefficients. This will
be discussed in detail below.

If $d=2k-1$ is an odd positive integer, then $c=2f$ or $2f+1$ is any positive integer. According to \cite[Chapitre I]{W01},
the set of the $F$-rational $G_m(F)$-adjoint orbits in
the stable orbit $\CO^\st_{[d^{c}1^{(2m-cd)}]}$ are parameterized by the pairs $([d^{c}1^{2m-cd}],(q_{d},q_1))$, where
$q_{d}$ runs through the set of all $F$-equivalence classes of the $c$-dimensional non-degenerate quadratic forms over $F$ and
$q_1$ runs through the set of all $F$-equivalence classes of the $(2m-cd)$-dimensional non-degenerate quadratic forms over $F$ such that
the quadratic forms $q_d\oplus q_1$ have the same $F$-anisotropic kernel as $q_V$. Put $\udl{q}:=(q_{d},q_1)$.
The stabilizer in $L_0$ of the character $\psi_{X_{\udl{q}}}$ is isomorphic to
$$
\SO_c(q_d)\times\SO_{2m-cd}(q_1).
$$
For any automorphic form $\varphi$ on $G_m(\BA)$, the $\psi_{X_{\udl{q}}}$-Fourier coefficient
$\CF^{\psi_{X_{\udl{q}}}}(\varphi)$ of $\varphi$, as defined in Section 4.2,
is automorphic on
$$
\SO_c(q_d,\BA)\times\SO_{2m-cd}(q_1,\BA).
$$

{\bf Odd Special Orthogonal Group Case.}\
Assume that $G_m$ is an $F$-split odd special orthogonal group $\SO_{2m+1}(q_V)$.

If $d=2k$ is an even positive integer, then
$c=2f$ must be an even positive integer. Following \cite[Chapitre I]{W01}, there is only one $F$-rational
$G_m(F)$-adjoint orbit in the given stable orbit $\CO^\st_{[d^c1^{(2m+1-cd)}]}$
associated to the partition $[d^c1^{(2m+1-cd)}]$, corresponding to the pair $([d^c1^{(2m+1-cd)}],(\cdot,q_1))$, where
$q_1$ is the $(2m+1-cd)$-dimensional $F$-split quadratic form over $F$. In this case, put $\udl{q}:=(\cdot,q_1)$.
The stabilizer in $L_0$ of the character $\psi_{X_{\udl{q}}}$ is isomorphic to
$$
\Sp_{c}\times\SO_{2m+1-cd}(q_1).
$$
Note that in this case, the space $\frak{g}_1$ in (4.1) is nonzero, which leads the consideration of the {\bf Fourier-Jacobi} coefficients in
the construction of {\bf automorphic kernel functions}, instead of using the $\psi_{X_{\udl{q}}}$-Fourier coefficients. This will
be discussed in detail below.

If $d=2k-1$ is an odd positive integer, then $c=2f$ or $2f+1$ is any positive integer. According to \cite[Chapitre I]{W01},
the set of the $F$-rational $G_m(F)$-adjoint orbits in
the stable orbit $\CO^\st_{[d^{c}1^{2m+1-cd}]}$ are parameterized by the pairs $([d^{c}1^{2m+1-cd}],(q_{d},q_1))$, where
$q_{d}$ runs through the set of all $F$-equivalence classes of the $c$-dimensional non-degenerate quadratic forms over $F$ and
$q_1$ runs through the set of all $F$-equivalence classes of the $(2m+1-cd)$-dimensional non-degenerate quadratic forms over $F$ such that
the quadratic forms $q_d\oplus q_1$ are $F$-split. Put $\udl{q}:=(q_{d},q_1)$.
The stabilizer in $L_0$ of the character $\psi_{X_{\udl{q}}}$ is isomorphic to
$$
\SO_c(q_d)\times\SO_{2m+1-cd}(q_1).
$$
For any automorphic form $\varphi$ on $G_m(\BA)$, the $\psi_{X_{\udl{q}}}$-Fourier coefficient
$\CF^{\psi_{X_{\udl{q}}}}(\varphi)$ of $\varphi$, as defined in Section 4.2,
is automorphic on
$$
\SO_c(q_d,\BA)\times\SO_{2m+1-cd}(q_1,\BA).
$$

{\bf Unitary Group Case.}\
Assume that $G_m$ is a $F$-quasisplit unitary group $\RU(m_V)=\RU_{E/F}(m_V)=\RU_{m_V}(q_V)$ defined
by an $m_V$-dimensional non-degenerate hermitian form $q_V$ with $m=[\frac{m_V}{2}]$ and assume that $\kappa$ is the sign of
endoscopy data $(\RU(m_V),\xi_{\chi_\kappa})$ of $\RR_{E/F}(m_V)$. By \cite[Chapitre I]{W01}, the set of the adjoint $G_m(F)$-orbits in
the stable orbit $\CO^\st_{[d^{c}1^{(m_V-cd)}]}$ are parameterized by the pairs $([d^{c}1^{(m_V-cd)}],(q_{d},q_1))$. Here
$\udl{q}:=(q_{d},q_1)$ can be more specified as follows:
\begin{enumerate}
\item When $d=2k-1$,
$q_{d}$ runs through the set of all $F$-equivalence classes of the $c$-dimensional non-degenerate hermitian forms over $F$ and
$q_1$ runs through the set of all $F$-equivalence classes of the $(m_V-cd)$-dimensional non-degenerate hermitian forms over $F$ such that
the hermitian forms $q_d\oplus q_1$ have the same $F$-anisotropic kernel as $q_V$; and
\item when $d=2k$,
$q_{d}$ runs through the set of all $F$-equivalence classes of the $c$-dimensional non-degenerate hermitian forms over $F$ and
$q_1$ is the $(m_V-cd)$-dimensional non-degenerate hermitian form over $F$ such that $q_1$ shares the same $F$-anisotropic kernel as $q_V$.
\end{enumerate}

The stabilizer in $L_0$ of the character $\psi_{X_{\udl{q}}}$ is isomorphic to
$$
\RU_c(q_d)\times\RU_{m_V-cd}(q_1).
$$
For any automorphic form $\varphi$ on $G_m(\BA)$, the $\psi_{X_{\udl{q}}}$-Fourier coefficient
$\CF^{\psi_{X_{\udl{q}}}}(\varphi)$ of $\varphi$, as defined in Section 4.2,
is automorphic on
$$
\RU_c(q_d,\BA)\times\RU_{m_V-cd}(q_1,\BA).
$$
Note that when $d=2k$, the space $\frak{g}_1$ in (4.1) is nonzero, which leads the consideration of the {\bf Fourier-Jacobi} coefficients in
the construction of {\bf automorphic kernel functions}, instead of using the $\psi_{X_{\udl{q}}}$-Fourier coefficients. This will
be discussed in more detail below.

Note that the notion of Fourier coefficients and related results hold for the metaplectic double cover $\tilsp_{2n}(\BA)$ of $\Sp_{2n}(\BA)$.

\subsection{Automorphic kernel functions: Case (d=2k-1)}\label{ec(1)}
From the constructions of the family of Fourier coefficients, this section starts the discussion on the possible constructions of
{\bf automorphic kernel functions}, which produce certain types of {\bf endoscopy correspondences} via integral transforms
as in Conjecture \ref{mcec}.

Recall that for any global Arthur parameter $\psi$, the automorphic $L^2$-packet $\wt{\Pi}_\psi(\epsilon_\psi)$ of the global Arthur
packet $\wt{\Pi}_\psi$ consists of all elements $\pi=\otimes_v\pi_v\in\wt{\Pi}_\psi$ whose characters are equal to $\epsilon_\psi$.
The $[d^c1^*]$-type Fourier coefficients are applied to the automorphic representations of $G_m(\BA)$ belonging to the automorphic
$L^2$-packet $\wt{\Pi}_{\psi_0}(\epsilon_{\psi_0})$ for a relatively simple global Arthur parameter $\psi_0$, which will be specified below in each case.
Of course, it is important to know that such an automorphic $L^2$-packet $\wt{\Pi}_{\psi_0}(\epsilon_{\psi_0})$ is not empty, which
essentially follows from the work of \cite{JLZ12} and the discussion in Section 5.1. Certain refined structures on the cuspidal members in $\wt{\Pi}_{\psi_0}(\epsilon_{\psi_0})$ have been studied
first in \cite{GJS12} and more recently in \cite{Liu13} and \cite{JL13} for some classical groups.
It is expected that those work can be extended to all
classical groups. The structures of Fourier coefficients (Conjectures \ref{cubmfc} and \ref{cmfc}) for such an automorphic $L^2$-packet
$\wt{\Pi}_{\psi_0}(\epsilon_{\psi_0})$ can be proved as discussed in Section 4.4.

Assume that $\tau\in\CA_\cusp(a)$ is self-dual with $a=2e$ even, and assume that
$$
d=a-1=2e-1
$$
when consider the $[d^c1^*]$-type Fourier coefficients. The other cases of the pairs $(a,d)$ are considered in a forthcoming
work of the author (\cite{Jn13}).

{\bf Symplectic Group Case.}\
Assume that $G_m$ is the symplectic group $\Sp_{2m}$.

Write $2m=a(b+c)$ with $a=2e\geq 2$ being even and $b\geq 0,c\geq 1$ being integers. For a self-dual $\tau\in\CA_\cusp(a)$,
if $\tau$ is of symplectic type, take $b+c$ to be even and hence $b$ and $c$ are in the same parity;
and if $\tau$ is of orthogonal type, take $b+c$ to be odd and hence $b$ and $c$ are in the
different parity. Consider a global Arthur parameter
$$
\psi_0:=(\tau,b+c)\boxplus(1,1)
$$
of $\Sp_{2m}$, which is the lift of
the global Arthur parameter $(\tau,b+c)$ of simple type for $\SO_{2m}$, via the natural embedding of the corresponding dual groups.

For the partition $[d^c1^{ab+c}]$ for $\Sp_{2m}$ with $d=a-1$ odd, one must have that $c=2f$ is even. Hence if $\tau$ is of
symplectic type,
then $b=2l$ is even; and if $\tau$ is of orthogonal type, then $b=2l+1$ is odd.
It is not difficult to check that the automorphic $L^2$-packet $\wt{\Pi}_{\psi_0}(\epsilon_{\psi_0})$ of $G_m$ is not empty.
In fact, it follows directly from Theorem 5.1 when $\tau$ is of symplectic type and $b+c$ is even. When $\tau$ is of orthogonal type and
$b+c$ is odd, one may also construct a nonzero residual representation in the automorphic $L^2$-packet $\wt{\Pi}_{\psi_0}(\epsilon_{\psi_0})$.
To do so, one needs to construct an irreducible cuspidal automorphic representation $\epsilon$ of $\Sp_{2e}(\BA)$ from $\tau$ which
has the global Arthur parameter
$$
\psi_\epsilon=(\tau,1)\boxplus(1,1).
$$
This is done by the composition of the automorphic descent from $\tau$ of $\GL_{2e}(\BA)$ to $\SO_{2e}(\BA)$ with the theta correspondence
from $\SO_{2e}(\BA)$ to $\Sp_{2e}(\BA)$.

Take a member $\Theta$ in the automorphic $L^2$-packet $\wt{\Pi}_{\psi_0}(\epsilon_{\psi_0})$ of $G_m$.
The $\psi_{X_{\udl{q}}}$-Fourier coefficients $\CF^{\psi_{X_{\udl{q}}}}(\varphi_\Theta)$ of $\varphi_\Theta\in\Theta$,
as discussed in Section 5.2, are automorphic as functions restricted to the stabilizer
$$
\Sp_c(\BA)\times\Sp_{ab+c}(\BA).
$$
Define the {\bf automorphic kernel functions} in this case to be
\begin{equation}\label{eksp}
\CK_{\varphi_\Theta}^{\psi_0}(h,g):=\CF^{\psi_{X_{\udl{q}}}}(\varphi_\Theta)(h,g)
\end{equation}
for $(h,g)\in\Sp_c(\BA)\times\Sp_{ab+c}(\BA)$. This family of the automorphic kernel functions
$\CK_{\varphi_\Theta}^{\psi_0}(h,g)$ should produce the endoscopy correspondence
$$
\SO_{ab}(\BA)\times\Sp_c(\BA)\rightarrow\Sp_{ab+c}(\BA)
$$
for {\bf Case (a=2e)} of symplectic groups. Note that $\Sp_{ab+c}$ is the symplectic group $\Sp_{2n}$ with
$2n=ab+c$ and $c=2f$ as discussed in Section 3.2; and $\SO_{ab}$ is an $F$-quasisplit even special orthogonal group.
Note that when $b=0$, this is just the identity transfer from $\Sp_c(\BA)$ to itself.

More precisely, take $\psi_1:=(\tau,b)$. It is clear that $\psi_1$ is a simple global Arthur parameter of $\SO_{ab}$. Now one can easily
see that the global Arthur parameter $\psi_0$ is built from the global Arthur parameter $\psi_1$ of $G_0=\SO_{ab}$ and the structures of
the groups $H=\Sp_c$ and $G=\Sp_{ab+c}$.
This construction of the automorphic kernel functions and their related endoscopy correspondence can be described by the following diagram:
$$
\begin{matrix}
&&&\Sp_{a(b+c)}&\\
&&&&\\
&&&\Theta\in\wt\Pi_{\psi_0}(\epsilon_{\psi_0})&\\
&&&&\\
&&\nearrow&\downarrow&\\
&&&&\\
&\GL_a,\tau&&\CK_{\varphi_\Theta}^{\psi_0}(h,g)&\\
&&&&\\
\swarrow&&&\downarrow&\\
&&&&\\
\SO_{ab}&\times&\Sp_{c}&\longleftrightarrow&\Sp_{ab+c}\\
&&&&\\
\wt\Pi_{(\tau,b)}(\epsilon_{(\tau,b)})&&\wt\Pi_{\psi_2}(\epsilon_{\psi_2})&&\wt\Pi_{(\tau,b)\boxplus\psi_2}(\epsilon_{(\tau,b)\boxplus\psi_2})
\\
\\
\end{matrix}
$$
where $a=2e$, $c=2f$, and $2n=ab+c$. When $\tau$ is of symplectic type., $b=2l$ is even; and when $\tau$ is of orthogonal type, $b=2l+1$ is odd.

The main conjecture (Conjecture \ref{mcec}) specializes to the following conjecture
for the current case.

\begin{conj}[$\Sp_{2n}$: {\bf Case (a=2e)}]\label{ecsp1}
For integers $a=2e\geq 2$, $b\geq 0$, $c=2f\geq 2$, and $d=a-1$, let $\tau\in\CA_\cusp(a)$ be self-dual and
$\psi_1=(\tau,b)$ be a simple global Arthur parameter of an $F$-quasisplit $\SO_{ab}$; and
let $\psi_2$ and $\psi:=\psi_1\boxplus\psi_2$ be global Arthur parameters of $\Sp_c$ and $\Sp_{ab+c}$, respectively.
For $\sigma\in\CA_2(\Sp_c)$ and $\pi\in\CA_2(\Sp_{ab+c})$,
if there exists an
automorphic member $\Theta\in\wt{\Pi}_{\psi_0}(\epsilon_{\psi_0})$, such that the following integral
$$
\int_{\Sp_c(F)\bks \Sp_c(\BA)}\int_{\Sp_{ab+c}(F)\bks \Sp_{ab+c}(\BA)}
\CK_{\varphi_\Theta}^{\psi_0}(h,g)\varphi_\sigma(h)\ovl{\varphi_\pi(g)}dhdg
$$
is nonzero for some choice of $\varphi_\Theta\in\Theta$, $\varphi_\sigma\in\sigma$, and $\varphi_\pi\in\pi$, assuming the
convergence of the integral, then $\sigma\in\wt{\Pi}_{\psi_2}(\epsilon_{\psi_2})$ if and only if
$\pi\in\wt{\Pi}_{\psi}(\epsilon_{\psi})$.
\end{conj}

Note that when $\tau$ is of orthogonal type and $b=1$ and when the global Arthur parameter $\psi_2$ is generic, the construction
of the endoscopy transfer for irreducible generic cuspidal automorphic representations was outlined in \cite{G08} and see also
\cite{G12}.

{\bf Example 1.}\ A more explicit explanation of this construction is given for a phototype example.
Take a standard parabolic subgroup $P_{c^{e-1}}$ of $\Sp_{a(b+c)}$ whose Levi decomposition is given by
$$
P_{c^{e-1}}=M_{c^{e-1}}V_{c^{e-1}}=(\GL_{c}^{\times(e-1)}\times\Sp_{ab+2c})V_{c^{e-1}}.
$$
The unipotent radical $V_{c^{e-1}}$, which is the same as $V_{X_{\udl{q}}}$, have the structure that
$$
V_{c^{e-1}}/[V_{c^{e-1}},V_{c^{e-1}}]
\cong
\Mat_{c}^{\oplus(e-2)}\oplus\Mat_{c,n}\oplus\Mat_{c}\oplus\Mat_{c,n},
$$
where $\Mat_{c}$ is the space of all $c\times c$-matrices and $\Mat_{c,n}$ is the space of all $c\times n$-matrices with $n=eb+f$.
The projection of elements $v\in V_{c^{e-1}}$ is given by
$$
(X_1,\cdots,X_{e-2},Y_1,Y_2,Y_3).
$$
For a nontrivial additive character $\psi_F$ of $F\bks\BA$, the corresponding character $\psi_{X_{\udl{q}}}$ on $V_{c^{e-1}}(\BA)$
is defined by
\begin{equation}\label{psi(2f)(e-1)}
\psi_{X_{\udl{q}}}(v)
:=
\psi_F(\tr(X_1+\cdots+X_{e-2}+Y_2)).
\end{equation}
Then the stabilizer of the character $\psi_{X_{\udl{q}}}$ in the Levi subgroup $M_{c^{e-1}}$ is
\begin{equation}\label{stab(1)}
\Sp_{c}\times\Sp_{ab+c}\mapsto (\GL_{c}^{\times(e-1)}\times\Sp_{ab+2c})
\end{equation}
which is given by the embedding
$$
(h,g)\mapsto(h^{\Delta(e-1)},\begin{pmatrix}g_1& &g_2\\ &h&\\ g_3& &g_4\end{pmatrix})
$$
where $h\mapsto h^{\Delta(e-1)}$ is the diagonal embedding from $\Sp_{c}$ into $\GL_{c}^{\times(e-1)}$, and
$g=\begin{pmatrix}g_1&g_2\\ g_3&g_4\end{pmatrix}\in\Sp_{ab+c}$.

{\bf Even Special Orthogonal Group Case.}\
Assume that $G_m$ is an $F$-quasisplit even special orthogonal group $\SO_{2m}(q_V)$
defined by the $2m$-dimensional non-degenerate quadratic form $q_V$.

Write $2m=a(b+c)$ with $a=2e\geq 2$ being even and $b,c\geq 1$ being any positive integers. For a self-dual $\tau\in\CA_\cusp(a)$,
take a simple global Arthur parameter
$$
\psi_0:=(\tau,b+c)
$$
for $\SO_{2m}(q_V)$. Hence if $\tau$ is of symplectic type, then $b+c$ is even and
hence $b$ and $c$ are in the same parity; and if $\tau$ is of orthogonal type, then $b+c$ is odd and hence $b$ and $c$ are in the
different parity.

Assume that $c=2f+1$ is odd. If $\tau$ is of symplectic type, then $b=2l+1$ is odd; and if $\tau$ is of orthogonal type, then $b=2l$ is even.
Consider the partition $[d^c1^{(ab+c)}]$ with $d=a-1$. Each $F$-rational adjoint $G_m(F)$-orbit in the stable orbit
$\CO^\st_{[d^{c}1^{(ab+c)}]}$ corresponds to a pair $([d^c1^{(ab+c)}], (q_d,q_1))$ as discussed in Section 5.2. Put $\udl{q}:=(q_d,q_1)$.
Following Theorem 5.1, the automorphic $L^2$-packet $\wt{\Pi}_{\psi_0}(\epsilon_{\psi_0})$ of $G_m$ is not empty.
Take a member $\Theta$ in $\wt{\Pi}_{\psi_0}(\epsilon_{\psi_0})$.
The $\psi_{X_{\udl{q}}}$-Fourier coefficients $\CF^{\psi_{X_{\udl{q}}}}(\varphi_\Theta)$ of $\varphi_\Theta\in\Theta$,
as discussed in Section 5.2, are automorphic as functions on
$$
\SO_c(q_d)(\BA)\times\SO_{ab+c}(q_1)(\BA).
$$
Define the {\bf automorphic kernel functions} in this case to be
\begin{equation}\label{eksoe1}
\CK_{\varphi_\Theta}^{\psi_0}(h,g):=\CF^{\psi_{X_{\udl{q}}}}(\varphi_\Theta)(h,g)
\end{equation}
for $(h,g)\in\SO_c(q_d,\BA)\times\SO_{ab+c}(q_1,\BA)$. This family of the automorphic kernel functions
$\CK_{\varphi_\Theta}^{(\tau,b+c)}(h,g)$ should produce the endoscopy correspondence
$$
\SO_{ab+1}(q_0,\BA)\times\SO_c(q_d,\BA)\rightarrow\SO_{ab+c}(q_1,\BA)
$$
for {\bf Case (a=2e)}. Note that for the endoscopy group in $\wt\CE_\simp(ab+c-1)$, one takes $\SO_{ab+c}(q_1)$ to be an
$F$-split odd special orthogonal group with $2n+1=ab+c$ and $c=2f+1$, and hence both $\SO_{ab+1}$ and $\SO_c$ are
$F$-split odd special orthogonal groups as discussed in Section 3.1.
Recall from Section 5.2 that the condition in general is that the quadratic form
$q_d\oplus q_1$ shares the anisotropic kernel with the quadratic form $q_V$ defining $\SO_{2m}$. Hence the formulation here is
more general than the elliptic endoscopy transfer for the $F$-quasisplit case, and is related to \cite[Chapter 9]{Ar12} for
the transfers for inner forms.

More precisely, take $\psi_1:=(\tau,b)$. It is clear that $\psi_1$ is a simple global Arthur parameter of $\SO_{ab+1}(q_0)$ and
the global Arthur parameter $\psi_0$ is built from the global Arthur parameter $\psi_1$ of $G_0=\SO_{ab+1}(q_0)$ and the structures of
the groups $H=\SO_c(q_d)$ and $G=\SO_{ab+c}(q_1)$. The main conjecture (Conjecture \ref{mcec}) specializes to the following conjecture
for the current case. It is clear that when $b=0$, it produces the identity transfer from $\SO_c$ to itself; and when $c=1$ and $b=1$, this
construction reduces to the automorphic descent from $\GL_{a}$ to $\SO_{a+1}$ as discussed in \cite{GRS11}. The more general case of
$c=1$ will be discussed in Section 6.

\begin{conj}[$\SO_{2n+1}$: {\bf Case (a=2e) and (c=2f+1)}]\label{ecsoo1co}
For integers $a=2e\geq 2$, $b\geq 0$, $c=2f+1\geq 1$, and $d=a-1$, let $\tau\in\CA_\cusp(a)$ be self-dual and
$\psi_1=(\tau,b)$ be a simple global Arthur parameter of $F$-split $\SO_{ab+1}(q_0)$; and
let $\psi_2$ and $\psi:=\psi_1\boxplus\psi_2$ be global Arthur parameters of $\SO_c(q_d)$ and $\SO_{ab+c}(q_1)$, respectively.
For $\sigma\in\CA_2(\SO_c(q_d))$ and $\pi\in\CA_2(\SO_{ab+c}(q_1))$,
if there exists an
automorphic member $\Theta\in\wt{\Pi}_{\psi_0}(\epsilon_{\psi_0})$, such that the following integral
$$
\int_{\SO_c(q_d,F)\bks \SO_c(q_d,\BA)}\int_{\SO_{ab+c}(q_1,F)\bks \SO_{ab+c}(q_1,\BA)}
\CK_{\varphi_\Theta}^{\psi_0}(h,g)\varphi_\sigma(h)\ovl{\varphi_\pi(g)}dhdg
$$
is nonzero for some choice of $\varphi_\Theta\in\Theta$, $\varphi_\sigma\in\sigma$, and $\varphi_\pi\in\pi$, assuming the
convergence of the integral, then $\sigma\in\wt{\Pi}_{\psi_2}(\epsilon_{\psi_2})$ if and only if
$\pi\in\wt{\Pi}_{\psi}(\epsilon_{\psi})$.
\end{conj}

This construction of the automorphic kernel functions and their related endoscopy correspondence can be described by the following diagram:
$$
\begin{matrix}
&&&\SO_{a(b+c)}(q_V)&\\
&&&&\\
&&&\Theta\in\wt\Pi_{\psi_0}(\epsilon_{\psi_0})&\\
&&&&\\
&&\nearrow&\downarrow&\\
&&&&\\
&\GL_a,\tau&&\CK_{\varphi_\Theta}^{\psi_0}(h,g)&\\
&&&&\\
\swarrow&&&\downarrow&\\
&&&&\\
\SO_{ab+1}(q_0)&\times&\SO_{c}(q_d)&\longleftrightarrow&\SO_{ab+c}(q_1)\\
&&&&\\
\wt\Pi_{(\tau,b)}(\epsilon_{(\tau,b)})&&\wt\Pi_{\psi_2}(\epsilon_{\psi_2})&&\wt\Pi_{(\tau,b)\boxplus\psi_2}(\epsilon_{(\tau,b)\boxplus\psi_2})
\\
\\
\end{matrix}
$$
where $a=2e$, $c=2f+1$, and $2n=ab+c-1$. When $\tau$ is of symplectic type, $b=2l+1$ is odd; when $\tau$ is of orthogonal type, $b=2l$ is even.

Note that when $\tau$ is of symplectic type and $b=1$ and when the global Arthur parameter $\psi_2$ is generic, the construction
of the endoscopy transfer for irreducible generic cuspidal automorphic representations was outlined in \cite{G08} and see also
\cite{G12}. Also when $\tau$ is of symplectic type and $b=c=1$, this can be viewed as a case of the automorphic descent from
$\GL_{2e}$ to $\SO_{2e+1}$ and was completely proved in \cite{GRS11}.
The general case of this conjecture is the subject discussed in \cite{JZ13}.

Assume now that $c=2f$ is even. If $\tau$ is of symplectic type, then $b=2l$ is even; and if $\tau$ is of orthogonal type, then $b=2l+1$ is odd.
Consider the partition $[d^c1^{(ab+c)}]$ with $d=a-1$. Each $F$-rational adjoint $G_m(F)$-orbit in the stable orbit
$\CO^\st_{[d^{c}1^{(ab+c)}]}$ corresponds to a pair $([d^c1^{(ab+c)}], (q_d,q_1))$ as discussed in Section 5.2.
Put $\udl{q}:=(q_d,q_1)$, and
for $\Theta\in\wt{\Pi}_{\psi_0}(\epsilon_{\psi_0})$ of $G_m$,
the $\psi_{X_{\udl{q}}}$-Fourier coefficients $\CF^{\psi_{X_{\udl{q}}}}(\varphi_\Theta)$ of $\varphi_\Theta\in\Theta$,
as discussed in Section 5.2, is automorphic on
$$
\SO_c(q_d)(\BA)\times\SO_{ab+c}(q_1)(\BA).
$$
Note that the automorphic $L^2$-packet $\wt{\Pi}_{\psi_0}(\epsilon_{\psi_0})$ is not empty by Theorem 5.1.
The {\bf automorphic kernel functions} in this case are defined to be
\begin{equation}\label{eksoe2}
\CK_{\varphi_\Theta}^{\psi_0}(h,g):=\CF^{\psi_{X_{\udl{q}}}}(\varphi_\Theta)(h,g)
\end{equation}
for $(h,g)\in\SO_c(q_d,\BA)\times\SO_{ab+c}(q_1,\BA)$. This family of the automorphic kernel functions
$\CK_{\varphi_\Theta}^{(\tau,b+c)}(h,g)$ should produce the endoscopy correspondence
$$
\SO_{ab}(q_0,\BA)\times\SO_c(q_d,\BA)\rightarrow\SO_{ab+c}(q_1,\BA)
$$
for {\bf Case (a=2e)}. Note that for the endoscopy group in $\wt\CE_\simp(ab+c)$ with $a=2e$ and $c=2f$,
one takes $\SO_{ab+c}(q_1)$ to be an $F$-quasisplit even special orthogonal group, and hence both $\SO_{ab}(q_0)$ and $\SO_c(q_d)$ are
$F$-quasisplit even special orthogonal groups. Recall from Section 5.2 that the condition in general is that the quadratic form
$q_d\oplus q_1$ shares the anisotropic kernel with the quadratic form $q_V$ defining $\SO_{2m}$, which implies that
$\eta_{q_d}\cdot\eta_{q_1}=\eta_{q_V}$, following the definition of character $\eta$ in \cite[Section 1.2]{Ar12}. Hence in order for
$\SO_{ab}(q_0)\times\SO_c(q_d)$ to be an elliptic endoscopy group of $\SO_{ab+c}(q_1)$, one must have
$$
\eta_{q_1}=\eta_{q_0}\cdot\eta_{q_d}.
$$
Hence one must put the condition that
$$
\eta_{q_V}=\eta_{q_d}\cdot\eta_{q_1}=\eta_{q_0}\cdot\eta_{q_d}^2=\eta_{q_0}
$$
to make $\SO_{ab}(q_0)\times\SO_c(q_d)$ an elliptic endoscopy group of $\SO_{ab+c}(q_1)$ as discussed in Section 3.3.
Hence the formulation here is more general than the elliptic endoscopy transfer for the $F$-quasisplit case,
and is related to \cite[Chapter 9]{Ar12} for the transfers for inner forms.

More precisely, take $\psi_1:=(\tau,b)$. It is clear that $\psi_1$ is a simple global Arthur parameter of $\SO_{ab}(q_0)$ and
the global Arthur parameter $\psi_0$ is built from the global Arthur parameter $\psi_1$ of $G_0=\SO_{ab}(q_0)$ and the structures of
the groups $H=\SO_c(q_d)$ and $G=\SO_{ab+c}(q_1)$. The main conjecture (Conjecture \ref{mcec}) specializes to the following conjecture
for the current case. It is clear that when $b=0$, it produces the identity transfer from $\SO_c(q_d)$ to itself; and when $c=0$, this
construction reduces the identity transfer from $\SO_{ab}(q_0)$ to itself.
The more general discussion of this endoscopy correspondence and its relation to the automorphic descents will be discussed in Section 6.

\begin{conj}[$\SO_{2n}$: {\bf Case (a=2e) and (c=2f)}]\label{ecsoo2co}
For integers $a=2e\geq 2$, $b\geq 1$, $c=2f\geq 2$, and $d=a-1$, let $\tau\in\CA_\cusp(a)$ be self-dual and
$\psi_1=(\tau,b)$ be a simple global Arthur parameter of $F$-split $\SO_{ab}(q_0)$; and
let $\psi_2$ and $\psi:=\psi_1\boxplus\psi_2$ be global Arthur parameters of $\SO_c(q_d)$ and $\SO_{ab+c}(q_1)$, respectively.
For $\sigma\in\CA_2(\SO_c(q_d))$ and $\pi\in\CA_2(\SO_{ab+c}(q_1))$,
if there exists an
automorphic member $\Theta\in\wt{\Pi}_{\psi_0}(\epsilon_{\psi_0})$, such that the following integral
$$
\int_{\SO_c(q_d,F)\bks \SO_c(q_d,\BA)}\int_{\SO_{ab+c}(q_1,F)\bks \SO_{ab+c}(q_1,\BA)}
\CK_{\varphi_\Theta}^{\psi_0}(h,g)\varphi_\sigma(h)\ovl{\varphi_\pi(g)}dhdg
$$
is nonzero for some choice of $\varphi_\Theta\in\Theta$, $\varphi_\sigma\in\sigma$, and $\varphi_\pi\in\pi$, assuming the
convergence of the integral, then $\sigma\in\wt{\Pi}_{\psi_2}(\epsilon_{\psi_2})$ if and only if
$\pi\in\wt{\Pi}_{\psi}(\epsilon_{\psi})$.
\end{conj}

This construction of the automorphic kernel functions and their related endoscopy correspondence can be described by the following diagram:
$$
\begin{matrix}
&&&\SO_{a(b+c)}(q_V)&\\
&&&&\\
&&&\Theta\in\wt\Pi_{\psi_0}(\epsilon_{\psi_0})&\\
&&&&\\
&&\nearrow&\downarrow&\\
&&&&\\
&\GL_a,\tau&&\CK_{\varphi_\Theta}^{\psi_0}(h,g)&\\
&&&&\\
\swarrow&&&\downarrow&\\
&&&&\\
\SO_{ab}(q_0)&\times&\SO_{c}(q_d)&\longleftrightarrow&\SO_{ab+c}(q_1)\\
&&&&\\
\wt\Pi_{(\tau,b)}(\epsilon_{(\tau,b)})&&\wt\Pi_{\psi_2}(\epsilon_{\psi_2})&&\wt\Pi_{(\tau,b)\boxplus\psi_2}(\epsilon_{(\tau,b)\boxplus\psi_2})
\\
\\
\end{matrix}
$$
where $a=2e$, $c=2f$, and $2n=ab+c$. When $\tau$ is of symplectic type, $b=2l$ is even; when $\tau$ is of orthogonal type, $b=2l+1$ is odd.

This is the case has been discussed in \cite{Jn11}. The details of the proof of this conjecture for the case when $b=2$ and the global
Arthur parameter $\psi_2$ is generic are given in the two papers in \cite{GJS13}. When $\tau$ is of orthogonal type and $b=1$ and when the global Arthur parameter $\psi_2$ is generic, the construction
of the endoscopy transfer for irreducible generic cuspidal automorphic representations was outlined in \cite{G08} and see also
\cite{G12}. The general case of this conjecture will be the subject discussed in \cite{JZ13}.

{\bf Odd Special Orthogonal Group Case.}\
Assume that $G_m$ is the $F$-split odd special orthogonal group $\SO_{2m+1}$.
According to the discussion in Section 5.2 with $d=2k-1$,
for any automorphic form $\varphi$ on $G_m(\BA)$, the $\psi_{X_{\udl{q}}}$-Fourier coefficient
$\CF^{\psi_{X_{\udl{q}}}}(\varphi)$ of $\varphi$ is automorphic on
$$
\SO_c(q_d,\BA)\times\SO_{2m+1-cd}(q_1,\BA).
$$
Note that the integers $c$ and $2m+1-cd$ are in different parity. The functorial meaning between
automorphic forms $\CA(\SO_c(q_d))$ and automorphic forms $\CA(\SO_{2m+1-cd}(q_1))$ is not so easily related to the endoscopy structure in
general and will be considered in the forthcoming work of the author (\cite{Jn13}).
For {\bf Case (c=2f)}, it could be translated to the metaplectic double cover of $\Sp_{2m}(\BA)$ and formulated in the framework of
endoscopy correspondence below. However, for {\bf Case (c=2f+1)}, it is more mysterious (see \cite{Jn13}), except the
case where $c=1$ and $b=1$, which can be formulated as a case of automorphic descents of Ginzburg, Rallis and Soudry (\cite{GRS11}) as
follows.

Assume that $c=1$ and $b=1$. Since $a=2e$ and $d=2e+1=a+1$, the partition is $[(2e+1)1^{(2e)}]$ with $m=a$. This is the case that the
automorphic descent takes $\tau\in\CA_\cusp(2e)$ of orthogonal type to an cuspidal automorphic representation of $\SO_{2e}(\BA)$.
More general discussion on automorphic descents will be reviewed in Section 6. The general situation of this case will be discussed discussion
in \cite{Jn13}.

{\bf Metaplectic Cover of Symplectic Group Case:}\
Take $G_m(\BA)$ to be the metaplectic double cover $\tilsp_{2m}(\BA)$ of $\Sp_{2m}(\BA)$.
Consider that $2m=a(b+c)$ with $c=2f\geq 2$. As in the case of
$G_m=\Sp_{2m}$, attached to the partition $[d^c1^{(ab+c)}]$ with $d=2e-1$, the $F$-rational $G_m(F)$-adjoint nilpotent orbit is $F$-stable. The
associated $\psi_{X_{\udl{q}}}$-Fourier coefficient of any automorphic form $\varphi$ on $\tilsp_{2m}(\BA)$,
$\CF^{\psi_{X_{\udl{q}}}}(\varphi)$, is automorphic on
$$
\tilsp_{c}(\BA)\times\tilsp_{ab+c}(\BA).
$$
The meaning in endoscopy correspondence of this construction can be explained as follows.

Let $\tau\in\CA_\cusp(a)$ be self dual. Assume that if $\tau$ is of symplectic type, the central value $L(\frac{1}{2},\pi)$ of the
standard $L$-function (one could also use the partial $L$-function here) of $\pi$ is nonzero. When $\tau$ is of symplectic type, take
$b=2l+1$; and when $\tau$ is of orthogonal type, take $b=2l$. Then consider a global Arthur parameter
$$
\psi_0:=(\tau,b+c)
$$
Note that the dual group of $\tilsp_{2m}(\BA)$ is considered to be $\Sp_{2m}$ and the theory of stable trace formula and the theory of
endoscopy for $\tilsp_{2m}(\BA)$ have been developing through a series of recent papers, including his PhD thesis of W.-W. Li (\cite{Liw11}),
although the complete theory is still in progress.
Following the arguments in \cite{JLZ12}, it is not hard to see that Theorem 5.1 holds for $\tilsp_{2m}(\BA)$ and hence the
automorphic $L^2$-packet $\wt{\Pi}_{\psi_0}(\epsilon_{\psi_0})$ is not empty (\cite{Liu13}). For $\Theta\in\wt{\Pi}_{\psi_0}(\epsilon_{\psi_0})$,
Define the {\bf automorphic kernel functions} in this case to be
\begin{equation}\label{ekspt}
\CK_{\varphi_\Theta}^{\psi_0}(h,g):=\CF^{\psi_{X_{\udl{q}}}}(\varphi_\Theta)(h,g)
\end{equation}
for $(h,g)\in\tilsp_c(\BA)\times\tilsp_{ab+c}(\BA)$.
The endoscopy correspondence is given by
$$
\SO_{ab+1}(q_0,\BA)\times\tilsp_c(\BA)\rightarrow\tilsp_{ab+c}(\BA),
$$
or by
$$
\tilsp_{ab}(\BA)\times\tilsp_c(\BA)\rightarrow\tilsp_{ab+c}(\BA).
$$
Note that this case can be regarded as a variant of the case $\SO_{2n+1}$ with {\bf (a=2e)} and $2n=ab+c$
as discussed in Section 3.3; and that the odd special orthogonal group $\SO_{ab+1}(q_0)$ is $F$-split.

More precisely, take $\psi_1:=(\tau,b)$. It is clear that $\psi_1$ is a simple global Arthur parameter of
$\SO_{ab+1}(q_0,\BA)$ or of $\tilsp_{ab}(\BA)$, and
the global Arthur parameter $\psi_0$ is built from the global Arthur parameter $\psi_1$ of $G_0=\SO_{ab+1}(q_0)$ or $\tilsp_{ab}(\BA)$,
and the structures of
the groups $H=\tilsp_c(\BA)$ and $G=\tilsp_{ab+c}(\BA)$.
The main conjecture (Conjecture \ref{mcec}) can also be specialized
to the following conjecture for the current case.

\begin{conj}[$\tilsp_{2n}(\BA)$: {\bf Case (a=2e) and (c=2f)}]\label{ecsptco}
Assume that integers $a=2e\geq 2$, $b\geq 1$, $c=2f\geq 2$, and $d=a-1$. Let $\tau\in\CA_\cusp(a)$ be self-dual with
the assumption that if $\tau$ is of symplectic type, then $L(\frac{1}{2},\pi)\neq 0$.
Let $\psi_1=(\tau,b)$ be a simple global Arthur parameter of $F$-split $\SO_{ab+1}(q_0,\BA)$ or of $\tilsp_{ab}(\BA)$; and
let $\psi_2$ and $\psi:=\psi_1\boxplus\psi_2$ be global Arthur parameters of $\tilsp_c(\BA)$ and $\tilsp_{ab+c}(\BA)$, respectively.
For $\sigma\in\CA_2(\tilsp_c(\BA))$ and $\pi\in\CA_2(\tilsp_{ab+c}(\BA))$,
if there exists an
automorphic member $\Theta\in\wt{\Pi}_{\psi_0}(\epsilon_{\psi_0})$, such that the following integral
$$
\int_{\Sp_c(F)\bks \tilsp_c(\BA)}\int_{\Sp_{ab+c}(F)\bks \tilsp_{ab+c}(\BA)}
\CK_{\varphi_\Theta}^{\psi_0}(h,g)\varphi_\sigma(h)\ovl{\varphi_\pi(g)}dhdg
$$
is nonzero for some choice of $\varphi_\Theta\in\Theta$, $\varphi_\sigma\in\sigma$, and $\varphi_\pi\in\pi$, assuming the
convergence of the integral, then $\sigma\in\wt{\Pi}_{\psi_2}(\epsilon_{\psi_2})$ if and only if
$\pi\in\wt{\Pi}_{\psi}(\epsilon_{\psi})$.
\end{conj}

This construction of the automorphic kernel functions and their related endoscopy correspondence can be described by the following diagram:
$$
\begin{matrix}
&&&\tilsp_{a(b+c)}(\BA)&\\
&&&&\\
&&&\Theta\in\wt\Pi_{\psi_0}(\epsilon_{\psi_0})&\\
&&&&\\
&&\nearrow&\downarrow&\\
&&&&\\
&\GL_a,\tau&&\CK_{\varphi_\Theta}^{\psi_0}(h,g)&\\
&&&&\\
\swarrow&&&\downarrow&\\
&&&&\\
\SO_{ab+1}(q_0,\BA)&\times&\tilsp_{c}(\BA)&\longleftrightarrow&\tilsp_{ab+c}(\BA)\\
&&&&\\
\wt\Pi_{(\tau,b)}(\epsilon_{(\tau,b)})&&\wt\Pi_{\psi_2}(\epsilon_{\psi_2})&&\wt\Pi_{(\tau,b)\boxplus\psi_2}(\epsilon_{(\tau,b)\boxplus\psi_2})
\\
\\
\end{matrix}
$$
where $a=2e$, $c=2f$, and $2n=ab+c$. When $\tau$ is of symplectic type, $b=2l+1$ is odd; when $\tau$ is of orthogonal type, $b=2l$ is even.

Some preliminary cases of this conjecture is considered in \cite{JL13}.

\subsection{Automorphic kernel functions: Case (d=2k)}\label{ec(2)}
The constructions of automorphic kernel functions with $a=2e+1\geq 1$ being odd integers are considered here.
In this case, the cuspidal automorphic representation $\tau$ must be of orthogonal type. Following Waldspurger (\cite[Chapitre I]{W01},
the structure of $F$-rational $G_m(F)$-adjoint (nilpotent) orbits in the $F$-stable orbit $\CO^\st_{[d^c1^*]}$ associated to
a partition of type $[d^c1^*]$ was discussed in Section 5.2. Note that when $d=2k$ is even, the subspaces $\frak{g}_{\pm 1}$ are
non-zero. Hence there exists a generalized Jacobi subgroup, containing $\Stab_{L_0}(\psi_{X_{\udl{q}}})$ as the semi-simple part,
of $G_m$ which stabilizes the character $\psi_{X_{\udl{q}}}$. More precisely, following \eqref{sl2}, define $V_1$ to be the unipotent
subgroup of $G_m$ whose Lie algebra is $\oplus_{i\geq 1}\frak{g}_i$. It is clear that $V_2=V_{X_{\udl{q}}}$ is a normal subgroup of $V_1$ and
the quotient $V_1/V_{X_{\udl{q}}}$, which is abelian, is isomorphic to $\frak{g}_1$. Let $\ker_{V_2}(\psi_{X_{\udl{q}}})$ be the
kernel of the character $\psi_{X_{\udl{q}}}$ in $V_2=V_{X_{\udl{q}}}$. It is easy to check that $V_1/\ker_{V_2}(\psi_{X_{\udl{q}}})$ is
the Heisenberg group of dimension $\dim\frak{g}_1+1$.

In order to define the Fourier-Jacobi coefficients of automorphic forms, let $\frak{g}_1:=\frak{g}_1^+\oplus\frak{g}_1^-$ be a
polarization of the symplectic space $\frak{g}_1$ via the Lie structure of $\frak{g}$ and modulo $\oplus_{i\geq 2}\frak{g}_i$. Let
$\theta^{\psi_{X_{\udl{q}}}}_\phi$ be the theta function of the (generalized) Jacobi group
$$
(\wt{\Stab_{L_0}(\psi_{X_{\udl{q}}})}\ltimes V_1/\ker_{V_2}(\psi_{X_{\udl{q}}}))(\BA)
$$
attached to a Bruhat-Schwartz function $\phi\in\CS(\frak{g}_1^+(\BA))$. For an automorphic form $\varphi$ on $G_m(\BA)$, the
$\psi_{X_{\udl{q}}}$-Fourier-Jacobi coefficient of $\varphi$ is defined by
\begin{equation}\label{fjc}
\fj^{\psi_{X_{\udl{q}}}}_\phi(\varphi)(x):=\int_{V_1(F)\bks V_1(\BA)}\varphi(vx)\ovl{\theta^{\psi_{X_{\udl{q}}}}_\phi(vx)}dv
\end{equation}
where $x\in \wt{\Stab_{L_0}(\psi_{X_{\udl{q}}})}(\BA)$.

Following \cite{PS83} and \cite{Ik94}, the basic theory of the Fourier-Jacobi
coefficients of automorphic forms asserts that for any irreducible automorphic representation $\Theta\in\CA(G_m)$,
the $\psi_{X_{\udl{q}}}$-Fourier coefficient, $\CF^{\psi_{X_{\udl{q}}}}(\varphi_\Theta)$ is nonzero for some $\varphi_\Theta\in\Theta$
if and only if the following (generalized) Fourier-Jacobi coefficient $\fj^{\psi_{X_{\udl{q}}}}_\phi(\varphi_\Theta)(x)$
is nonzero for certain choice of data. It is this Fourier-Jacobi coefficient that helps the construction of the {\bf automorphic kernel functions}
for {\bf Case (d=2k)}.

{\bf Symplectic group case.}\
Assume that $G_m$ is the symplectic group $\Sp_{2m}$. The partition is $[d^c1^{(2m-cd)}]$ with $d=2e$. It follows that $c=2f\geq 2$ or
$c=2f+1\geq 1$. Take $2m+1=a(b+c)$. Since $\tau\in\CA_\cusp(a)$ is of orthogonal type, $b+c$ must be odd. Consider the global
Arthur parameter
$$
\psi_0:=(\tau,b+c).
$$
By Theorem 5.1, the automorphic $L^2$-packet $\wt\Pi_{\psi_0}(\epsilon_{\psi_0})$ is not empty.
Take a member $\Theta\in\wt\Pi_{\psi_0}(\epsilon_{\psi_0})$ of $G_m(\BA)$. The Fourier-Jacobi coefficients
$\fj^{\psi_{X_{\udl{q}}}}_\phi(\varphi_\Theta)(h,g)$ of $\varphi_\Theta\in\Theta$ define automorphic functions on
$$
\wt\SO_c(q_d,\BA)\times\tilsp_{2n}(\BA)
$$
with $2n=ab+c-1$. Note that this pair forms a reductive dual pair in the sense of R. Howe in the theory of theta correspondence and hence
it is known that if $c=2f$, the covering splits and one has
$$
\SO_{2f}(q_d,\BA)\times\Sp_{2n}(\BA);
$$
and if $c=2f+1$, the covering splits over $\wt\SO_c(q_d,\BA)$ and one has
$$
\SO_{2f+1}(q_d,\BA)\times\tilsp_{2n}(\BA).
$$
Note that by Section 3.2, the special orthogonal group $\SO_c(q_d)$ may not be $F$-quasisplit.

Define the {\bf automorphic kernel functions} in this case to be
\begin{equation}\label{akspao}
\CK^{\psi_0}_{\varphi_\Theta}(h,g):=\fj^{\psi_{X_{\udl{q}}}}_\phi(\varphi_\Theta)(h,g)
\end{equation}
for $(h,g)\in\wt\SO_c(q_d,\BA)\times\tilsp_{2n}(\BA)$, accordingly. This family of the automorphic kernel functions
$\CK^{\psi_0}_{\varphi_\Theta}(h,g)$ should produce the endoscopy correspondence
$$
G_0(\BA)\times \wt{\SO}_c(q_d,\BA)\rightarrow\tilsp_{2n}(\BA).
$$

When $c=2f$, which implies that $b=2l+1$, one has
\begin{equation}
\Sp_{ab-1}(\BA)\times \SO_{2f}(q_d,\BA)\rightarrow\Sp_{2n}(\BA).
\end{equation}
Take $\psi_1:=(\tau,b)$ with $b=2l+1$, which is a global Arthur parameter for $G_0=\Sp_{ab-1}$.
This construction of the automorphic kernel functions and their related endoscopy correspondence can be described by the following diagram:
$$
\begin{matrix}
&&&\Sp_{a(b+c)-1}&\\
&&&&\\
&&&\Theta\in\wt\Pi_{\psi_0}(\epsilon_{\psi_0})&\\
&&&&\\
&&\nearrow&\downarrow&\\
&&&&\\
&\GL_a,\tau&&\CK_{\varphi_\Theta}^{\psi_0}(h,g)&\\
&&&&\\
\swarrow&&&\downarrow&\\
&&&&\\
\Sp_{ab-1}&\times&\SO_{c}(q_d)&\longleftrightarrow&\Sp_{ab+c-1}\\
&&&&\\
\wt\Pi_{(\tau,b)}(\epsilon_{(\tau,b)})&&\wt\Pi_{\psi_2}(\epsilon_{\psi_2})&&\wt\Pi_{(\tau,b)\boxplus\psi_2}(\epsilon_{(\tau,b)\boxplus\psi_2})
\\
\\
\end{matrix}
$$
where $a=2e+1$, $c=2f$, $b=2l+1$, $2n=ab+c-1$, and $\tau$ is of orthogonal type.
The main conjecture (Conjecture \ref{mcec}) specializes to the following conjecture
for the current case. Note that if $c=2f=0$, it is the identity transfer for $\Sp_{ab-1}(\BA)$.

\begin{conj}[$\Sp_{2n}$: {\bf Case (a=2e+1) and (c=2f)}]\label{ecsp1co}
Assume that integers $a=2e+1\geq 1$ with $d=a-1$, and $b=2l+1\geq 1$ and $c=2f\geq 0$.
Let $\tau\in\CA_\cusp(a)$ be of orthogonal type.  Let $\psi_1=(\tau,b)$ be a simple global Arthur parameter of for $\Sp_{ab-1}$,
$\psi_2$ be a global Arthur parameters of $\SO_c(q_d)$, and $\psi:=\psi_1\boxplus\psi_2$ be the global Arthur parameter for
$G_n(\BA):=\Sp_{2n}(\BA)$, respectively.
For $\sigma\in\CA_2(\SO_c(q_d))$ and $\pi\in\CA_2(G_n)$,
if there exists an
automorphic member $\Theta\in\wt{\Pi}_{\psi_0}(\epsilon_{\psi_0})$, such that the following integral
$$
\int_{\SO_c(q_d,F)\bks \SO_c(q_d,\BA)}\int_{G_n(F)\bks G_n(\BA)}
\CK_{\varphi_\Theta}^{\psi_0}(h,g)\varphi_\sigma(h)\ovl{\varphi_\pi(g)}dhdg
$$
is nonzero for some choice of $\varphi_\Theta\in\Theta$, $\varphi_\sigma\in\sigma$, and $\varphi_\pi\in\pi$, assuming the
convergence of the integral, then $\sigma\in\wt{\Pi}_{\psi_2}(\epsilon_{\psi_2})$ if and only if
$\pi\in\wt{\Pi}_{\psi}(\epsilon_{\psi})$.
\end{conj}

Note that when $b=1$ and the global Arthur parameter $\psi_2$ is generic, an outline of discussion of this conjecture for
endoscopy transfer of irreducible generic cuspidal automorphic representations was given in
\cite{G08} and see also \cite{G12}.

When $c=2f+1$, which implies that $b=2l$, one has
\begin{equation}
\SO_{ab+1}(q_0,\BA)\times \SO_{2f+1}(q_d,\BA)\rightarrow\tilsp_{2n}(\BA).
\end{equation}
This is exactly the endoscopy transfer introduced by Wen-Wei Li in \cite{Liw11} when both $\SO_{ab+1}(q_0)$ and $\SO_{2f+1}(q_d)$
are $F$-split. Recall from Section 3.2 that in general $\SO_{ab+1}(q_0)$ and $\SO_{2f+1}(q_d)$ may not be $F$-quasisplit. This
general case will be discussed elsewhere.

Take $\psi_1:=(\tau,b)$, which is a global Arthur parameter for $\SO_{ab+1}(q_0)$.
The construction of the automorphic kernel functions and their related endoscopy correspondence in this case
can be described by the following diagram:
$$
\begin{matrix}
&&&\Sp_{a(b+c)-1}&\\
&&&&\\
&&&\Theta\in\wt\Pi_{\psi_0}(\epsilon_{\psi_0})&\\
&&&&\\
&&\nearrow&\downarrow&\\
&&&&\\
&\GL_a,\tau&&\CK_{\varphi_\Theta}^{\psi_0}(h,g)&\\
&&&&\\
\swarrow&&&\downarrow&\\
&&&&\\
\SO_{ab+1}(q_0)&\times&\SO_{c}(q_d)&\longleftrightarrow&\tilsp_{ab+c-1}(\BA)\\
&&&&\\
\wt\Pi_{(\tau,b)}(\epsilon_{(\tau,b)})&&\wt\Pi_{\psi_2}(\epsilon_{\psi_2})&&\wt\Pi_{(\tau,b)\boxplus\psi_2}(\epsilon_{(\tau,b)\boxplus\psi_2})
\\
\\
\end{matrix}
$$
where $a=2e+1$, $c=2f+1$, $b=2l$, $2n=ab+c-1$, and $\tau$ is of orthogonal type.

The main conjecture (Conjecture \ref{mcec}) specializes to the following conjecture
for the current case of endoscopy correspondences.

\begin{conj}[$\tilsp_{2n}$: {\bf Case (a=2e+1) and (c=2f+1)}]\label{ecspt2co}
Assume that integers $a=2e+1\geq 1$ with $d=a-1$, and $b=2l\geq 2$ and $c=2f+1\geq 1$, and assume that
$\tau\in\CA_\cusp(a)$ is of orthogonal type.  Let $\psi_1=(\tau,b)$ be a simple global Arthur parameter for $\SO_{ab+1}(q_0)$,
$\psi_2$ be a global Arthur parameter of $\SO_c(q_d)$, and $\psi:=\psi_1\boxplus\psi_2$ be the global Arthur parameter for
$G_n(\BA):=\tilsp_{2n}(\BA)$, respectively.
For $\sigma\in\CA_2(\SO_c(q_d))$ and $\pi\in\CA_2(G_n)$,
if there exists an
automorphic member $\Theta\in\wt{\Pi}_{\psi_0}(\epsilon_{\psi_0})$, such that the following integral
$$
\int_{\SO_c(q_d,F)\bks \SO_c(q_d,\BA)}\int_{G_n(F)\bks G_n(\BA)}
\CK_{\varphi_\Theta}^{\psi_0}(h,g)\varphi_\sigma(h)\ovl{\varphi_\pi(g)}dhdg
$$
is nonzero for some choice of $\varphi_\Theta\in\Theta$, $\varphi_\sigma\in\sigma$, and $\varphi_\pi\in\pi$, assuming the
convergence of the integral, then $\sigma\in\wt{\Pi}_{\psi_2}(\epsilon_{\psi_2})$ if and only if
$\pi\in\wt{\Pi}_{\psi}(\epsilon_{\psi})$.
\end{conj}

Some preliminary cases of this conjecture will be treated in \cite{JL13}.

{\bf Example 2.}\ A more explicit explanation of the construction of the Fourier-Jacobi coefficients is given
for a phototype example.
Let $\tau\in\CA_\cusp(a)$ with $a=2e+1$ be self-dual. Then $\tau$ must be of orthogonal type. Also
$b=2l+1$ is odd and $c=2e$ is even. Assume that $2n=c+ab-1$.

Consider a standard parabolic subgroup $P_{c^e}$ of $\Sp_{a(c+b)-1}$, whose Levi decomposition is
$$
Q_{c^e}=M_{c^e}V_{c^e}
=(\GL_{c}^{\times(e)}\times\Sp_{ab+c-1})V_{c^e}
$$
with $2n=ab+c-1$.
The elements in the unipotent radical $V_{c^e}$, which is equal to $V_1$, are of form
\begin{equation}
v=
\begin{pmatrix}
I_{c}&X_1&&&&&&\cdots\\
&I_{c}&&&&&*&\cdots\\
&&\cdots&&&&&\cdots\\
&&&I_{c}&X_{e-1}&*&*&\cdots\\
&&&&I_{c}&Y&X_0&\cdots\\
&&&&&I_{2n}&Y^*&\cdots\\
&&&&&&I_{c}&\cdots\\
&&&&&&&\cdots
\end{pmatrix}.
\end{equation}
The subgroup consisting of all elements with $Y=0$ is a normal subgroup of $V_{c^e}$, which is equal to $V_{X_{\udl{q}}}$.
The character $\psi_{X_{\udl{q}}}$ of $V_{X_{\udl{q}}}$ is given by
\begin{equation}
\psi_{X_{\udl{q}}}(v)
:=
\psi_F(\tr(X_1+\cdots+X_{e-1})+\tr(S_{\udl{q}}X_0)),
\end{equation}
where $S_{\udl{q}}$ is a non-degenerate $c\times c$-symmetric matrix defining $\udl{q}$.
The stabilizer of the character $\psi_{X_{\udl{q}}}$ in
the Levi subgroup $\GL_{c}^{\times(e)}\times\Sp_{2n}$ is isomorphic to
$$
\SO_c(S_{\udl{q}})\times\Sp_{2n}.
$$
The elements $(h,g)$ of $\SO_c(S_{\udl{q}})\times\Sp_{2n}$ is embedded into $\GL_{c}^{\times(e)}\times\Sp_{2n}$ by
$$
(h,g)\mapsto (h^{\Delta (e)},g)
$$
where $h\mapsto h^{\Delta (e)}$ is the diagonal embedding $\SO_c(S_{\udl{q}})$ into $\GL_{c}^{\times(e)}$.

{\bf Even special orthogonal group case.}\
Assume that $G_m$ is an $F$-quasisplit even special orthogonal group $\SO_{2m}(q_V)$
defined by the $2m$-dimensional non-degenerate quadratic form $q_V$.

For the partition $[d^c1^{(2m-cd)}]$ for $\SO_{2m}(q_V)$,
$d=a-1=2e$ is even and hence one must have that $c=2f$ is even. In this case, take that $2m=a(b+c)+1$ with $b+c$ odd and
hence $b=2l+1$ must be odd. For $\tau\in\CA_\cusp(a)$ to be of orthogonal type, take the global Arthur parameter
$$
\psi_0:=(\tau,b+c)\boxplus(\chi,1)
$$
for $\SO_{2m}(q_V)$, which is the lift of the simple global Arthur parameter for $\Sp_{a(b+c)-1}$ to $\SO_{a(b+c)+1}(q_V)$ via the
theta correspondence, based on the theta lift from $\Sp_{a-1}$ to $\SO_{a+1}(q)$ with $q$ and $q_V$ share the same anisotropic kernel.

Following \cite{JLZ12}, one can produce a nonzero residual representation in the automorphic $L^2$-packet $\wt{\Pi}_{\psi_0}(\epsilon_{\psi_0})$,
with cuspidal support
$$
(\GL_a^{\times(l+f)}\times\SO_{a+1}(q),\tau^{\otimes(l+f)}\otimes\epsilon)
$$
where $\epsilon$ is the image of $\tau$ under the composition of the automorphic descent of $\GL_a$ to $\Sp_{a-1}$ with the
theta correspondence from $\Sp_{a-1}$ to $\SO_{a+1}$.
Take a member $\Theta$ in the set $\wt{\Pi}_{\psi_0}(\epsilon_{\psi_0})$ of $G_m$. The Fourier-Jacobi coefficients
$\fj^{\psi_{X_{\udl{q}}}}_\phi(\varphi_\Theta)(h,g)$ of $\varphi_\Theta\in\Theta$ define automorphic functions on
$$
\Sp_c(\BA)\times\SO_{ab+c+1}(q_1,\BA).
$$
Define the {\bf automorphic kernel functions} in this case to be
\begin{equation}\label{eksoe3}
\CK_{\varphi_\Theta}^{\psi_0}(h,g):=\fj^{\psi_{X_{\udl{q}}}}_\phi(\varphi_\Theta)(h,g)
\end{equation}
for $(h,g)\in\Sp_c(\BA)\times\SO_{ab+c+1}(q_1,\BA)$. This family of the automorphic kernel functions
$\CK_{\varphi_\Theta}^{\psi_0}(h,g)$ should produce the endoscopy correspondence
$$
\Sp_{ab-1}(\BA)\times\Sp_c(\BA)\rightarrow\SO_{ab+c+1}(q_1,\BA).
$$
Note that $\SO_{ab+c+1}(q_1)$ is the $F$-quasisplit even special orthogonal group defined by $q_1$, which shares the same $F$-anisotropic
kernel with $q_V$.

The main conjecture (Conjecture \ref{mcec}) specializes to the following conjecture
for the current case.

\begin{conj}[$\SO_{2n}$: {\bf Case (a=2e+1) and (c=2f)}]\label{ecsoeaoce}
Assume that integers $a=2e+1\geq 1$ with $d=a-1$, and $b=2l+1\geq 1$ and $c=2f\geq 2$.
Let $\tau\in\CA_\cusp(a)$ be of orthogonal type.  Let $\psi_1=(\tau,b)$ be a simple global Arthur parameter of for $\Sp_{ab-1}$.
Let $\psi_2$ be a global Arthur parameters of $\Sp_c$ and let $\psi:=\psi_1\boxplus\psi_2$ be the global Arthur parameter for
$\SO_{2n}(\BA)$ with $2n=ab+c+1$. For $\sigma\in\CA_2(\Sp_c)$ and $\pi\in\CA_2(\SO_{2n})$,
if there exists an
automorphic member $\Theta\in\wt{\Pi}_{\psi_0}(\epsilon_{\psi_0})$, such that the following integral
$$
\int_{\Sp_c(F)\bks \Sp_c(\BA)}\int_{\SO_{2n}(F)\bks\SO_{2n}(\BA)}
\CK_{\varphi_\Theta}^{\psi_0}(h,g)\varphi_\sigma(h)\ovl{\varphi_\pi(g)}dhdg
$$
is nonzero for some choice of $\varphi_\Theta\in\Theta$, $\varphi_\sigma\in\sigma$, and $\varphi_\pi\in\pi$, assuming the
convergence of the integral, then $\sigma\in\wt{\Pi}_{\psi_2}(\epsilon_{\psi_2})$ if and only if
$\pi\in\wt{\Pi}_{\psi}(\epsilon_{\psi})$.
\end{conj}

This construction of the automorphic kernel functions and their related endoscopy correspondence can be described by the following diagram:
$$
\begin{matrix}
&&&\SO_{a(b+c)+1}(q_V)&\\
&&&&\\
&&&\Theta\in\wt\Pi_{\psi_0}(\epsilon_{\psi_0})&\\
&&&&\\
&&\nearrow&\downarrow&\\
&&&&\\
&\GL_a,\tau&&\CK_{\varphi_\Theta}^{\psi_0}(h,g)&\\
&&&&\\
\swarrow&&&\downarrow&\\
&&&&\\
\Sp_{ab-1}&\times&\Sp_{c}&\longleftrightarrow&\SO_{ab+c+1}(q_1,\BA)\\
&&&&\\
\wt\Pi_{(\tau,b)}(\epsilon_{(\tau,b)})&&\wt\Pi_{\psi_2}(\epsilon_{\psi_2})&&\wt\Pi_{(\tau,b)\boxplus\psi_2}(\epsilon_{(\tau,b)\boxplus\psi_2})
\\
\\
\end{matrix}
$$
where $a=2e+1$, $c=2f$, $b=2l+1$, $2n=ab+c+1$, and $\tau$ is of orthogonal type.

Note that when $b=1$ and the global Arthur parameter $\psi_2$ is generic, an outline of discussion of this conjecture was given
for the endoscopy transfer of irreducible generic cuspidal automorphic representations in \cite{G08} and see also
\cite{G12}. The general case is the subject that will be discussed in \cite{JZ13}.

{\bf Odd special orthogonal group case.}\
Assume that $G_m$ is $F$-split odd special orthogonal group $\SO_{2m+1}$. The partition is $[d^c1^{(2m+1-cd)}]$ with $d=a-1=2e$, and
hence $c=2f$. Take $2m=a(b+c)$ with $b=2l$. Consider the simple global Arthur parameter for $\SO_{2m+1}$:
$$
\psi_0:=(\tau,b+c)
$$
with $\tau\in\CA_\cusp(a)$ being of orthogonal type. Take $\Theta\in\wt{\Pi}_{\psi_0}(\epsilon_{\psi_0})$ of $G_m$.
The Fourier-Jacobi coefficients $\fj^{\psi_{X_{\udl{q}}}}_\phi(\varphi_\Theta)(h,g)$ of $\varphi_\Theta\in\Theta$ define automorphic functions on
$$
\tilsp_c(\BA)\times\SO_{ab+c+1}(q_1,\BA).
$$
Define the {\bf automorphic kernel functions} in this case to be
\begin{equation}\label{eksoo}
\CK_{\varphi_\Theta}^{\psi_0}(h,g):=\fj^{\psi_{X_{\udl{q}}}}_\phi(\varphi_\Theta)(h,g)
\end{equation}
for $(h,g)\in\Sp_c(\BA)\times\SO_{ab+c+1}(q_1,\BA)$. This family of the automorphic kernel functions
$\CK_{\varphi_\Theta}^{\psi_0}(h,g)$ should produce the endoscopy correspondence
$$
\tilsp_{ab}(\BA)\times\tilsp_c(\BA)\rightarrow\SO_{ab+c+1}(q_1,\BA).
$$
Note that $\SO_{ab+c+1}(q_1)$ is the $F$-split odd special orthogonal group defined by $q_1$. It is possible to consider
this transfer for more general $q_V$ and $q_1$.
The main conjecture (Conjecture \ref{mcec}) specializes to the following conjecture
for the current case of endoscopy correspondences.

\begin{conj}[$\SO_{2n+1}$: {\bf Case (a=2e+1) and (c=2f)}]\label{ecsooaoce}
Assume that integers $a=2e+1\geq 1$ with $d=a-1$, and $b=2l, c=2f\geq 2$.
Let $\tau\in\CA_\cusp(a)$ be of orthogonal type.  Let $\psi_1=(\tau,b)$ be a simple global Arthur parameter of for $\tilsp_{ab}(\BA)$.
Let $\psi_2$ be a global Arthur parameters of $\tilsp_c(\BA)$ and let $\psi:=\psi_1\boxplus\psi_2$ be the global Arthur parameter for
$\SO_{2n+1}(\BA)$ with $2n=ab+c$. For $\sigma\in\CA_2(\tilsp_c(\BA))$ and $\pi\in\CA_2(\SO_{2n+1})$,
if there exists an
automorphic member $\Theta\in\wt{\Pi}_{\psi_0}(\epsilon_{\psi_0})$, such that the following integral
$$
\int_{\Sp_c(F)\bks \tilsp_c(\BA)}\int_{\SO_{2n+1}(F)\bks\SO_{2n+1}(\BA)}
\CK_{\varphi_\Theta}^{\psi_0}(h,g)\varphi_\sigma(h)\ovl{\varphi_\pi(g)}dhdg
$$
is nonzero for some choice of $\varphi_\Theta\in\Theta$, $\varphi_\sigma\in\sigma$, and $\varphi_\pi\in\pi$, assuming the
convergence of the integral, then $\sigma\in\wt{\Pi}_{\psi_2}(\epsilon_{\psi_2})$ if and only if
$\pi\in\wt{\Pi}_{\psi}(\epsilon_{\psi})$.
\end{conj}

The construction of the automorphic kernel functions and their related endoscopy correspondence in this case
can be described by the following diagram:
$$
\begin{matrix}
&&&\SO_{a(b+c)+1}(q_V)&\\
&&&&\\
&&&\Theta\in\wt\Pi_{\psi_0}(\epsilon_{\psi_0})&\\
&&&&\\
&&\nearrow&\downarrow&\\
&&&&\\
&\GL_a,\tau&&\CK_{\varphi_\Theta}^{\psi_0}(h,g)&\\
&&&&\\
\swarrow&&&\downarrow&\\
&&&&\\
\tilsp_{ab}(\BA)&\times&\tilsp_{c}(\BA)&\longleftrightarrow&\SO_{ab+c+1}(q_1,\BA)\\
&&&&\\
\wt\Pi_{(\tau,b)}(\epsilon_{(\tau,b)})&&\wt\Pi_{\psi_2}(\epsilon_{\psi_2})&&\wt\Pi_{(\tau,b)\boxplus\psi_2}(\epsilon_{(\tau,b)\boxplus\psi_2})
\\
\\
\end{matrix}
$$
where $a=2e+1$, $c=2f$, $b=2l$, $2n=ab+c$, and $\tau$ is of orthogonal type.

Note that the discussion here assumes that $d=a-1$, and see \cite{Jn13} for other possible choices of the pairs $(a,d)$.

\subsection{Unitary group case}
Assume that $G_m$ is an $F$-quasisplit unitary group $\RU_{E/F}(m_V)=\RU(m_V)=\RU_{m_V}(q_V)$
defined by the $m_V$-dimensional non-degenerate hermitian form $q_V$, with $m=[\frac{m_V}{2}]$, and assume that $\kappa$ is the sign of
in the endoscopy data
$$
(\RU(m_V),\xi_{\chi_\kappa})=(\RU(m_V),{\chi_\kappa})
$$
of $\RR_{E/F}(m_V)$. Further, assume that $d=a-1$ and that $m_V=a(b+c)$ or $m_V=a(b+c)+1$. Recall from Section 5.2 that the Fourier coefficients
associated to the partition $[d^c1^{(m_V-cd)}]$ produce the stabilizer
$$
\RU_c(q_d)\times\RU_{m_V-cd}(q_1).
$$
Assume in the rest of this section that $\RU_c(q_d)\times\RU_{m_V-cd}(q_1)$ is $F$-quasisplit and hence is denoted by
$$
\RU(c)\times\RU(m_V-cd).
$$
The more general cases will be considered elsewhere.

The (standard or twisted) elliptic endoscopy is to consider the following transfer
\begin{equation}\label{etug}
\RU(m_V-ac)\times\RU(c)\rightarrow\RU(m_V-dc).
\end{equation}
The signs $\kappa_{m_V-ac}$, $\kappa_c$ and $\kappa_{m_V-dc}$ of the $L$-embedding for $\RU(m_V-ac)$, $\RU(c)$ and $\RU(m_V-dc)$,
respectively, can be determined as follows. It is clear that $\RU(m_V-dc)$ is part of the Levi subgroup of $\RU(m_V)$, and hence
$$
\kappa_{m_V-dc}=\kappa=\kappa_{m_V}.
$$
By \cite[Section 2.4]{Mk12}, one must have that
$$
\kappa_{m_V-ac}(-1)^{m_V-ac-1}=\kappa_c(-1)^{c-1}=\kappa(-1)^{m_V-dc-1}.
$$
It follows that
\begin{eqnarray}\label{etsgn}
\kappa_{m_V-ac}&=&\kappa(-1)^c\nonumber\\
\kappa_c&=&\kappa(-1)^{m_V-ac}.
\end{eqnarray}
Therefore, the elliptic endoscopy considered here is to consider the transfer
\begin{equation}\label{eesgn}
(\RU(m_V-ac)\times\RU(c),(\kappa(-1)^c,\kappa(-1)^{m_V-ac}))\rightarrow(\RU(m_V-dc),\kappa).
\end{equation}

The {\bf automorphic kernel functions} are constructed via automorphic representations $\Theta$ belonging to the automorphic
$L^2$-packet $\wt\Pi_{\psi_0}(\epsilon_{\psi_0})$ associated to the global Arthur parameter
$$
\psi_0=
\begin{cases}
(\tau,b+c)&\text{if}\ m_V=a(b+c);\\
(\tau,b+c)\boxplus(\chi,1)&\text{if}\ m_V=a(b+c)+1,
\end{cases}
$$
where $\tau$ is a conjugate self-dual member in $\CA_\cusp(a)$ and $\chi$ is either the trivial character or $\omega_{E/F}$ the
quadratic character attached to the quadratic field extension $E/F$ by the global classfield theory, depending the parity of the
data.

The construction and the theory in Sections 5.3 and 5.4 work for unitary groups. The details will be presented in a forthcoming work
of Lei Zhang and the author in \cite{JZ13} and are omitted here.

It is worthwhile to mention that the constructions outlined here may work for groups which are not necessarily $F$-quasisplit, which
leads the endoscopy correspondences for inner forms of classical groups (\cite[Chapter 9]{Ar12}). The constructions discussed here
essentially assume that $d=a-1$ and the remaining possible cases are discussed in the forthcoming work of the author (\cite{Jn13}).

\section{Automorphic Descents and Automorphic Forms of Simple Type}

As discussed in previous sections, the main idea to construct explicit endoscopy correspondences for classical groups is to produce
{\bf automorphic kernel functions} $\CK_{\varphi_\Theta}^{\psi_0}(h,g)$ on $H(\BA)\times G(\BA)$ by means of automorphic forms of simple or
simpler type on an ample group $G_m(\BA)$. The integral transforms with the constructed automorphic kernel functions
$\CK_{\varphi_\Theta}^{\psi_0}(h,g)$ produce the endoscopy correspondences: the {\bf endoscopy transfers} from $H(\BA)$ to $G(\BA)$ and
the {\bf endoscopy descents} from $G(\BA)$ to $H(\BA)$. Hence a refined explicit construction theory of endoscopy correspondences
requires refined knowledge about the automorphic $L^2$-packets of simple or simpler type. The automorphic descent method has been developed
by Ginzburg, Rallis and Soudry in 1999 and completed in their book in 2011 (\cite{GRS11}) for $F$-quasisplit classical groups, and
its extension to the spinor groups is given by J. Hundley and E. Sayag in \cite{HS12}. The generalization of the automorphic descent method by
Ginzburg, Soudry and the author (\cite{GJS12}), and through the work of Baiying Liu and the author
(\cite{JL13}) and the work of Lei Zhang and the author (\cite{JZ13})
provides an effective way to establish refined structures of automorphic $L^2$-packets of simple type. The connections with the topics
in the classical theory of automorphic forms: the Saito-Kurokawa conjecture and its relations to the Kohnan space and Maass space of
classical modular forms, with the Ikeda liftings and the
Andrianov conjecture in classical theory of automorphic forms will also be briefly discussed and more work in this aspect
will be addressed in a future work.

\subsection{Automorphic descents}
The method of automorphic descents and its connection with the Langlands functoriality between automorphic representations of classical groups
and general linear groups was discovered by Ginzburg, Rallis and Soudry in 1999, through their serious of papers published since
then. Their book (\cite{GRS11}) gives a complete account of their method for all $F$-quasisplit classical groups. The idea to
construct such automorphic descents may go back to the classical examples in earlier work of H. Maass (\cite{Ms79}),
of A. Andrianov (\cite{An79}), of I. Piatetski-Shapiro
(\cite{PS83}), of M. Eichler and D. Zagier {\cite{EZ85}), and of N. Skoruppa and D. Zagier (\cite{SZ88}), all of which are related to
the understanding of the classical example of non-tempered cuspidal automorphic forms of $\Sp(4)$ by H. Saito and N. Kurokawa in 1977.
The method in \cite{GRS11} extends to the great
generality the argument of Piatetski-Shapiro in \cite{PS83}.

In terms of the general formulation of the constructions of the {\bf automorphic kernel functions} in Section 5, the automorphic
descents of Ginzburg, Rallis and Soudry (\cite{GRS11}) is to deal with the special cases corresponding to the partitions of type
$[d^11^*]$ with $c=1$ and $b=1$. The following discussion will be focused on two cases of automorphic descents and concerned with
their extensions and refined properties as developed in the work of Ginzburg, Soudry and the author (\cite{GJS12}),
the PhD. thesis of Liu (\cite{Liu13}) and a more recent work of Liu and the author (\cite{JL13}).

Recall from Section 5.4 the formulation of the Fourier-Jacobi coefficients of automorphic forms in the special case corresponding to
the partition of type $[(2k)1^{(2m-2k)}]$ with $G_m(\BA)=\Sp_{2m}(\BA)$ or $\tilsp_{2m}(\BA)$, and $c=1$ and $d=2k$. The Fourier-Jacobi
coefficient defines a mapping $\wt{\fj}^{(2m)}_{(2m-2k)}$ from automorphic forms on $\Sp_{2m}(\BA)$ to automorphic forms on $\tilsp_{2m-2k}(\BA)$
and a mapping $\fj^{(2m)}_{(2m-2k)}$ from automorphic forms on $\tilsp_{2m}(\BA)$ to automorphic forms on $\Sp_{2m-2k}(\BA)$. In general,
one does not expect that those mappings carry any functorial meaning for automorphic forms in the sense of Langlands. However,
they do yield the functorial transfers when they are restricted to the automorphic forms of simple type in the sense of Arthur (\cite{Ar12}).

The discussion will be focused on one case when $\tau\in\CA_\cusp(2e)$ has the property that it is of symplectic type and $L(\frac{1}{2},\tau)\neq0$
for the standard $L$-function $L(s,\tau)$. Following \cite{GRS11}, the automorphic descent construction can be reformulated in
terms of the Arthur classification theory of discrete spectrum and described in the following diagram:
\begin{equation}
\begin{matrix}
          &               &                  & \Sp_{4e}(\BA)  & \CE_\tau\in\wt\Pi_{(\tau,2)\boxplus(1,1)}(\epsilon_{(\tau,2)\boxplus(1,1)})\\
                    &                              &{^{RES}\nearrow }  &                             &                   \\
 \tau&  \GL_{2e}(\BA) &                  &\downarrow\vcenter{\rlap{FJ}}             &                   \\
                   &                               &{_{LFT}\nwarrow }  &                             &                    \\
                     &               &                  & \tilsp_{2e}(\BA)      & \tilpi\in\wt\Pi_{(\tau,1)}(\epsilon_{(\tau,1)})\end{matrix}
\end{equation}

This diagram of constructions can be reformulated as
$$
\begin{matrix}
          &               &                  & \Sp_{4e}    & (\tau,2)\boxplus(1,1)\\
                    &                              &\swarrow   &                             &                   \\
(\tau,1) & \tilsp_{2e}   &                  &\uparrow               &                   \\
                   &                               &\searrow   &                             &                    \\
                     &               &                  & \Sp_{0}            & (1,1)\end{matrix}
$$
which can be extended to the general version of the constructions:
\begin{equation}\label{bt1}
\begin{matrix}
          &               &                  & \Sp_{4el+4e}    & (\tau,2l+2)\boxplus(1,1)\\
                    &                              &\swarrow   &                             &                   \\
(\tau,2l+1) & \tilsp_{4el+2e}   &                  &\uparrow               &                   \\
                   &                               &\searrow   &                             &                    \\
                     &               &                  & \Sp_{4el}            & (\tau,2l)\boxplus(1,1)\end{matrix}
\end{equation}
This diagram is called a {\bf basic triangle} of constructions. The dual triangle is the following:
\begin{equation}\label{bt2}
\begin{matrix}
(\tau,2l+1)        & \tilsp_{4el+2e}           &                  &           &                   \\
                    &                               &\searrow   &                            &                    \\
                    &\uparrow                &                  & \Sp_{4el}            & (\tau,2l)\boxplus(1,1)\\
                    &                              &\swarrow   &                             &                   \\
(\tau,2l-1)        & \tilsp_{4el-2l}           &      & & \end{matrix}
\end{equation}
Putting the two basic triangles of constructions, one obtains the following general diagram with $\tau$ symplectic, $L(\frac{1}{2},\tau)\neq0$:
\begin{equation}\label{diag1}
\begin{matrix}
\vdots           &\uparrow                &                  & \Sp_{4el}    & (\tau,2l)\boxplus(1,1)\\
                    &                              &\swarrow   &                             &                   \\
(\tau,2l-1) & \tilsp_{4el-2e}   &                  &\uparrow               &                   \\
                   &                               &\searrow   &                             &                    \\
                     &\uparrow                &                  & \Sp_{4el-4e}            & (\tau,2l-2)\boxplus(1,1)\\
                     &                              &\swarrow   &                             &                   \\
\vdots          &\vdots                    &\vdots        &\vdots                  &  \vdots         \\
                    &\uparrow                &                  & \Sp_{8e}            & (\tau,4)\boxplus(1,1)\\
                    &                              &\swarrow   &                            &                    \\
(\tau,3)        & \tilsp_{6e}           &                  &\uparrow              &                   \\
                    &                               &\searrow   &                            &                    \\
                    &\uparrow                &                  & \Sp_{4e}            & (\tau,2)\boxplus(1,1)\\
                    &                              &\swarrow   &                             &                   \\
(\tau,1)        & \tilsp_{2e}           &      & & \end{matrix}
\end{equation}
The matching of the global Arthur parameters in this diagram under the Fourier-Jacobi mappings has been verified in \cite{JL13}.

It is proved in \cite{GJS12} that those two basic triangles are commuting diagrams for distinguished members in the packets. It is clear that
for distinguished members of the packets, the vertical arrow-up is constructed by taking residue of the Eisenstein series with the given cuspidal
datum. However, if one considers general members of the packets, the vertical arrow-up can be constructed by taking a special case of
the endoscopy transfer considered in Section 5. Of course, the existence of such transfers is confirmed in \cite{Ar12} (and also \cite{Mk12})
as one step of the induction argument of Arthur. The commutativity of those basic triangles in general will be discussed in \cite{JL13}.

\subsection{Finer structures of automorphic $L^2$-packets of simple type}
Diagram \eqref{diag1} and its properties can be used to establish finer structures of automorphic $L^2$-packets of simple type in terms of
of the structure of Fourier coefficients in the spirit of Conjectures \ref{cubmfc} and \ref{cmfc}. This is tested in \cite{Liu13} and \cite{JL13}.

The following basic triangle was first considered in \cite{GJS12}:
\begin{equation}\label{bt3}
\begin{matrix}
&(\tau,3)        & \tilsp_{6e}           &                  &           &                 &  \\
     &               &                               &\searrow   &                            &    &                \\
    &                &\uparrow                &                  & \Sp_{4e}            & (\tau,2)& \\
   &                 &                              &\swarrow   &                             &       &            \\
 &(\tau,1)        & \tilsp_{2e}           &      & & &\end{matrix}
\end{equation}
The motivation was to understand more precise structures of the transfers between $\tilsp_{2e}(\BA)$ and $\Sp_{4e}(\BA)$. From the diagram,
the transfer $\Phi$ from $\Sp_{4e}(\BA)$ to $\tilsp_{2e}(\BA)$ is given by the corresponding Fourier-Jacobi coefficient. On the other
hand, the composition of the transfer from $\tilsp_{2e}(\BA)$ to $\tilsp_{6e}(\BA)$ with the transfer from $\tilsp_{6e}(\BA)$ to
$\Sp_{4e}(\BA)$ given by the corresponding Fourier-Jacobi coefficient yields a transfer $\Psi$ directly from $\tilsp_{2e}(\BA)$ to $\Sp_{4e}(\BA)$.
The question was to characterize the image under the transfer $\Psi$ from $\tilsp_{2e}(\BA)$ to $\Sp_{4e}(\BA)$ of the
automorphic $L^2$-packet $\wt\Pi_{(\tau,1)}(\epsilon_{(\tau,1)})$ in automorphic $L^2$-packet
$\wt\Pi_{(\tau,2)\boxplus(1,1)}(\epsilon_{(\tau,2)\boxplus(1,1)})$. The main result of \cite{GJS12} can be reformulated in terms of the Arthur
classification as follows.

Recall that the Fourier-Jacobi coefficients defined in Section 5.4 depends on an $F$-rational datum $\udl{q}$. In the current case here,
this $\udl{q}$ reduces to the square classes of the number field $F$. For simplicity in notation here, the square classes are denoted by $\alpha$.
This $\alpha$ defines a nontrivial additive character of $F$, denoted by $\psi_\alpha$, i.e. $\psi_\alpha(x):=\psi_F(\alpha x)$. To be more
precise, the global Arthur packet $\wt\Pi_\psi$ for $\tilsp_{2e}(\BA)$ depends on a choice of additive character of $F$. It is reasonable
to denote by $\wt\Pi_\psi^\alpha(\epsilon_\psi)$ the automorphic $L^2$-packet within the global Arthur packet $\wt\Pi_\psi^\alpha$, whose
members are $\psi_\alpha$-functorial transferred to the Arthur representation $\pi_\psi$ of $\GL_{2e}(\BA)$ as discussed in Section 2.2.
When $\psi$ is a generic global Arthur parameter, $\psi_\alpha$-functorial transfer from $\tilsp_{2e}(\BA)$ to $\GL_{2e}(\BA)$ was
established in \cite{GRS11} for the generic members in $\wt\Pi_\psi^\alpha(\epsilon_\psi)$. Of course, the whole theory on the
Arthur classification of the discrete spectrum for $\tilsp_{2e}(\BA)$ is not known yet, and is expected to be achieved through
series of publications of Wen-Wei Li started from his PhD. thesis \cite{Liw11}.

Define $\wt\Pi_{(\tau,2)\boxplus(1,1)}^{(-\alpha)+}(\epsilon_{(\tau,2)\boxplus(1,1)})$ to be the subset of the automorphic $L^2$-packet
$\wt\Pi_{(\tau,2)\boxplus(1,1)}(\epsilon_{(\tau,2)\boxplus(1,1)})$ consisting of all members whose $\psi_{-\alpha}$-Fourier-Jacobi descent
as discussed above from $\wt\Pi_{(\tau,2)\boxplus(1,1)}(\epsilon_{(\tau,2)\boxplus(1,1)})$ to $\wt\Pi_{(\tau,1)}^\alpha(\epsilon_{(\tau,1)})$ are non-zero. Hence one has
the following disjoint union decomposition
\begin{equation}
\wt\Pi_{(\tau,2)\boxplus(1,1)}(\epsilon_{(\tau,2)\boxplus(1,1)})
=
\wt\Pi_{(\tau,2)\boxplus(1,1)}^{(-\alpha)+}(\epsilon_{(\tau,2)\boxplus(1,1)})
\cup
\wt\Pi_{(\tau,2)\boxplus(1,1)}^{(-\alpha)-}(\epsilon_{(\tau,2)\boxplus(1,1)})
\end{equation}
where $\wt\Pi_{(\tau,2)\boxplus(1,1)}^{(-\alpha)-}(\epsilon_{(\tau,2)\boxplus(1,1)})$ is defined to be the complimentary subset of
the subset $\wt\Pi_{(\tau,2)\boxplus(1,1)}^{(-\alpha)+}(\epsilon_{(\tau,2)\boxplus(1,1)})$ in the automorphic $L^2$-packet
$\wt\Pi_{(\tau,2)\boxplus(1,1)}(\epsilon_{(\tau,2)\boxplus(1,1)})$. The main result of \cite{GJS12} can be reformulated as follows in terms
of the Arthur classification theory.

\begin{thm}[\cite{GJS12}]\label{gjsijm}
With notation as above, the transfers $\Phi^\alpha$ and $\Psi^\alpha$ establish one-to-one correspondence between the set
$\wt\Pi_{(\tau,1)}^\alpha(\epsilon_{(\tau,1)})$ and the set $\wt\Pi_{(\tau,2)\boxplus(1,1)}^{(-\alpha)+}(\epsilon_{(\tau,2)\boxplus(1,1)})$.
\end{thm}

In order to prove this theorem, one has to combine the basic triangle \eqref{bt3} with the basic triangle \eqref{bt1} with $l=1$.
This leads to the following diagram of transfers:
\begin{equation}
\begin{matrix}
                    &              &                  & \Sp_{8e}            & (\tau,4)\boxplus(1,1)\\
                    &                              &\swarrow   &                            &                    \\
(\tau,3)        & \tilsp_{6e}           &                  &\uparrow              &                   \\
                    &                               &\searrow   &                            &                    \\
                    &\uparrow                &                  & \Sp_{4e}            & (\tau,2)\boxplus(1,1)\\
                    &                              &\swarrow   &                             &                   \\
(\tau,1)        & \tilsp_{2e}           &      & & \end{matrix}
\end{equation}
The commutativity of the transfers in this diagram plays essential role in the proof of Theorem \ref{gjsijm}.

It is worthwhile to make some remarks on the relation of this theorem to previous known works in automorphic representation theory and
in classical theory of automorphic forms.
First of all, this theorem can be regarded as a refinement of the main results of \cite{GRS05}, where the authors tried to construct,
with a condition of the structure of certain Fourier coefficients,
some cuspidal automorphic representations of $\Sp_{2n}(\BA)$ or $\tilsp_{2n}(\BA)$ that are expected to belong to one of the global Arthur
packets of simple type as discussed above. Those conditions will be discussed below related to a work of Liu and the author \cite{JL13}.

Secondly, when $e=1$, this theorem is a refinement of the main result of Piatetski-Shapiro's pioneer paper (\cite{PS83}) on the
representation-theoretic approach to the Saito-Kurokawa examples of non-tempered cuspidal automorphic forms on $\Sp_4$; and can also be
viewed as the representation-theoretic version of a theorem of Skoruppa and Zagier (\cite{SZ88}), which completed the earlier work of
H. Maass (\cite{Ms79}) and of A. Andrianov (\cite{An79}). For a complete account of this classical topic, the readers are referred to
the book (\cite{EZ85}) and the survey paper (\cite{Z81}).

Finally, it will be an interesting problem to work out the explicit version in the classical theory of automorphic forms of this theorem,
which will be the natural extension of the Skoruppa and Zagier theorem (\cite{SZ88}) to classical automorphic forms of higher genus.

One of the motivations of the work of Liu and the author (\cite{Liu13} and \cite{JL13})
is to extend Theorem \ref{gjsijm} to the cases with general global
Arthur parameters of simple type for $\Sp_{4el}(\BA)$ and $\tilsp_{4el+2e}(\BA)$, respectively.
More precisely, one defines similarly the following decomposition
\begin{equation}\label{alphapmsp}
\wt\Pi_{\psi}(\epsilon_{\psi})
=
\wt\Pi_{\psi}^{(-\alpha)+}(\epsilon_{\psi})
\cup
\wt\Pi_{\psi}^{(-\alpha)-}(\epsilon_{\psi})
\end{equation}
for $\Sp_{4el}(\BA)$ with $\psi=(\tau,2l)\boxplus(1,1)$, where
$\wt\Pi_{\psi}^{(-\alpha)+}(\epsilon_{\psi})$ consists of elements of
$\wt\Pi_{\psi}(\epsilon_{\psi})$ with a non-zero $\psi_{-\alpha}$-Fourier-Jacobi descent from
$\Sp_{4el}(\BA)$ to $\tilsp_{4el-2e}(\BA)$, and $\wt\Pi_{\psi}^{(-\alpha)-}(\epsilon_{\psi})$
is the complement of $\wt\Pi_{\psi}^{(-\alpha)+}(\epsilon_{\psi})$ in
$\wt\Pi_{\psi}(\epsilon_{\psi})$.

Similarly, one defines
\begin{equation}\label{alphapmspt}
\wt\Pi_{(\tau,2l+1)}(\epsilon_{(\tau,2l+1)})
=
\wt\Pi_{(\tau,2l+1)}^{\alpha+}(\epsilon_{(\tau,2l+1)})
\cup
\wt\Pi_{(\tau,2l+1)}^{\alpha-}(\epsilon_{(\tau,2l+1)})
\end{equation}
for $\tilsp_{4el+2e}(\BA)$, where the subset
$\wt\Pi_{(\tau,2l+1)}^{\alpha+}(\epsilon_{(\tau,2l+1)})$ consists of elements in
$\wt\Pi_{(\tau,2l+1)}(\epsilon_{(\tau,2l+1)})$ with a non-zero $\psi_{\alpha}$-Fourier-Jacobi descent from
$\tilsp_{4el+2e}(\BA)$ to $\Sp_{4el}(\BA)$, and the subset $\wt\Pi_{(\tau,2l+1)}^{\alpha-}(\epsilon_{(\tau,2l+1)})$ is defined to be its
complement.

It is clear that further understanding of the structures of those sets and their mutual relations are
the key steps to understand the structures of Fourier coefficients of the automorphic representations within
$\wt\Pi_{(\tau,2l)\boxplus(1,1)}(\epsilon_{(\tau,2l)\boxplus(1,1)})$ and $\wt\Pi_{(\tau,2l+1)}(\epsilon_{(\tau,2l+1)})$, respectively.
In \cite{Liu13}, the cuspidal part of those sets has been more carefully studied, which is
an extension of the main result of \cite{GJS12}.

The problems have been formulated for other classical groups according to the functorial type of $\tau$. The details of the progress in
this aspect will be considered elsewhere. Finally, it is worthwhile to remark briefly that the theory developed here
will yields explicit constructions for the Duke-Imamoglu-Ikeda lifting from certain cuspidal automorphic forms on $\wt\SL_2(\BA)$ to
cuspidal automorphic forms on $\Sp_{4l}(\BA)$, the construction via explicit Fourier expansion of this transfer was obtain by
T. Ikeda in \cite{Ik01}. If combined with the constructions of endoscopy correspondences discussed in Section 5,
this method is able to construct the conjectural lifting of Andrianov (\cite{An01}). See also the work of Ikeda (\cite{Ik06} and \cite{Ik08})
for relevant problems.

\section{Theta Correspondence and $(\chi,b)$-Theory}

The theta correspondence for reductive dual pairs is a method, introduced by R. Howe 1979 (\cite{Hw79}),
to construct automorphic forms using the theta functions built from the Weil representation on the Schr\"odinger model. The method to
use classical theta functions to construct automorphic forms goes back at least to the work of C. Siegel (\cite{Sg66}) and
the work of A. Weil (\cite{Wl64} and \cite{Wl65}).
It is the fundamental work of Howe (\cite{Hw79}) which started the representation-theoretic approach in
this method guided by the idea from Invariant Theory.

The theta correspondence gave the representation-theoretic framework of the classical theory of modular forms of half-integral weight and
its relation to the the theory of modular forms of integral weight, i.e. the classical Shimura correspondence (\cite{Sm73}). A series of papers of
J.-L. Waldspurger (\cite{W80}, \cite{W84} and \cite{W91}) gave a complete account of this theory, which was beautifully summarized by
Piatetski-Shapiro in \cite{PS84}.

One of the early contributions of this method to the theory of automorphic forms was the discovery of families of the counter-examples of
the generalized Ramanujan conjecture for algebraic groups different from the general linear groups. The first of such examples was
constructed by Howe and Piatetski-Shapiro (\cite{HPS79}) and then more examples of similar kind were constructed by Piatetski-Shapiro in 1980's
on $\GSp(4)$ including the example of Saito-Kurokawa and by many other people for more general groups. The existence of those non-tempered cuspidal automorphic representations enriched
the spectral theory of automorphic forms for general reductive algebraic groups (\cite{Ar84} and \cite{Ar04}).

In order to understand the spectral structures of such striking examples constructed by the theta correspondence, S. Rallis launched,
after his series of work joint with G. Schiffmann (\cite{RS77}, \cite{RS78} and \cite{RS81}) and the work of Waldspurger cited above,
a systematic investigation (sometimes called the {\bf Rallis program} which was summarized later in (\cite{Rl91})) of the spectral structure of the theta correspondences in family
(the so-called {\bf Rallis Tower Property} (\cite{Rl84})), the relation between theta correspondence and the Langlands functoriality (\cite{Rl84}),
and the characterization of the non-vanishing property of such constructions (\cite{Rl87}).
The Rallis program leads to a series of the fundamental contributions in the modern theory of automorphic forms:
the {\bf Rallis Inner Product Formula} (\cite{Rl87}), the {\bf doubling method} for the standard $L$-functions of classical groups by Piatetski-Shapiro and Rallis (\cite{GPSR87}), and the {\bf Regularized Siegel-Weil formula} of S. Kudla and S. Rallis (\cite{KR88}, \cite{KR94}
and \cite{Kd08}).

The regularized Siegel-Weil formula has been extended to all classical groups through the work of T. Ikeda (\cite{Ik96}),
C. Moeglin (\cite{Mg97}); A. Ichino (\cite{Ic01}, \cite{Ic04} and \cite{Ic07});  D. Soudry and
the author (\cite{JS07}); W.-T. Gan and S. Takeda (\cite{GT11}; S, Yamana
(\cite{Ym11} and \cite{Ym12}); and W.-T. Gan, Y.-N. Qiu, and S. Takeda (\cite{GQT12}.  Furthermore, the work of W.-T. Gan extends
the formula to some exceptional groups (\cite{Gn00}, \cite{Gn08} and \cite{Gn11}).

Applications of the Rallis program to number theory and
arithmetic are broad and deep, including arithmeticity of special values of certain automorphic $L$-functions
(\cite{HLSk05}, \cite{HLSk06}, and \cite{Pr09} for instance);
the nonvanishing of cohomology groups of certain degree over Shimura varieties (\cite{Li92} and \cite{BMM11}); the
Kudla program on special cycles and a generalized Gross-Zagier formula (\cite{Kd02} and \cite{YZZ12}), and see also relevant work (\cite{Hr94}, \cite{Hr07},
\cite{Hr08}, and \cite{HLSn11}), to mention a few.

\subsection{Rallis program}
The main idea and method in the Rallis program are explained here only with
one case of reductive dual pairs, which is completely the same for all other reductive dual pairs. Take $(\Sp_{2n},\SO_{m})$, which
forms a reductive dual pair in $\Sp_{2mn}$. The theta function $\Theta_\phi^{\psi_\alpha}(x)$ attached to a global Bruhat-Schwartz function
$\phi$ in the Schr\"odinger module $\CS(\BA^{mn})$ of the Weil representation associated to the additive character $\psi_\alpha$ is an automorphic
function of moderate growth on $\tilsp_{2mn}(\BA)$. When restricted to the subgroup $\tilsp_{2n}(\BA)\times\wt{\SO}_{m}(\BA)$,
the metaplectic double splits. For $\sigma\in\CA(\SO_{m})$ and $\pi\in\CA(\Sp_{2n})$, if the following integral
\begin{equation}\label{tc1}
\int_{\Sp_{2n}(F)\bks\Sp_{2n}(\BA)}\int_{\SO_{2m}(F)\bks\SO_{2m}(\BA)}
\Theta_\phi^{\psi_\alpha}(g,h)\varphi_\pi(g)\ovl{\varphi_\sigma(h)}dgdh
\end{equation}
converges and is non-zero, then this family of the integrals with all $\phi\in\CS(\BA^{mn})$ yields a correspondence between automorphic
forms on $\Sp_{2n}(\BA)$ and automorphic forms on $\SO_{m}(\BA)$. For instance, letting the integer $m$ run, Rallis got the following
diagram (called the Witt tower) of theta correspondences:
\begin{equation}\label{ttw}
\begin{matrix}
                              &                        &                                    &  \vdots                         &        \\
                              &                        &                                    &                                     & \\
                              &                        &                                   & \SO_{m}  & \\
                              &                        &                                    &  \vdots                         &        \\
                              &                         &                                    &                                     & \\
                              &                        &      \nearrow               &  \uparrow                     &                   \\
                              &                         &                                    &                                     & \\
  &  \Sp_{2n} &        \longrightarrow    &\SO_{4}          &\\
                              &                         &                                    &                                     & \\
                              &                       &       \searrow                & \uparrow                     &                   \\
                              &                         &                                    &                                     & \\
                              &                        &                                   & \SO_{2}   &  \end{matrix}
\end{equation}

The Rallis tower property of theta correspondence asserts that the first occurrence of the theta correspondence of cuspidal automorphic
representation $\pi$ of $\Sp_{2n}(\BA)$ in this tower is always cuspidal.
Here the first occurrence means that for a given automorphic representation $\pi$ of $\Sp_{2n}(\BA)$,
the smallest integer $m$ such that the integral \eqref{tc1} is non-zero for some automorphic representation $\sigma$ of $\SO_{m}(\BA)$ with
a certain choice of the Bruhat-Schwartz function $\phi$. Hence the spectral theory of the theta correspondence is to ask the following
two basic questions:
\begin{enumerate}
\item {\bf First Occurrence Problem:} For a given irreducible cuspidal automorphic representation $\pi$ of $\Sp_{2n}(\BA)$, how
to determine the first occurrence of the theta correspondence in the Witt tower in terms of the basic invariants attached to $\pi$?
\item {\bf Spectral Property of the First Occurrence:} For a given irreducible cuspidal automorphic representation $\pi$ of $\Sp_{2n}(\BA)$,
what kind of cuspidal automorphic representations are expected to be constructed at the first occurrence?
\end{enumerate}

For the {\bf First Occurrence Problem}, Rallis calculated the $L^2$-inner product of the image of the theta correspondence, which leads to
the well-known {\bf Rallis Inner Product Formula}. This formula expresses the $L^2$-inner product of the image of the theta correspondence
in terms of a basic invariant of the original given $\pi$, that is, the special value or residues of the standard automorphic $L$-function
of $\pi$, up to certain normalization. The whole program has been completely carried out through the work by many people, notably,
Rallis (\cite{Rl87} and \cite{Rl91}); Kudla and Rallis (\cite{KR94}); Moeglin (\cite{Mg97}); Lapid and Rallis (\cite{LR05});
Soudry and the author (\cite{JS07}); Ginzburg, Soudry and the author (\cite{GJS09}); Gan and Takeda (\cite{GT11}); Yamana (\cite{Ym12}); and
Gan, Qiu and Takeda (\cite{GQT12}).

The {\bf Spectral Property of the First Occurrence} is mainly to understand the relation of the transfer of automorphic representations
given by theta correspondence with the Langlands functoriality.

As suggested by the global Howe duality conjecture (\cite{Hw79}), the irreducibility of the first occurrence
of the theta correspondences was proved by Moeglin (\cite{Mg97}) for even orthogonal group case, by Soudry and the author (\cite{JS07})
for the odd orthogonal group case, and by Chenyan Wu (a current postdoc of the author) in \cite{Wu12} for the unitary group case.

The work of Rallis in 1984 (\cite{Rl84}) started the investigation of the relation between the theta correspondence in family and
the local Langlands parameters at unramified local places (see also the work of Kudla \cite{Kd86}). It seems that the unramified local Langlands
parameters can not be well preserved through the theta correspondence in family. It was J. Adams who first observed that
the theta correspondence may transfer the Arthur parameters in a way compatible with the Arthur conjecture. This is known as the Adams
conjecture (\cite{Ad89}), which was discussed explicitly for unitary group theta correspondences by Harris, Kudla and Sweet in
 \cite{HKS96}. More recent progress in this aspect is referred to the work of Moeglin (\cite{Mg11}). The nature of theta correspondence and
 the Adams conjecture suggest that the theory of theta correspondence and related studies in automorphic forms are to detect the
 global Arthur parameter of type $(\chi,b)$ with a quadratic character $\chi$ of the relevant automorphic representations, and hence
 it seems reasonable to call such studies the $(\chi,b)$-theory of automorphic representations.

\subsection{$(\chi,b)$-theory}
Let $\pi$ be an irreducible cuspidal automorphic representation of $\Sp_{2n}(\BA)$ and
let $\chi$ be a quadratic character of $F^\times\bks\BA^\times$. The (twisted) standard (partial) $L$-function $L^S(s,\pi\times\chi)$
was studied through the doubling method of Piatetski-Shapiro and Rallis (\cite{GPSR87}). A combination of the work of Lapid and Rallis (\cite{LR05})
and Yamana (\cite{Ym12}) completes the basic theory of this automorphic $L$-function.

Following Kudla and Rallis (\cite[Theorem 7.2.5]{KR94}),
This partial $L$-function $L^S(s,\pi\times\chi)$ may have a simple pole at
$$
s_0\in\{1,2,\cdots,[\frac{n}{2}]+1\}.
$$
If $s=s_0$ is a right-most pole of $L^S(s,\pi\times\chi)$, then $\pi$ has the first occurrence $\sigma$ at $\SO_{2n-2s_0+2}$, which
is an irreducible cuspidal automorphic representation of $\SO_{2n-2s_0+2}(\BA)$. According to the local calculation (\cite{Kd86} and
\cite{HKS96}}, the Adams conjecture (\cite{Ad89} and \cite{Mg11}), and the relevant global results (\cite{GRS97},
\cite{JS07}, \cite{GJS09} and \cite[(2)]{GJS11}),
it is expected that the global Arthur parameter of $\pi$ contains $(\chi,b)$ as a simple formal summand for some quadratic character $\chi$,
which should be decided through the choice of the additive character for the theta correspondence. In principle, one needs also to consider
the Witt tower of nonsplit orthogonal groups in order to detect more precisely the simple Arthur parameter $(\chi,b)$ in the
global Arthur parameter of $\pi$, as suggested by a conjecture of Ginzburg and Gurevich in \cite{GG06}.

When $L^S(s,\pi\times\chi)$ is holomorphic for $\Re(s)\geq 1$, the first occurrence of $\pi$ at $\SO_{m}$ with $2n+2\leq m\leq 4n$ is
detected by the nonvanishing of $L(s,\pi\times\chi)$ at $s=\frac{m-2n}{2}$, in addition to the local occurrence conditions following
the recent work of Gan, Qiu and Takeda (\cite{GQT12}). If $L^S(s,\pi\times\chi)$ is holomorphic for $\Re(s)\geq 1$ for all quadratic
characters $\chi$, then the global Arthur parameter of $\pi$ contains no formal summand of simple parameters of type $(\chi,b)$.
In this case, the relation between $(\chi,b)$ and the global Arthur parameter of $\pi$ should be detected through a version of
the Gan-Gross-Prasad conjecture (\cite{GGP12}). This also suggests that
one needs to consider a theory to detect the simple parameter $(\tau,b)$ for $\pi$ with $\tau\in\CA_\cusp(a)$ and $a>1$.
This is exactly the reason to suggest the $(\tau,b)$-theory of automorphic representations, which will be discussed below with more
details.

\section{Endoscopy Correspondence and $(\tau,b)$-Theory}

As indicated in the $(\chi,b)$-theory discussed above, for a general $F$-quasisplit classical group $G_n$, the relation between the structure
of global Arthur parameters and the structure of the poles of the partial tensor product $L$-functions $L^S(s,\pi\times\tau)$ can be described as follows, where $\pi\in\CA_\cusp(G_n)$ and $\tau\in\CA_\cusp(a)$. This is the key reason to suggest the $(\tau,b)$-theory for
automorphic forms of $F$-quasisplit classical groups.

\subsection{Global Arthur parameters and poles of partial tensor product $L$-functions}
Let $G$ be a classical group belonging to $\wt{\CE}_\simp(N)$.  Each $\psi\in\wt{\Psi}_2(G)$ can be written as
$$
\psi:=\psi_1\boxplus\psi_2\boxplus\cdots\boxplus\psi_r
$$
where each simple global Arthur parameter $\psi_i=(\tau_i,b_b)$ with $\tau_i\in\CA_\cusp(a_i)$ being self-dual and $b_i\geq 1$ being integers
that satisfy certain parity condition with $\tau_i$, for $i=1,2,\cdots,r$. Now, consider the partial tensor product $L$-functions
$L^S(s,\pi\times\tau)$ for some self-dual member in $\CA_\cusp(a)$. According to the Arthur classification of the discrete spectrum of $G$,
this family of partial $L$-functions can be defined and has meromorphic continuation to the whole complex plane $\BC$.

When $\psi$ is a generic parameter, i.e. $b_i=1$ for all $i$, then those partial $L$-functions can be completed through
the Arthur-Langlands transfer, that is, the ramified local $L$-functions can also be defined (\cite{Ar12}).
However, when the global parameter $\psi$ is not generic, then
it is still a problem to define the ramified local $L$-functions for this family of global $L$-functions.
On the other hand, the Rankin-Selberg method is expected to define the ramified local $L$-functions for this family of the global $L$-functions.
The recent work of Zhang and the author extends the work of Ginzburg, Piatetski-Shapiro and Rallis (\cite{GPSR97}) from orthogonal groups
to the classical groups of hermitian type (\cite{JZ13}) and the skew-hermitian case is also considered by Zhang and the author,
including the work of Xin Shen (\cite{Sn12} and \cite{Sn13}),
extending the work of Ginzburg, Rallis, Soudry and the author on the symplectic groups and
metaplectic groups (\cite{GJRS11}).

In the following discussion, the partial tensor product $L$-functions are used as invariants to make connections to the global Arthur parameters.
For any $\pi\in\wt{\Pi}_\psi(\epsilon_\psi)$, the automorphic $L^2$-packet associated to the global Arthur parameter $\psi$,
define the set:
\begin{equation}\label{taui}
\frak{T}_\pi:=\{\tau\in\cup_{a\geq 1}\wt\CA_\cusp(a)\ |\ L^S(s,\pi\times\tau)\ \text{has a pole at}\ s\geq\frac{1}{2}\},
\end{equation}
where $\wt\CA_\cusp(a)$ denotes the subset of self-dual members in $\CA_\cusp(a)$. It is a non-empty finite set by the Arthur classification
theory. Write
$$
\frak{T}_\pi=\{\tau_1,\tau_2,\cdots,\tau_t\}
$$
with $\tau_i\in\wt\CA_\cusp(a_i)$ and $\tau_i\not\cong\tau_j$ if $i\neq j$.

For $\tau\in\frak{T}_\pi$ and $b\geq 1$ integer, define $[\tau,b]_\pi$, which may be called a Jordan block of $\pi$ following Moeglin,
to be a pair attached to $\pi$ if the following properties hold:
\begin{enumerate}
\item the partial tensor product $L$-function $L^S(s,\pi\times\tau)$ is holomorphic for the real part of $s$ greater than $\frac{b+1}{2}$; and
\item $L^S(s,\pi\times\tau)$ has a simple pole at $s=\frac{b+1}{2}$.
\end{enumerate}
Note that the right-most pole of the partial tensor product $L$-function $L^S(s,\pi\times\tau)$ is at most simple. This
follows from the structure of global Arthur parameters in $\wt{\Psi}_2(G)$ or the Rankin-Selberg integral method (see \cite{JZ12} for instance).

For each $i$, take $\pi_{(\tau_i,b_i)}$ to be the Arthur representation of $\GL_{a_ib_i}(\BA)$, which are realized as the Speh residues as
automorphic representations. It is an easy exercise to check that $N=\sum_{i=1}^ta_ib_i$ if and only if the following quotient of
partial $L$-functions
\begin{equation}\label{ql}
\frac{L^S(s,\pi\times\tau_i)}{L^S(s,\pi_{(\tau_i,b_i)}\times\tau_i)}
\end{equation}
is holomorphic for $\Re(s)\geq\frac{1}{2}$ for all $i$. If there is an $i$ such that the quotient \eqref{ql} has a pole at
$s=\frac{b_{i,1}+1}{2}$ for some $b_{i,1}\geq 1$ and is holomorphic for $\Re(s)>\frac{b_{i,1}+1}{2}$. One can define $[\tau_i,b_{i,1}]$.
By the structure of global Arthur parameter $\psi$, $b_{i,0}:=b_i$ and $b_{i,1}$ have the same parity. Repeating this argument, one obtains
the following sequence of $\tau_i$-types of $\pi$:
\begin{equation}\label{tautype}
[\tau_i,b_{i,0}]_\pi, [\tau_i,b_{i,1}]_\pi,\cdots,[\tau_i,b_{i,r_i}]_\pi,
\end{equation}
with the property that $b_{i,0}>b_{i,1}>\cdots>b_{i,r_i}>0$ and the quotient
$$
\frac{L^S(s,\pi\times\tau_i)}{\prod_{j=0}^{r_i}L^S(s,\pi_{(\tau_i,b_{i,j})}\times\tau_i)}
$$
is holomorphic for $\Re(s)\geq\frac{1}{2}$. Finally the originally given global Arthur parameter $\psi$ can be expressed as
$$
\boxplus_{i=1}^t\boxplus_{j=0}^{r_i}(\tau_i,b_{i,j}).
$$

At this point, the $(\tau,b)$-theory is to detect, for a given $\pi\in\wt{\Pi}_\psi(\epsilon_\psi)$, the pairs
\begin{equation}\label{mtautypes}
[\tau_1,b_{i,0}]_\pi, [\tau_2,b_{2,0}]_\pi,\cdots,[\tau_t,b_{t,0}]_\pi,
\end{equation}
by the constructions provided in previous sections. In other words, it is to detect the "maximal" simple factors of $\psi$ in the sense that
a simple factor $(\tau,b)$ of $\psi$ is maximal if
$$
b=\max\{\beta\in\BZ\ |\ (\tau,\beta)\ \text{is a simple summand of }\ \psi\}.
$$
From the discussion below, it will be much involved in dealing with the simple factors of $\psi$, which are not maximal.

\subsection{$(\tau,b)$-towers of endoscopy correspondences}
From the constructions discussed in Section 5, if let the integer $b$ vary, one gets the following diagrams of constructions of
endoscopy correspondences for classical groups.

Consider the case of $a=2e$ and $2n+1=ab+c$ with $c=2f+1$. When $\tau$ is of symplectic type, take $b=2l+1$ to be odd. Conjecture 5.4 can
be more specified at this case by the following diagram of endoscopy correspondences:
\begin{equation}\label{tbtw1}
\begin{matrix}
                       &           &               &  \vdots           & \vdots        \\
                       &           &               &                   & \\
                       &           &               & \SO_{c+3a}&\psi_{\SO_{c}}\boxplus(\tau,3)\\
                       &           &               &                   & \\
                       &           &\nearrow       &  \uparrow         &                   \\
                       &           &               &                   & \\
 \psi_{\SO_{c}}&\SO_{c} &\longrightarrow&\SO_{c+a}     &\psi_{\SO_{c}}\boxplus(\tau,1)\\
                       &           &               &                   & \\
                       &           &\searrow       & \uparrow          &                   \\
                       &           &               &                   & \\
                       &           &               &\SO_{c-a} & \psi_{\SO_{c}}\boxminus(\tau,1)            \\
                       &           &               &                   & \\
                       &           &               &    \vdots                      & \vdots
                                      \end{matrix}
\end{equation}
Here the parameter $\psi_{\SO_{c}}\boxminus(\tau,1)$ is understood in the sense of the Grothendieck group generated by the global Arthur parameters,
which also gives the meaning for general $\psi_1\boxminus\psi_2$.
When $\tau$ is of orthogonal type, take $b=2l$ to be even. Conjecture 5.4 can be more specified by the following diagram
of endoscopy correspondences:
\begin{equation}\label{tbtw2}
\begin{matrix}
                       &           &               &  \vdots           & \vdots        \\
                       &           &               &                   & \\
                       &           &               & \SO_{c+2a}&\psi_{\SO_{c}}\boxplus(\tau,2)\\
                       &           &               &                   & \\
                       &           &\nearrow       &  \uparrow         &                   \\
                       &           &               &                   & \\
 \psi_{\SO_{c}}&\SO_{c} &\longrightarrow&\SO_{c}     &\psi_{\SO_{c}}\boxplus(\tau,1)\\
                       &           &               &                   & \\
                       &           &\searrow       & \uparrow          &                   \\
                       &           &               &                   & \\
                       &           &               &\SO_{c-2a} & \psi_{\SO_{c}}\boxminus(\tau,2)            \\
                       &           &               &                   & \\
                       &           &               &    \vdots                      & \vdots
                                      \end{matrix}
\end{equation}

The constructions discussed in Section 5.3 will produce similar diagrams of endoscopy correspondences for other cases as above.
The details will be omitted here.

It is important to note that there are also two types of {\bf basic triangles} of constructions
here in Diagrams \eqref{tbtw1} and \eqref{tbtw2}, which is a generalization of the two types of the basic triangles of constructions discussed in
Section 6. The properties of these basic triangles of constructions are the key to establish the
functorial transfer properties of the general diagrams of constructions in \eqref{tbtw1} and \eqref{tbtw2}.

Next, consider the constructions discussed in Section 5.4.

Take $a=2e+1$ with $c=2f+1$ and $2n=ab+c-1$. In this case, $\tau$ must be of orthogonal type and $b=2l$ must be even.
A combination of the constructions discussed in Conjectures 5.10 and 5.12 can be more specified by the following diagram
of endoscopy correspondences:
\begin{equation}\label{tbtw3}
\begin{matrix}
                       &           &               &  \vdots           & \vdots        \\
                       &           &               &                   & \\
                       &           &               & \tilsp_{2f+2a}&\psi_{\SO_{2f+1}}\boxplus(\tau,2)\\
                       &           &               &                   & \\
                       &           &\nearrow       &  \uparrow         &                   \\
                       &           &               &                   & \\
 \psi_{\SO_{2f+1}}&\SO_{2f+1} &\longrightarrow&\tilsp_{2f}     &\psi_{\SO_{2f+1}}\boxplus(\tau,0)\\
                       &           &               &                   & \\
                       &           &\searrow       & \uparrow          &                   \\
                       &           &               &                   & \\
                       &           &               &\tilsp_{2f-2a} &\psi_{\SO_{2f+1}}\boxminus(\tau,2)             \\
                       &           &               &                   & \\
                       &           &               &    \vdots                      & \vdots
                                      \end{matrix}
\end{equation}
Of course, a different combination of
the constructions discussed in Conjectures 5.10 and 5.12 can be more specified by the following diagram
of endoscopy correspondences:
\begin{equation}\label{tbtw4}
\begin{matrix}
                       &           &               &  \vdots           & \vdots        \\
                       &           &               &                   & \\
                       &           &               & \SO_{2f+1+2a}&\psi_{\tilsp_{2f}}\boxplus(\tau,2)\\
                       &           &               &                   & \\
                       &           &\nearrow       &  \uparrow         &                   \\
                       &           &               &                   & \\
 \psi_{\tilsp_{2f}}&\tilsp_{2f} &\longrightarrow&\SO_{2f+1}     &\psi_{\tilsp_{2f}}\boxplus(\tau,0)\\
                       &           &               &                   & \\
                       &           &\searrow       & \uparrow          &                   \\
                       &           &               &                   & \\
                       &           &               &\SO_{2f+1-2a} &\psi_{\tilsp_{2f}}\boxminus(\tau,2)             \\
                       &           &               &                   & \\
                       &           &               &    \vdots                      & \vdots
                                      \end{matrix}
\end{equation}
Note that if take $a=1$ and $(\tau,b)=(\chi,b)$, then the diagram \eqref{tbtw3} gives one of the Witt towers for the theta correspondences.
Since the {\bf automorphic kernel functions} used for the endoscopy correspondences are defined in a different way, it is necessary to
check that both the endoscopy correspondences constructed here and the theta correspondences do produce the same correspondences when
the towers coincide. This is in fact a way to verify the conjecture of Adams on the compatibility of theta correspondences and
the Arthur parameters as discussed in the $(\chi,b)$-theory in Section 7.

Similarly, various combinations of the constructions in Conjectures 5.9 and 5.11 produce diagrams of endoscopy correspondences which may
also specialized to the theta correspondences for reductive dual pairs $(\SO_{2f},\Sp_{2n})$ and pairs of unitary groups.

\subsection{$(\tau,b)$-theory of automorphic forms}
The $(\tau,b)$-theory of automorphic forms is to establish the spectral property of the diagrams of endoscopy correspondences as
discussed in the previous section and verify the conjectures stated in Sections 5.3 and 5.4; and detect the structure of
global Arthur parameters for a given irreducible cuspidal automorphic representation $\pi$ of an $F$-quasisplit classical group in terms of the
basic invariants attached to $\pi$. In short, this is to establish a version of the Rallis program in this generality and look for
the analogy of the Rallis inner product formula for the endoscopy correspondences and for a replacement of
the Siegel-Weil formula for this general setting. An important application of those constructions has been found in a recent work
of Lapid and Mao (\cite{LM12}, which can be viewed as part of the $(\tau,b)$-theory.

\section{Final Remarks}

There are a few closely related topics, which will be briefly discussed here and the more details of which may be considered elsewhere.

\subsection{On generic global Arthur parameters}
A global Arthur parameter $\psi\in\Psi_2(G)$ is called {\it generic} if it has form
$$
\psi=(\tau_1,1)\boxplus(\tau_2,1)\boxplus\cdots\boxplus(\tau_r,1)
$$
with $\tau_i\in\CA_\cusp(a_i)$ being conjugate self-dual. The special cases where $b=1$ and $\psi_2$ is generic of conjectures in Sections
5.3 and 5.4 yield the constructions of endoscopy correspondences for the decomposition of the generic global Arthur parameter
$
\psi=(\tau,1)\boxplus\psi_2.
$
For $G\in\wt\CE_\simp(N)$, if one only concerns the endoscopy transfers for irreducible generic cuspidal automorphic representations of
$G(\BA)$, there is an indirect way to prove the existence of such a transfer, by combining the general constructions of automorphic descents
(\cite{GRS11}) with the Langlands functorial transfer from the classical group $G$ to the general linear group $\GL(N)$
for irreducible generic cuspidal automorphic representations of $G(\BA)$ via the method of Converse Theorem (\cite{CKPSS04}
and \cite{CPSS11}),
as explained in the book (\cite{GRS11}). In 2008, Ginzburg in \cite{G08} gave a framework to constructing endoscopy transfers
for generic cuspidal automorphic representations for $F$-split classical groups, and sketched proofs for some cases, which may not be
applied to all the cases because of technical reasons. The constructions described in this paper extends the framework of \cite{G08} to
cover all possible cases. The on-going work of Lei Zhang and the author is dealing with the case of $F$-quasisplit classical groups of
hermitian type (\cite{JZ13}).

For irreducible generic cuspidal automorphic representations $\pi$ of $G(\BA)$, it is a basic question to know the structure of their
global Arthur parameters in terms of the basic invariants attached to $\pi$. It is clear that one needs invariants in addition to the
Fourier coefficients. One choice is to use the set $\frak{T}_\pi$ as introduced in (8.1), which is to use the poles of the
tensor product $L$-functions $L^S(s,\pi\times\tau)$ with $\tau\in\CA_\cusp(a)$. Note that the poles occur at $s=1$ since $\pi$ is generic.

From the point of view of the Langlands Beyond Endoscopy (\cite{L04}), it is better to have a characterization of the endoscopy
structure in terms of the order of the pole at $s=1$ of a family of automorphic $L$-functions
$L^S(s,\pi,\rho)$ for some complex representations $\rho$ of the Langlands dual group $^LG$.
Such a characterization is fully discussed by the author in \cite{Jn06}
for irreducible generic cuspidal automorphic representation of $G=\SO_{2n+1}$ with $\rho$ going through all the fundamental
representations of the dual group $\Sp_{2n}(\BC)$. The results can be naturally extended to all classical groups.

Another characterization of the endoscopy structure in terms of periods of automorphic forms is formulated as
conjectures in a recent work of Ginzburg and the author in \cite{GJ12} for $G=\SO_{2n+1}$, where some lower rank cases have been verified.
The formulation of the conjectures works for all classical groups. The details are omitted here.

\subsection{On the local theory}
Just as in the theory of theta correspondences, it is interesting to consider the local theory of theta correspondences and
the local-global relations of the theta correspondences for automorphic forms. For the $(\tau,b)$-theory, one may consider the
local $(\tau,b)$-theory with $\tau$ an irreducible, essentially tempered representation of $\GL_a$ over a local field. This local
theory should include the constructions of local endoscopy correspondences and other related topics as the extension of the
current local theory of theta correspondences and related topics. In this aspect, one may consider the theory of local descent as
a preliminary step of the local $(\tau,b)$-theory, which has been discussed in the work of
Soudry and the author (\cite{JS03}, \cite{JS04}, and \cite{JS12}), in the work of C.-F. Nien, Y.-J. Qin and the author (\cite{JNQ08},
\cite{JNQ10}, and \cite{JNQ12}), and the work of Baiying Liu (\cite{Liu11}) and of C. Jantzen and Liu (\cite{JaL12}) over the $p$-adic local fields.
Over the archimedean local fields, the recent work of A. Aizenbud, D. Gourevitch, and S. Sahi (\cite{AGS11},
of Gourevitch and Sahi in \cite{GS12}, and of B. Harris in \cite{H12} are relevant to the theme of this paper. However,
the local $(\tau,b)$-theory for archimedean local fields is open.

\subsection{Constructions of different types}
When the Weil representation is replaced by the minimal representation of exceptional groups,
the local theory of theta correspondences has been extended to
the so-called exceptional theta correspondences, which works both locally and globally. The local theory started from the work of
D. Kazhdan and G. Savin in \cite{KzS90} and the global theory is referred to the work of Ginzburg, Rallis and Soudry in \cite{GRS97}.
Among many contributors to the theory locally and globally are
Kazhdan, Savin, Rubenthaler, Ginzburg, Rallis, Soudry, J.-S. Li, Gross, Gan, Gurevich, the author and others. It should be
mentioned that important applications of the exceptional theta correspondences include the work of Ginzburg, Rallis and Soudry on
the theory of cubic Shimura correspondences (\cite{GRS97}), the work of Gan, Gurevich and the author on cuspidal spectrum of $G_2$ with
high multiplicity (\cite{GGJ02}), and the work of Gross and Savin on the existence of
motivic Galois group of type $G_2$ (\cite{GrS98}), to mention a few. The $(\tau,b)$-theory for exceptional groups is essential not known.

Finally, it is to mention the work of Ginzburg, Soudry and the author on the partial automorphic descent constructions for $\GL_2$
in \cite{GJS11}, which is expected to be in a different nature comparing to the whole theory discussed
here. A more general consideration of the constructions of this nature will be discussed in a future work.

\end{document}